\documentclass[11pt]{article}
\usepackage{epsfig,epsfig}
\usepackage{amsmath}
\usepackage{mathrsfs}
\usepackage{bm}
\usepackage{amssymb}
\usepackage{booktabs}
\usepackage{color}
\usepackage{multirow}
\usepackage{graphicx}
\usepackage{subfigure}
\usepackage{float}
\usepackage{appendix}
\usepackage{tikz}

\usepackage{srcltx,graphicx}

\newtheorem{thm}{Theorem}[section]

\newtheorem{proof}[thm]{Proof}

\newtheorem{example}[thm]{Example}

\usepackage{hyperref}

\numberwithin{equation}{section} \topmargin=-2cm \oddsidemargin=1cm
\begin{document}
\title{\textbf{An adaptive moving mesh discontinuous Galerkin method for the radiative transfer equation}\footnote{The research was partly supported by China NSAF grant U1630247, NSFC grants 11471049 and 11871111, and Science Challenge Project, No. TZ2016002.}}
\author{Min Zhang\footnote{School of Mathematical Sciences, Xiamen University,
Xiamen, Fujian 361005, China. E-mail: minzhang2015@stu.xmu.edu.cn
},
~Juan Cheng\footnote{Institute of Applied Physics and Computational Mathematics, Beijing 100094, China.
E-mail: cheng\_juan@iapcm.ac.cn
},
~Weizhang Huang\footnote{Department of Mathematics, University of Kansas, Lawrence, Kansas 66045, USA. E-mail: whuang@ku.edu
},
~and Jianxian Qiu\footnote{School of Mathematical Sciences and Fujian Provincial Key Laboratory of Mathematical Modeling and High-Performance Scientific Computing, Xiamen University, Xiamen, Fujian 361005, China.
E-mail: jxqiu@xmu.edu.cn
}
}

\date{}
\maketitle
\noindent\textbf{Abstract:}
The radiative transfer equation models the interaction of radiation with scattering and absorbing media and has important applications in various fields in science and engineering. It is an integro-differential equation involving time, space and angular variables and contains an integral term in angular directions while being hyperbolic in space. The challenges for its numerical solution include the needs to handle with its high dimensionality, the presence of the integral term, and the development of discontinuities and sharp layers in its solution along spatial directions. Its numerical solution is studied in this paper using an adaptive moving mesh discontinuous Galerkin method for spatial discretization together with the discrete ordinate method for angular discretization. The former employs a dynamic mesh adaptation strategy based on moving mesh partial differential equations to improve computational accuracy and efficiency. Its mesh adaptation ability, accuracy, and efficiency are demonstrated in a selection of one- and two-dimensional numerical examples.


\vspace{5pt}

\noindent\textbf{The 2010 Mathematics Subject Classification:} 65M50, 65M60, 65M70, 65R05, 65.75

\vspace{5pt}

\noindent\textbf{Keywords:}
adaptive moving mesh,
discontinuous Galerkin method,
unsteady radiative transfer equation,
high order accuracy,
high resolution

\newcommand{\h}{\hspace{1.cm}}
\newcommand{\hh}{\hspace{2.cm}}
\newtheorem{yl}{\hspace{1.cm}Lemma}
\newtheorem{dl}{\hspace{1.cm}Theorem}
\renewcommand{\sec}{\section*}
\renewcommand{\l}{\langle}
\renewcommand{\r}{\rangle}
\newcommand{\be}{\begin{eqnarray}}
\newcommand{\ee}{\end{eqnarray}}

\normalsize \vskip 0.2in
\newpage
\section{Introduction}
The radiative transfer equation (RTE) models the interaction of radiation with scattering and absorbing media,
which has important applications in fields such as astrophysics, high energy density physics,
nuclear physics, inertial confinement fusion, heat transfer, stellar atmospheres, optical molecular imaging,
infrared and visible light in space and the atmosphere, and biomedicine.
The RTE is an integro-differential equation with seven independent variables (for time, space, and angles)
for a time-dependent three-spatial-dimensional problem.
Containing an integral term and with its high dimensionality, the RTE
presents a challenge in the development of efficient numerical algorithms. On the other hand,
the efficient solution of the RTE plays an important role in the study of radiation hydrodynamics
where the RTE is often coupled with the Euler equations, the energy equation, and the equation of state.

In the past, a number of methods have been developed for the numerical solution of the RTE.
Those methods can be divided roughly into two categories: stochastic and deterministic approaches.
The Monte Carlo method is a widely used method in the former category \cite{MC-1, MC-2}.
On the other hand, deterministic approaches involve discrete approximations of the variables in the RTE.
In particular, the discretization needs to be applied to all coordinates in space and angles.
For angular coordinates, the $P_N$ method, first introduced in \cite{PN} and
also known as the spherical harmonics method, uses an orthogonal, harmonic basis to approximate the solution.
Another approach called the discrete ordinate method (DOM) \cite{BG-Carlson, Sn_Tn} employs
spectral collocation and the Legendre-Chebyshev quadrature to discretize the integro-differential
equation in angular coordinates. DOM is widely used for the numerical solution of the transport
equation \cite{SI, JC-Ragusa} due to its high accuracy, flexibility, and relatively low computational cost.
The angle-discretized RTE forms a system of linear hyperbolic equations with source terms,
which can be discretized in space using a standard method such as a finite difference, finite volume,
or finite element method. The discontinuous Galerkin (DG) method is employed for this purpose in the current work.

The DG method is known to be a particularly powerful numerical tool for the simulation of hyperbolic
transport problems. It was first used for the RTE by Reed and Hill \cite{reed} and theoretically studied
by Lesaint and Raviart \cite{LeSaint}. The method was later extended to nonlinear conservation laws
by Cockburn and Shu \cite{bs1,bs2,bs3,bs4}. The DG method has the advantages of high-order accuracy,
geometric flexibility, suitability for handling $h$- and $p$-adaptivity, extremely local data structure,
high parallel efficiency, and a good theoretical foundation for stability and error estimates. Over the last
few decades, the DG method has been used widely in scientific and engineering computation.

The objective of the current work is to study an adaptive moving mesh DG method (for spatial discretization)
combined with DOM (for angular discretization) for the numerical solution of the unsteady RTE.
Due to its hyperbolic nature, the solution of the RTE can develop discontinuities or sharp layers along
spatial directions, which makes mesh adaptation
an indispensable tool for use in improving computational accuracy and efficiency.
Mesh adaptation methods can be classified roughly into three groups.
The first one is $h$-methods, which generate a new mesh by adding or removing points to an existing mesh.
Typically, mesh points are added in regions where the solution variation or error is large,
and mesh points are removed in regions where the solution is smooth.
The second group is $p$-methods where the order of polynomial approximation varies from
place to place according to a certain error estimate or indicator.
The third group is $r$-methods, also called moving mesh methods, which relocate
mesh point positions while maintaining the total number of mesh points and the mesh connectivity.

Several works have been done in using mesh adaptation for the numerical solution of the RTE.
In \cite{RTE2}, an adaptive mesh refinement (AMR)
algorithm has been formulated and implemented for the RTE by minimizing the spatial discretization error.
In \cite{RTE1}, $hp$-adaptive DG methods have been developed for the numerical solution of
simplified $P_N$ approximations of radiative transfer in non-grey semitransparent media and
it has been found that it is possible to approximate the radiative field with a significantly lower
computational cost than solving the equations using the conventional finite element method.
In \cite{RTE3}, $hp$-adaptive methods have been develop based on a short-characteristics approach
embedded in the discontinuous finite element framework.
No work has been done so far in using moving mesh methods for the numerical solution of the RTE.

We here study a moving mesh method based on moving mesh partial differential equations
(MMPDEs) \cite{Huang2011, Huang1994-1, Huang1994-2} for this purpose.
An MMPDE moves the mesh continuously in time and
orderly in space and is formulated as the gradient flow equation of a meshing functional.
We use a newly developed discrete approach \cite{Huang2015} that
makes the implementation of the MMPDE method not only significantly simpler
but also much more reliable in the sense that there is a theoretical guarantee for mesh nonsingularity
at semi and fully discrete levels \cite{Huang2018}.
The MMPDE method determines the mesh adaptivity for the size, shape,
and orientation of mesh elements using a metric tenor (also called a monitor function) which is matrix-valued
function defined throughout the physical domain. 

The full discretization of the RTE includes the discretization in angular directions using DOM,
in space using an adaptive moving mesh DG method, and in time using the backward Euler scheme.
The DOM discretization leads to radiative intensity functions for different angular directions
which can have discontinuities or sharp layers at different location in space. To take this into account,
we compute a metric tensor based on the Hessian for each of these functions and then combine
them via matrix intersection (see its definition in Section~\ref{sec:mesh})
into a single metric tensor to be used with the MMPDE method.
Numerical results in Section~\ref{sec:numerical} show that the adaptive moving mesh DG method
with this strategy works well for problems with single or multiple sharp layers and discontinuities in the sense that
it is able to automatically concentrate the mesh points in the regions of discontinuities or steep transition
layers and is more efficient than its fixed mesh counterpart.


An outline of the paper is as follows.
In Section \ref{sec:RTE}, the unsteady RTE and DOM are described.
The adaptive moving mesh DG method for solving the two-dimensional unsteady RTE is presented
in Section \ref{sec:2d}. The generation of adaptive moving meshes using a new implementation
of the MMPDE method is discussed in Section \ref{sec:mesh}. A selection of one- and two-dimensional examples
are presented in Section \ref{sec:numerical} to demonstrate the mesh adaptation ability
of the adaptive moving mesh DG method and its accuracy and efficiency.
Finally, Section \ref{sec:conclusion} contains conclusions and further comments.

\section{The unsteady radiative transfer equation}
\label{sec:RTE}

The RTE is an integro-differential equation modeling the conservation of photons \cite{gc}.
We consider a case with one-group, isotropically scattering radiative transfer. The governing equation
for this case reads as
\begin{equation}\label{s2.1}
\begin{split}
\frac{1}{c}\frac{\partial I(\bm{r},\Omega,t)}{\partial t}
+\Omega \cdot \nabla I(\bm{r},\Omega,t)+\sigma_tI(\bm{r},\Omega,t)
=\frac{\sigma_s}{4\pi}\int_S I(\bm{r},\tilde{\Omega},t)d\tilde{\Omega}+q(\bm{r},\Omega,t),
\end{split}
\end{equation}
where $c$ is the speed of photons, $\bm{r}$ is the spatial variable,
$\nabla$ is the gradient operator with respect to $\bm{r}$,
$\Omega$ is the unit angular variable,
$S$ is the unit sphere, $t$ is time,
$I(\bm{r},\Omega,t)$ is the radiative intensity in the direction $\Omega$,
$\sigma_s \geq 0 $ is the scattering coefficient of the medium,
$ \sigma_t \geq \sigma_s $ is the extinction coefficient of the medium due to both absorption and scattering,
and $ q(\bm{r}, \Omega,t)$ is a given source term.
The vector $\bm{r}$ is described by the Cartesian coordinates $x, y, z$
while $\Omega$ is usually described by a polar angle $\beta$ measured with respect to a fixed axis
in space (such as the $z$ axis) and a corresponding azimuthal angle $\varphi$.
Letting $\mu = \cos\beta$,
then
\begin{equation*}
 d\bm{r} = dxdydz, \quad d\Omega = \sin\beta d\beta d\varphi = -d\mu d\varphi.
\end{equation*}

In this work we consider the numerical solution of \eqref{s2.1} in one and two spatial dimensions.
Since the numerical algorithm is similar in one and two dimensions, we describe
it only in two dimensions. The equation \eqref{s2.1} reads in two dimensions as
\begin{equation}\label{s3.1}
\begin{split}
&\frac{1}{c}\frac{\partial I(x,y,\Omega,t)}{\partial t}
+\Omega \cdot \nabla I(x,y,\Omega,t)+\sigma_tI(x,y,\Omega,t)
\\
&\qquad \qquad = \frac{\sigma_s}{4\pi}\int_S I(x,y,\tilde{\Omega},t)d\tilde{\Omega} + q(x,y,\Omega,t),
\quad (x,y)\in \mathbb{D},
\end{split}
\end{equation}
where $\Omega = (\zeta,\eta)$ and
\[
\zeta = \sin\beta \cos\varphi=\sqrt{1-\mu^2} \cos\varphi \in [-1, 1],
\quad \eta =\sin\beta \sin\varphi=\sqrt{1-\mu^2} \sin\varphi \in [-1, 1].
\]
Denote by $\bm{n}(x,y)$ the unit outward normal vector of the domain boundary $\partial\mathbb{D}$
at the point $(x,y)$
and define $\partial \mathbb{D}_{in} = \{(x,y)\in \partial\mathbb{D} \; |\;  \bm{n}(x,y) \cdot\Omega<0\}$.
Then, the boundary condition can be expressed as
\begin{equation}\label{s3.2b}
I(x,y,\Omega,t) = g(x,y,\Omega,t) ,\quad (x,y)\in \partial \mathbb{D}_{in},
\end{equation}
and the initial condition is
\begin{equation}\label{s3.2i}
I(x,y,\Omega,0) = I_0(x,y,\Omega) , \quad (x,y)\in \mathbb{D},
\end{equation}
where $g(x,y,\Omega,t)$ and $I_0(x,y,\Omega)$ are given functions.
Note that no boundary condition is needed in $\Omega$ directions.

The RTE \eqref{s3.1} needs to be discretized in angular, spatial, and time variables.
A challenge for solving \eqref{s3.1} is due to its high dimensionality: it has two dimensions
totally in both angular and spatial coordinates for a one-spatial-dimensional problem and four dimensions
for a two-spatial-dimensional problem. To tackle this challenge, a common strategy is to use
a high accuracy discretization in angular coordinates. We use the discrete-ordinate method \cite{BG-Carlson}
for this purpose. DOM will be described later in this section.
Another challenge for solving \eqref{s3.1} is due to the fact that it is hyperbolic in space.
This means that its solution can develop discontinuities and sharp layers across the physical domain,
which requires a spatial discretization that can handle those structures and mesh adaptation that
can provide high resolution in regions around them.  An adaptive mesh DG method to be presented
in the next section will be used to tackle this challenge.

DOM \cite{BG-Carlson} is a spectral collocation-type method \cite{Canuto1988}. Indeed,
\eqref{s3.1} is collocated for a finite number of angular directions while the integral
in the angular variable is approximated by the Legendre-Chebyshev quadrature \cite{Sn_Tn}
where the nodes in $\mu$ and $\varphi $ are chosen as the roots of Legendre and
Chebyshev polynomials, respectively. Specifically, the discrete-ordinate approximation
of \eqref{s3.1} is given by
\begin{equation}\label{s3.3}
\begin{split}
&\frac{1}{c}\frac{\partial I_m(x,y,t)}{\partial t}
+\Omega_m \cdot \nabla I_m(x,y,t)+\sigma_tI_m(x,y,t)
\\
& \qquad \qquad = \sigma_s\sum_{m'=1}^{N_a}w_{m'}I_{m'}(x,y,t)+q_m(x,y,t),
\quad m=1,\cdots, N_a,
\end{split}
\end{equation}
where $\Omega_m = (\zeta_m,\eta_m),\; m=1,\cdots,N_a$, are the discrete angular directions,
$I_m(x,y,t)$ is an approximation of $I(x,y,\Omega_m,t)$,
$q_m(x,y,t) = q(x,y,\Omega_m,t)$,
and $\sum_{m=1}^{N_a}w_{m}I_{m}(x,y,t)$
is a Legendre-Chebyshev quadrature rule with weights $w_m>0$
for $({1}/{4\pi})\int_S I(x,y,\tilde{\Omega},t)d\tilde{\Omega}$.

\section{An adaptive moving mesh DG method for the two-dimensional unsteady DOM RTEs}\label{sec:2d}

We notice that \eqref{s3.3} is a system of hyperbolic equations and its solution can develop discontinuities
and sharp layers. The DG method has been known to be a powerful numerical tool for the simulation of hyperbolic
problems with discontinuous solutions \cite{LeSaint,reed}.
Mesh adaptation is also crucial to provide accurate resolution of discontinuities
and sharp layers in the solution and improve computational efficiency.
We describe a DG method for \eqref{s3.3} on a general adaptive moving mesh
in this section and the adaptive mesh movement in the next section.

Specifically, we consider time instants
\[
t_0 = 0 < t_1 < \cdots < t_{n} < t_{n+1} < \cdots .
\]
For the moment, we assume that a triangular mesh, which consists of non-overlapping triangles
covering $\mathbb{D}$ completely and whose vertices depend on $t$, is known at the time instants,
i.e., $\mathscr{T}_{h}^n$, $n = 0, 1, ...$, are given. We also assume that the mesh
keeps the same connectivity and the same number of elements and vertices for the whole time period.
(The position of the vertices is the only thing that changes with time.)
The generation of such a moving mesh is discussed
in Section~\ref{sec:mesh}.

For $t\in[t_n,t_{n+1}]$, the coordinates and velocities of the vertices of the mesh are defined as
\begin{equation}\label{xy}
\begin{split}
&x_{j}(t)= x^{n}_{j} \frac{t_{n+1}-t}{\Delta t_n}+x^{n+1}_{j}\frac{t-t_{n}}{\Delta t_n},
\quad \dot{x}_{j}(t)=\frac{x^{n+1}_{j}- x^{n}_{j}}{\Delta t_n}, \quad j=1,\cdots,N_v,
\\
&y_{j}(t)= y^{n}_{j} \frac{t_{n+1}-t}{\Delta t_n}+y^{n+1}_{j}\frac{t-t_{n}}{\Delta t_n},
\quad\dot{y}_{j}(t)=\frac{y^{n+1}_{j}- y^{n}_{j}}{\Delta t_n},\quad j=1,\cdots,N_v,
\end{split}
\end{equation}
where $N_v$ is the number of the vertices.
The corresponding mesh is denoted by $\mathscr{T}_{h}(t)$.

We now describe the DG discretization of \eqref{s3.3} on $\mathscr{T}_{h}(t)$.
For any element $K\in \mathscr{T}_h(t)$, denote its vertices by
$(x^K_1,y^K_1)$, $(x^K_2,y^K_2)$, $(x^K_3,y^K_3)$ and its area by $|K|$.
Consider a set of local orthogonal polynomials of up to degree $k$  in $K$,
\begin{equation}
P^k(K,t)=\text{span}\{\phi^{(K)}_0(x,y,t),\; \phi^{(K)}_1(x,y,t),\; \cdots,\; \phi^{(K)}_{L-1}(x,y,t)\},
\end{equation}
where $L = {(k+1)(k+2)}/{2}$ is the dimension of $P^k(K,t)$.
Then the associated DG finite element space can be defined as
\begin{equation}
V_h^k(t)=\{I_m^h(x,y,t)\in L^2(\mathbb{D})\; :\; I_m^h(x,y,t)|_{K}\in P^k(K,t),\quad \forall K \in \mathscr{T}_h(t) \}.
\end{equation}
Notice that any function in this space can be expressed as
\begin{equation}\label{2dm.1}
I_m^h(x,y,t) =\sum_{p=0}^{L-1} I_{m,K}^{[p]}(t)\phi^{(K)}_p(x,y,t), \quad (x,y)\in K,
\end{equation}
where $I_{m,K}^{[p]}$'s are the degrees of freedom. Moreover, its time derivative can be written as
\begin{align}\label{2dm.2}
\frac{\partial I_m^h(x,y,t)}{\partial t} & =\sum_{p=0}^{L-1} \Big{(}\frac{dI_{m,K}^{[p]}(t)}{dt}\phi^{(K)}_p(x,y,t)+I_{m,K}^{[p]}(t)\frac{\partial \phi^{(K)}_p(x,y,t)}{\partial t}\Big{)} .
\end{align}
It is not difficult (e.g., see \cite{deforming}) to show  that
\begin{equation}\label{2dm.2+1}
\frac{\partial \phi^{(K)}_p(x,y,t)}{\partial t}=-\Pi_1(x,y,t) \cdot \nabla \phi^{(K)}_p(x,y,t),
\end{equation}
where $\Pi_1(x,y,t)=(\dot{X},\dot{Y})$ is the piecewise linear interpolation of the nodal mesh velocities, i.e.,
\begin{equation}\label{2dm.4}
\begin{split}
&\dot{X}=
 \frac{1}{3}(\dot{x}^K_1+\dot{x}^K_2+\dot{x}^K_3)\phi^{(K)}_0
-\frac{\sqrt{2}}{4}(\dot{x}^K_1-2\dot{x}^K_2+\dot{x}^K_3)\phi^{(K)}_1
-\frac{\sqrt{2}}{2}(\dot{x}^K_1-\dot{x}^K_3)\phi^{(K)}_2,
\\&\dot{Y}=
 \frac{1}{3}(\dot{y}^K_1+\dot{y}^K_2+\dot{y}^K_3)\phi^{(K)}_0
-\frac{\sqrt{2}}{4}(\dot{y}^K_1-2\dot{y}^K_2+\dot{y}^K_3)\phi^{(K)}_1
-\frac{\sqrt{2}}{2}(\dot{y}^K_1-\dot{y}^K_3)\phi^{(K)}_2 .
\end{split}
\end{equation}
Combining (\ref{2dm.2}) and (\ref{2dm.2+1}), we get
\begin{align}
\label{2dm.5}
\frac{\partial I_m^h(x,y,t)}{\partial t} & = \sum_{p=0}^{L-1}\frac{dI_{m,K}^{[p]}(t)}{dt}\phi^{(K)}_p(x,y,t)
-\Pi_1(x,y,t)\cdot \nabla I_{m}^h(x,y,t), \quad (x,y)\in K .
\end{align}
Multiplying \eqref{s3.3} by a test function $\phi(x,y,t)\in V_h^k(t)$, integrating
the resulting equation over $K$, replacing $I_m(x,y,t)$ with its approximation $I_m^h(x,y,t)$,
and using \eqref{2dm.5}, we have
\begin{equation}\label{2dm.7}
\begin{split}
&\int_{K}\Big{(}\frac{1}{c}\sum_{p=0}^{L-1}\frac{dI_{m,K}^{[p]}(t)}{dt}\phi^{(K)}_p\Big{)}\phi dxdy
+\int_{K}\sigma_tI_{m}^h(x,y,t)\phi dxdy
\\&+\int_{K}(\Omega_m -\frac{1}{c}\Pi_1(x,y,t))\cdot \nabla I_{m}^h(x,y,t)\; \phi dxdy
\\=&\int_{K} \sigma_s\Psi_K(x,y,t)\phi dxdy+\int_{K}q_m(x,y,t)\phi dxdy,
\end{split}
\end{equation}
where
\[
\Psi_K(x,y,t)=\sum_{m'=1}^{N_a}w_{m'}I^h_{m', K}(x,y,t)
\]
and $I^h_{m', K}(x,y,t)$ denotes the restriction of $I^h_{m'}(x,y,t)$ on $K$.
Applying the divergence theorem on the third term, we get
\begin{equation}\label{2dm.8}
\begin{split}
&\int_{K}\Big{(}\frac{1}{c}\sum_{p=0}^{L-1}\frac{dI_{m,K}^{[p]}(t)}{dt}\phi^{(K)}_p\Big{)}\phi dxdy
+\int_{K}\sigma_tI_{m}^h(x,y,t)\phi dxdy
\\&-\int_{K}\left ( \nabla \cdot \Big( \phi (\Omega_m -\frac{1}{c}\Pi_1(x,y,t))\Big) \right ) I_{m}^h(x,y,t)) dxdy
\\&+ \int_{\partial K} \bm{n}_K\cdot \Big(\Omega_m-\frac{1}{c} \Pi_1(x,y,t) \Big) I_{m}^h(x,y,t)\phi dxdy
\\=&\int_{K} \sigma_s\Psi_K(x,y,t)\phi dxdy+\int_{K}q_m(x,y,t)\phi dxdy ,
\end{split}
\end{equation}
where $\bm{n}_{K}$ is the outward unit normal to the boundary $\partial K$.
In the above equation, $I_{m}^h(x,y,t)$ is discontinuous across the cell boundaries
in general and its value thereon is not well defined.
To specify the value, we define the outflow boundary $\partial K^{m+}$ and
the inflow boundary $\partial K^{m-}$ of the cell $K$ by
\begin{equation}
\label{boundary-red}
\begin{split}
&\partial K^{m+}=\{(x,y)\in \partial K\; |\; \Big(\Omega_m -\frac{1}{c} \Pi_1(x,y,t) \Big) \cdot \bm{n}_{K}(x,y) \geq 0\},\\
&\partial K^{m-}=\{(x,y)\in \partial K\; |\; \Big(\Omega_m - \frac{1}{c} \Pi_1(x,y,t) \Big) \cdot \bm{n}_{K}(x,y) < 0\}.
\end{split}
\end{equation}
It is useful to point out that, in practice, $(1/c) \Pi_1(x,y,t)$ is much smaller than $\Omega_m$ for most
situations and $\partial K^{m+}$ and $\partial K^{m-}$ can be computed using the simpler formulas
\begin{equation}
\label{boundary}
\begin{split}
&\partial K^{m+}=\{(x,y)\in \partial K\; |\; \Omega_m\cdot \bm{n}_{K}(x,y) \geq 0\},\\
&\partial K^{m-}=\{(x,y)\in \partial K\; |\; \Omega_m\cdot \bm{n}_{K}(x,y) < 0\} .
\end{split}
\end{equation}
Since each interior edge is shared by two triangular elements, the value
of $I_{m}^h(x,y,t)$ on any edge of $K$ can be defined based on its value in $K$
or in the other element sharing the common edge with $K$.
These values are denoted by $I_{m}^h(int(K),t)$) and $I_{m}^h(ext(K),t)$, respectively.
For the upwind numerical flux, we use $I_{m}^h(int(K),t)$ for the outflow boundary and
$I_{m}^h(ext(K),t)$ for the inflow boundary. Thus, we can rewrite \eqref{2dm.8} into
\begin{equation}\label{2dm.9}
\begin{split}
&\int_{K}\Big{(}\frac{1}{c}\sum_{p=0}^{L-1}\frac{dI_{m,K}^{[p]}(t)}{dt}\phi^{(K)}_p\Big{)}\phi dxdy
+\int_{K}\sigma_tI_{m}^h(x,y,t)\phi dxdy
\\&-\int_{K}\left ( \nabla \cdot \Big( \phi (\Omega_m -\frac{1}{c}\Pi_1(x,y,t))\Big) \right ) I_{m}^h(x,y,t) dxdy
\\&+ \int_{\partial K^{m-}} \bm{n}_K\cdot \Big(\Omega_m-\frac{1}{c} \Pi_1(x,y,t) \Big) I_{m}^h(ext(K),t)\phi dxdy
\\&+ \int_{\partial K^{m+}} \bm{n}_K\cdot \Big(\Omega_m-\frac{1}{c} \Pi_1(x,y,t) \Big) I_{m}^h(int(K),t)\phi dxdy
\\=&\int_{K} \sigma_s\Psi_K(x,y,t)\phi dxdy+\int_{K}q_m(x,y,t)\phi dxdy .
\end{split}
\end{equation}

Explicit time stepping can cause extremely small time steps due to the high photon speed.
To avoid this difficulty, we use the backward Euler scheme for \eqref{2dm.9}, i.e.,
\begin{equation}\label{2dm.10}
\begin{split}
&\int_{K}\Big{(}\frac{1}{c}\sum_{p=0}^{L-1}\frac{I_{m,K}^{[p]}(t_{n+1})-I_{m,K}^{[p]}(t_n)}
{\Delta t_n}\phi^{(K)}_p\Big{)}\phi dxdy
+\int_{K}\sigma_tI_{m}^h(x,y,t_{n+1})\phi dxdy
\\&-\int_{K}\left ( \nabla \cdot \Big( \phi (\Omega_m -\frac{1}{c}\Pi_1(x,y,t_{n+1}))\Big) \right ) I_{m}^h(x,y,t_{n+1}) dxdy
\\&+ \int_{\partial K^{m-}} \bm{n}_K\cdot \Big(\Omega_m-\frac{1}{c} \Pi_1(x,y,t_{n+1}) \Big) I_{m}^h(ext(K),t_{n+1})\phi dxdy
\\&+ \int_{\partial K^{m+}} \bm{n}_K\cdot \Big(\Omega_m-\frac{1}{c} \Pi_1(x,y,t_{n+1}) \Big) I_{m}^h(int(K),t_{n+1})\phi dxdy
\\=&\int_{K} \sigma_s\Psi_K(x,y,t_{n+1})\phi dxdy+\int_{K}q_m(x,y,t_{n+1})\phi dxdy ,
\\& \qquad \qquad \qquad \quad \qquad \qquad \qquad \quad \forall \phi \in V_h^k(t_{n+1}), \quad K \in \mathscr{T}_h(t_{n+1}), \quad m = 1, ..., N_a.
\end{split}
\end{equation}
The above equations form a coupled system for the unknown functions
$I_{m}^h(x,y,t_{n+1})$, $m = 1, ..., N_a$ since the function $\Psi_K(x,y,t_{n+1})$ contains all of them.
To decouple these functions from the equations, a functional-type iteration called
the source iteration (SI) \cite{SI,CJ} (also referred to as the grid sweeping algorithm)
has been widely used for solving the system in a Gauss-Seidel-like manner.
To be specific, assuming that the $\ell$-th iteration solutions $I_{m,K}^{h(\ell)}(x,y,t_{n+1})$ (for $m=1,\cdots,N_a$ and $K \in \mathscr{T}_h(t_{n+1})$) are known, we compute the new approximations $I^{h(\ell+1)}_{m, K}(x,y,t_{n+1})$
element by element in a sweeping direction \cite{Sweep} and through all angular directions $m=1,\cdots,N_a$ for each given element. Thus, for $K \in \mathscr{T}_h(t_{n+1})$, we have
\begin{equation}\label{2dm.11}
\begin{split}
&\int_{K}\Big{(}\frac{1}{c}\sum_{p=0}^{L-1}\frac{I_{m,K}^{[p] (\ell+1)}(t_{n+1})-I_{m,K}^{[p]}(t_n)}
{\Delta t_n}\phi^{(K)}_p\Big{)}\phi dxdy
+\int_{K}\sigma_tI_{m,K}^{h(\ell+1)}(x,y,t_{n+1})\phi dxdy
\\&-\int_{K}\left ( \nabla \cdot \Big( \phi (\Omega_m -\frac{1}{c}\Pi_1(x,y,t_{n+1}))\Big) \right )
I_{m,K}^{h(\ell+1)}(x,y,t_{n+1}) dxdy
\\&+ \int_{\partial K^{m-}} \bm{n}_K\cdot \Big(\Omega_m-\frac{1}{c} \Pi_1(x,y,t_{n+1}) \Big)
I_{m,K}^{h(\ell+1)}(ext(K),t_{n+1})\phi dxdy
\\&+ \int_{\partial K^{m+}} \bm{n}_K\cdot \Big(\Omega_m-\frac{1}{c} \Pi_1(x,y,t_{n+1}) \Big)
I_{m,K}^{h(\ell+1)}(int(K),t_{n+1})\phi dxdy
\\=&\int_{K} \sigma_s\Psi_K^{*}(x,y,t_{n+1})\phi dxdy+\int_{K}q_m(x,y,t_{n+1})\phi dxdy ,
\\& \qquad \qquad \qquad \quad \qquad \qquad \qquad \quad \forall \phi \in V_h^k(t_{n+1}), \quad m = 1, ..., N_a,
\end{split}
\end{equation}
where
\[
\Psi^{*}_K(x,y,t_{n+1})=\sum_{m'=1}^{N_a}w_{m'}I^{*}_{m',K}(x,y,t_{n+1}),
\]
\begin{equation*}\label{2dm.12}
I^{*}_{m',K}(x,y,t_{n+1})=\left\{\begin{array}{ll}
I_{m',K}^{h(\ell+1)}(x,y,t_{n+1}),\quad &\text{when avialiable}, \\
I_{m',K}^{h(\ell)}(x,y,t_{n+1}), \quad &\text{otherwise}.
\end{array}\right.
\end{equation*}
The iteration is stopped when the difference between two consecutive iterates is smaller than a given
tolerance. In our computation, we use $\max_{m} \|I_{m}^{h(\ell+1)}-I_{m}^{h(\ell)}\|_{\infty} \leq 10^{-12}$.
The source iteration is very effective, taking only a few iterations to achieve convergence for most of the problems
tested.

The time integration alternates between solving the physical equation and generating the mesh.
Starting with the current mesh $\mathscr{T}^{n}_h$ and a solution $I^h_m(x,y,t_n),\; m=1,\cdots, N_a$,
a new mesh $\mathscr{T}^{n+1}_h$ is generated using the MMPDE moving mesh strategy to be described
in the next section. Then, the DOM-DG scheme \eqref{2dm.11} at $t_{n+1}$ is solved
for the new solution approximation $I^h_m(x,y,t_{n+1}),\; m=1,\cdots, N_a$.

\section{The MMPDE moving mesh method on triangular meshes}
\label{sec:mesh}

In this section we describe the generation of $\mathscr{T}^{n+1}_{h}$ based on $\mathscr{T}^{n}_{h}$ and numerical solution $I^h_{m}(x,y,t_n),\; m=1,\cdots,N_a$ using the MMPDE moving mesh method
\cite{Huang2011,Huang1994-1,Huang1994-2}. The method utilizes a metric tensor
(or called a monitor function) to provide the information of the size, shape, orientation of mesh elements
throughout the domain that is needed for mesh adaptation.
We use here a new implementation of the method proposed in \cite{Huang2015}.
A unique feature for the numerical solution of the RTE is that the functions $I^h_{m}(x,y,t_n),\; m=1,\cdots,N_a$
which correspond to the radiative intensity at angular directions $\Omega_m$, $m=1,\cdots,N_a$
may have discontinuities and sharp layers at different locations in space. 
To take this into account, we first compute the metric tensor for each function and
then combine all of the metric tensors into a single one. 

We start with noting that $\mathscr{T}^{n+1}_{h}$ and $\mathscr{T}^{n}_{h}$ have the same number
of the elements $(N)$, the same number of the vertices $(N_v)$, and the same connectivity.
They differ only in the location of the vertices. We assume that a reference computational mesh
$\hat{\mathscr{T}}_c = \{(\hat \xi_j, \hat \eta_j),\;  j=1,\cdots N_v \}$, which also has the same connectivity
and the same numbers of vertices and elements as $\mathscr{T}^n_{h}$, has been chosen.
In our computation, we take it as a uniform mesh (in the Euclidean metric) defined on domain $\mathbb{D}$.
$\hat{\mathscr{T}}_c$  stays fixed in the computation.  The generation of $\mathscr{T}^{n+1}_{h}$ is
through a computational mesh $\mathscr{T}_c= \{(\xi_j,\eta_j), j=1,\cdots N_v \}$
which serves as an intermediate variable.

A key idea of the MMPDE moving mesh method is to view any nonuniform mesh as a uniform one
in some metric $\mathbb{M}$ \cite{Huang2006, Huang2011}.
The metric tensor $\mathbb{M} = \mathbb{M}(x,y,t)$ is a symmetric and uniformly positive definite
matrix-valued function defined on $\mathbb{D}$. It provides the magnitude and direction information
needed for determining the size, shape, and orientation of the mesh elements throughout the domain.
Various metric tensors have been proposed; e.g., see \cite{Huang2011, Huang2003}.
We here use a metric tensor based on the Hessian of the computed solution.
Let $H_K(I^h_{m}(t_n))$  be the Hessian or a recovered Hessian of $I^h_{m}(t_n)$ on $K$.
Let the eigen-decomposition of $H_K(I^h_{m}(t_n))$ be
\[
H_K(I^h_{m}(t_n)) = Q\hbox{diag}(\lambda_1,\lambda_2)Q^T .
\]
Denote
\[
|H_K(I^h_{m}(t_n))| = Q\hbox{diag}(|\lambda_1|,|\lambda_2|)Q^T.
\]
The metric tensor is then defined as
\begin{equation}\label{mer}
\mathbb{M}_{K,m} =\det \big{(}\mathbb{I}+\frac{1}{\alpha_{h,m}}|H_K(I^h_{m}(t_n))|\big{)}^{-\frac{1}{6}}
\big{(}\mathbb{I}+\frac{1}{\alpha_{h,m}}|H_K(I^h_{m}(t_n))|\big{)},
\quad \forall K \in \mathscr{T}_h
\end{equation}
where $\mathscr{T}_h$ denotes a general physical mesh,
$\mathbb{I}$ is the identity matrix, $\det(\cdot)$ is the determinant of a matrix,
and $\alpha_{h,m}$ is a regularization parameter defined through the algebraic equation
\[
\sum_{K\in\mathscr{T}_h}|K|\, \hbox{det}(\mathbb{M}_{K,m})^{\frac{1}{2}}=2\sum_{K\in\mathscr{T}_h}|K|\,
\hbox{det}(|H_K(I^h_{m}(t_n))|)^{\frac{1}{3}}.
\]
The metric tensor (\ref{mer}) is known \cite{Huang2003} to be optimal for the $L^2$-norm of
linear interpolation error. 

Notice that for each $m$ ($ 1 \le m \le N_a$), $\mathbb{M}_{K,m}$ provides the mesh adaptation
information only for function $I^h_{m}(x,y,t_n)$. To account for all of the functions in the mesh adaptation,
we need to combine the metric tensors into a single one.  We define
\[
\mathbb{M}_{K} = \mathbb{M}_{K,1}\cap \mathbb{M}_{K,2}\cap \cdots \cap \mathbb{M}_{K,N_a},
\]
where ``$\cap$'' stands for the intersection of symmetric and positive definite matrices which is defined as
follows. Let $A$ and $B$ be two symmetric and positive definite matrices. There exists a nonsingular matrix
$P$ such that $P A P^T = \mathbb{I}$ and $P B P^T = \text{diag}(b_1, b_2)$. The intersection of $A$ and $B$
is then defined as $A\cap B = P^{-1} \text{diag}(\max(1, b_1), \max(1, b_2)) P^{-T}$.
Define
\[
\mathbb{D}_C = \{ (x,y):\; (x, y) C (x,y)^T < 1\}
\]
for any symmetric and positive definite matrix $C$.
It is not difficult to show that $\mathbb{D}_{A\cap B} \subseteq \mathbb{D}_{A} \cap \mathbb{D}_{B}$,
which gives the meaning of the name ``intersection''.
Notice that the definition is not optimal in the sense that $\mathbb{D}_{A\cap B}$
is not necessarily the biggest ellipse inscribed in $\mathbb{D}_{A}$ and $\mathbb{D}_{B}$.

It is known \cite{Huang2006, Huang2011} that if $\mathscr{T}_h$ is uniform
in the metric $\mathbb{M}$ in reference to the computational  mesh $\mathscr{T}_c$, it satisfies
\begin{equation}\label{ec}
\begin{split}
|K|\sqrt{\det(\mathbb{M}_K)}=\frac{\sigma_h|K_c|}{|\mathbb{D}_c|},\quad \forall K\in \mathscr{T}_h,
\end{split}
\end{equation}
\begin{equation}\label{al}
\begin{split}
\frac{1}{2}\hbox{tr}\big{(} (F'_K)^{-1}\mathbb{M}_K^{-1}(F'_K)^{-T}\big{)}
= \hbox{det}\big{(}
(F'_K)^{-1}\mathbb{M}_K^{-1}(F'_K)^{-T}\big{)}^{\frac{1}{2}},\quad \forall K \in \mathscr{T}_h,
\end{split}
\end{equation}
where $F'_K$ is the Jacobian matrix of the affine mapping: $F_K: K_c\in \mathscr{T}_c \rightarrow K$,
$\mathbb{M}_K$ is the average of $\mathbb{M}$ over $K$, $\hbox{tr}(\cdot)$ denotes the trace of a matrix, and
\[
|\mathbb{D}_c|=\sum\limits_{K_c\in\mathscr{T}_c}|K_c|,\quad \sigma_h=\sum\limits_{K\in\mathscr{T}_h}|K|\hbox{det}(\mathbb{M}_K)^{\frac{1}{2}} .
\]
The condition \eqref{ec}, called the equidistribution condition, determines the size of elements through
the metric tensor $\mathbb{M}$. On the other hand,  \eqref{al}, referred to as the alignment condition,
determines the shape and orientation of elements through $\mathbb{M}_K$ and shape of $K_c$.
An energy function associated with these conditions is given by
\begin{equation}\label{fl}
\begin{split}
\mathcal{I}_h(\mathscr{T}_h,\mathscr{T}_c) =&\frac{1}{3}\sum_{K\in\mathscr{T}_h}|K|\hbox{det}(\mathbb{M}_K)^{\frac{1}{2}}\big{(}
\hbox{tr}((F'_K)^{-1}\mathbb{M}^{-1}_K(F'_K)^{-T})\big{)}^{2}
\\&+\frac{4}{3}\sum_{K\in\mathscr{T}_h}|K|\hbox{det}(\mathbb{M}_K)^{\frac{1}{2}}
\left (\hbox{det}(F'_K) \hbox{det}(\mathbb{M}_K)^{\frac{1}{2}} \right )^{-2},
\end{split}
\end{equation}
which is actually a Riemann sum of a continuous functional developed in \cite{Huang2001}
based on mesh equidistribution and alignment.

Note that $\mathcal{I}_h(\mathscr{T}_h,\mathscr{T}_c)$ is a function of the vertices $\bm{\xi}_j = (\xi_j, \eta_j),
\; j=1,\cdots,N_v$, of the computational mesh $\mathscr{T}_c$ and the vertices $\bm{x}_j = (x_j, y_j),\;
j=1,\cdots,N_v$, of the physical mesh $\mathscr{T}_h$. A straight way of solving the minimization problem
is to take $\mathscr{T}_c$ as $\hat{\mathscr{T}}_c$ and then solve the minimization problem of
$\mathcal{I}_h(\mathscr{T}_h,\hat{\mathscr{T}}_c)$ for the new physical mesh $\mathscr{T}_h^{n+1}$.
However, $\mathcal{I}_h(\mathscr{T}_h,\hat{\mathscr{T}}_c)$ is highly nonlinear in
$\bm{\xi}_j = (\xi_j, \eta_j), \; j=1,\cdots,N_v$. The fact that $\mathbb{M}$ is a function of $\bm{x}$
and thus $\mathbb{M}_K$ is a function of the coordinates of the physical vertices makes the situation more difficult.
Here, we adopt an indirect approach, i.e., to take $\mathscr{T}_h$
as $\mathscr{T}_h^n$, minimize $\mathcal{I}_h(\mathscr{T}_h^n,\mathscr{T}_c)$ with respect to $\mathscr{T}_c$,
and then obtain the new physical mesh using the relation between $\mathscr{T}_h^n$ and newly obtained
$\mathscr{T}_c$. The minimization is carried out by integrating the mesh equation which is defined as the gradient
system of the energy function (the MMPDE approach), viz.,
\begin{equation}\label{MM}
\begin{split}
\frac{d \bm{\xi}_j }{dt}=-\frac{\hbox{det}(\mathbb{M}(\bm{x_j}))^{\frac{1}{2}} }{\tau}\Big{(}\frac{\partial \mathcal{I}_h(\mathscr{T}^n_h,\mathscr{T}_c)}{\partial \pmb{\xi}_j}\Big{)}^T,\quad j=1,\cdots,N_v,
\end{split}
\end{equation}
where ${\partial \mathcal{I}_h }/{\partial \bm{\xi}_j}$ is considered as a row vector, $\tau>0$ is a parameter used to adjust the response time of mesh movement to the changes in $\mathbb{M}$.
Using the notion of scalar-by-matrix differentiation \cite{Huang2015}, we can rewrite \eqref{MM} as
\begin{equation}\label{xim}
\begin{split}
\frac{d\bm{\xi}_j}{dt}=\frac{\hbox{det}(\mathbb{M}(\bm{x_j}))^{\frac{1}{2}} }{\tau}\sum_{K\in\omega_j}|K|\bm{v}^K_{j_K},\quad j=1,\cdots,N_v,
\end{split}
\end{equation}
where $\omega_j$ is the element patch associated with the vertex $\bm{x}_j$, $j_K$ is the local index of $\bm{x}_j$ in $K$, and $\bm{v}^K_{j_K}$ is the local velocity contributed by $K$.
Denote the edge matrices of $K$ and $K_c$ by $E_K=[\bm{x}_1^K-\bm{x}_0^K,\; \bm{x}_2^K-\bm{x}_0^K]$
and $E_{K_c}=[\bm{\xi}_1^K - \bm{\xi}_0^K,\; \bm{\xi}_2^K -\bm{\xi}_0^K ]$, respectively.
Let $\mathbb{J}=(F'_K)^{-1} = E_{K_c}E_K^{-1}$ and define
\[
G(\mathbb{J},\hbox{det}(\mathbb{J}))=
\frac{1}{3}\hbox{det}(\mathbb{M}_{K})^{\frac{1}{2}}
(\hbox{tr}(\mathbb{J}\mathbb{M}_{K}^{-1}\mathbb{J}^T))^{2}
+\frac{4}{3}\hbox{det}(\mathbb{M}_{K})^{\frac{1}{2}}
\left (\frac{\hbox{det}(\mathbb{J})}{\hbox{det}(\mathbb{M}_{K})^{\frac{1}{2}}} \right )^{2}.
\]
It is not difficult \cite{Huang2015} to find the derivatives of $G$ with respect to $\mathbb{J}$ and $\hbox{det}(\mathbb{J})$ as
\begin{equation}
\begin{split}
&\frac{\partial G}{\partial\mathbb{J}}=
\frac{4}{3}\hbox{det}(\mathbb{M}_{K})^{\frac{1}{2}}
(\hbox{tr}(\mathbb{J}\mathbb{M}_K^{-1}\mathbb{J}^T))\mathbb{M}_K^{-1}\mathbb{J}^T,
\\&\frac{\partial G}{\partial \hbox{det}(\mathbb{J})}=
\frac{8}{3}\hbox{det}(\mathbb{M}_{K})^{-\frac{1}{2}}\hbox{det}(\mathbb{J}).
\end{split}
\end{equation}
Then, the local velocities are expressed as
\begin{equation}
\begin{split}
\left[
  \begin{array}{c}
    ( \bm{v}_1^K  )^T  \\
    ( \bm{v}_2^K  )^T   \\
   \end{array}
 \right]
=
- E_K^{-1}\frac{\partial G }{\partial \mathbb{J} }
- \frac{\partial G}{\partial \hbox{det}(\mathbb{J})}\frac{ \hbox{det}( E_{K_c} )}{\hbox{det}(E_K)} E_{K_c}^{-1},
\quad \bm{v}^K_0=-\bm{v}^K_1-\bm{v}^K_2 .
\end{split}
\end{equation}
It is noted that the velocities for boundary nodes should be modified properly so that they either stay fixed
(such as corner points) or slide on the boundary.

The mesh equation \eqref{xim} can be integrated from $t^n$ to $t^{n+1}$, starting with the reference
computational mesh $\hat{\mathscr{T}}_c$ as the initial mesh. The obtained new mesh is denoted by
$\mathscr{T}_c^{n+1}$. Note that $\mathscr{T}_h^n$ is kept fixed during the integration
and it and $\mathscr{T}_c^{n+1}$ form a correspondence, i.e., $\mathscr{T}_h^n=\Phi_h(\mathscr{T}_c^{n+1})$.
Then the new physical mesh $\mathscr{T}_h^{n+1}$ is defined as $\mathscr{T}_h^{n+1}=\Phi_h(\hat{\mathscr{T}}_c)$
which can be readily computed using linear interpolation.

To conclude this section we would like to point out that a number of other moving mesh methods
have been developed in the past and there is a vast literature in the area. The interested reader
is referred to the books/review articles \cite{Bai94a,Baines-2011,BHR09,Huang2011,Tan05}
and references therein. Also see recent applications \cite{wise2017,Fei-Zhang,zhang2017}.

\section{Numerical examples}\label{sec:numerical}

In this section we present numerical results obtained with the moving mesh DG method described
in the previous sections for a number of one- and two-dimensional examples for the RTE.
Unless otherwise stated, we use in the discrete-ordinate approximation
the Gauss-Legendre $P_8$ and the Legendre-Chebyshev $P_8$-$T_8$
rules for one- and two-dimensional problems, respectively, and
take the final time as $T = 0.1$.
In the computation, $\Delta t=10^{-3}$ is used, which is sufficiently small so that the spatial error dominates
the total error. For mesh movement, we take $\tau=0.1$ for smooth examples such as
Example \ref{Ex1-1d} and Example \ref{Ex1-2d} and $\tau=0.01$ for the others.
For the cases having an exact solution, the error in the computed solution is measured
in the (global) $L^1$, $L^2$, and $L^\infty$ norm, i.e.,
\begin{align*}
 \int_0^T \|e_h(\cdot,t)\|_{L^1}dt,
 \quad
 \int_0^T \|e_h(\cdot,t)\|_{L^2}dt,
 \quad
 \int_0^T \|e_h(\cdot,t)\|_{L^\infty}dt.
\end{align*}

\begin{example}\label{Ex1-1d}
(An accuracy test of the one-dimensional unsteady RTE for the absorbing-scattering model.)
\end{example}

\begin{figure}[H]
\centering
\subfigure[$P^1$-DG]{
\includegraphics[width=0.45\textwidth]{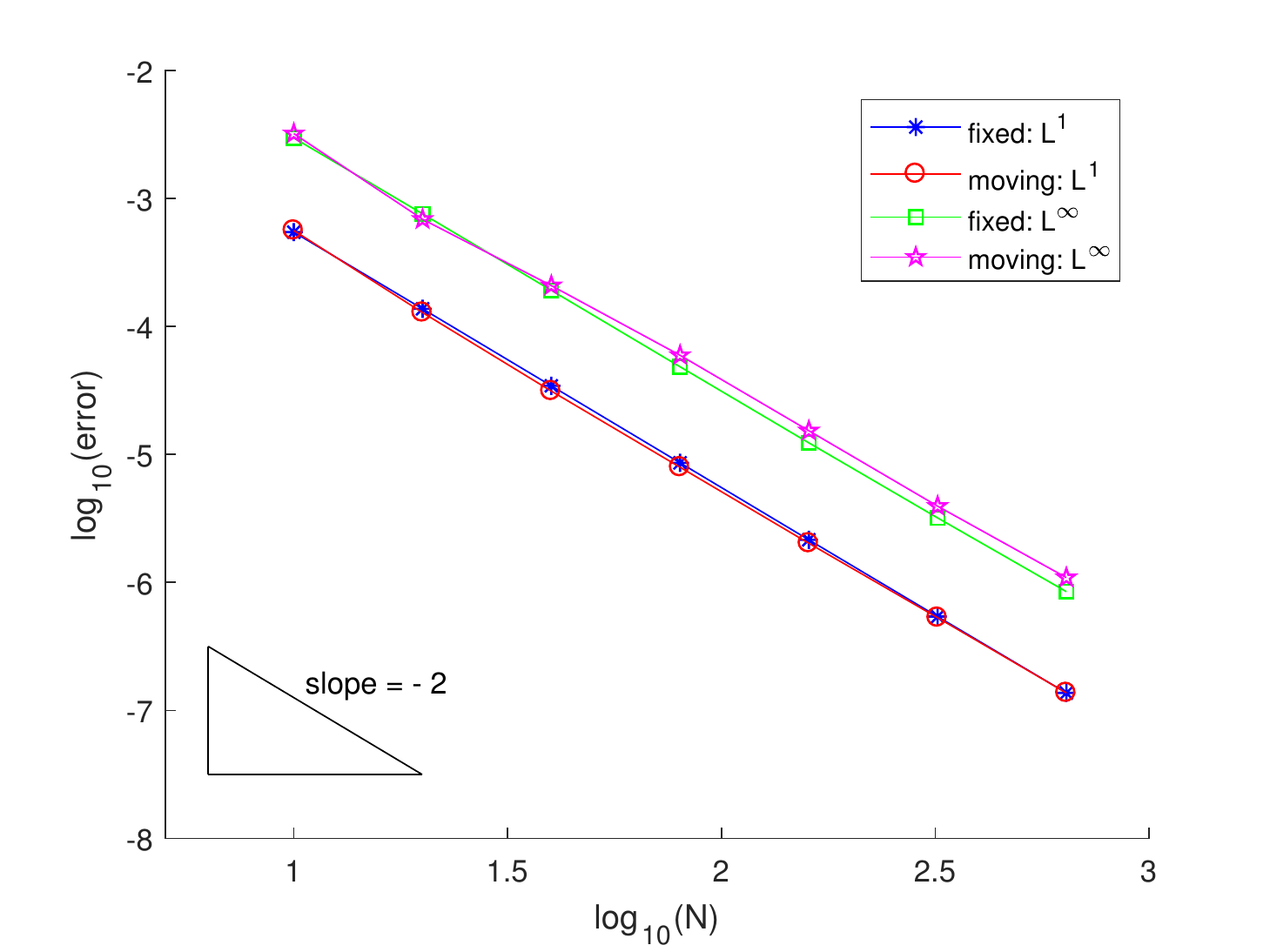}}
\subfigure[ $P^2$-DG]{
\includegraphics[width=0.45\textwidth]{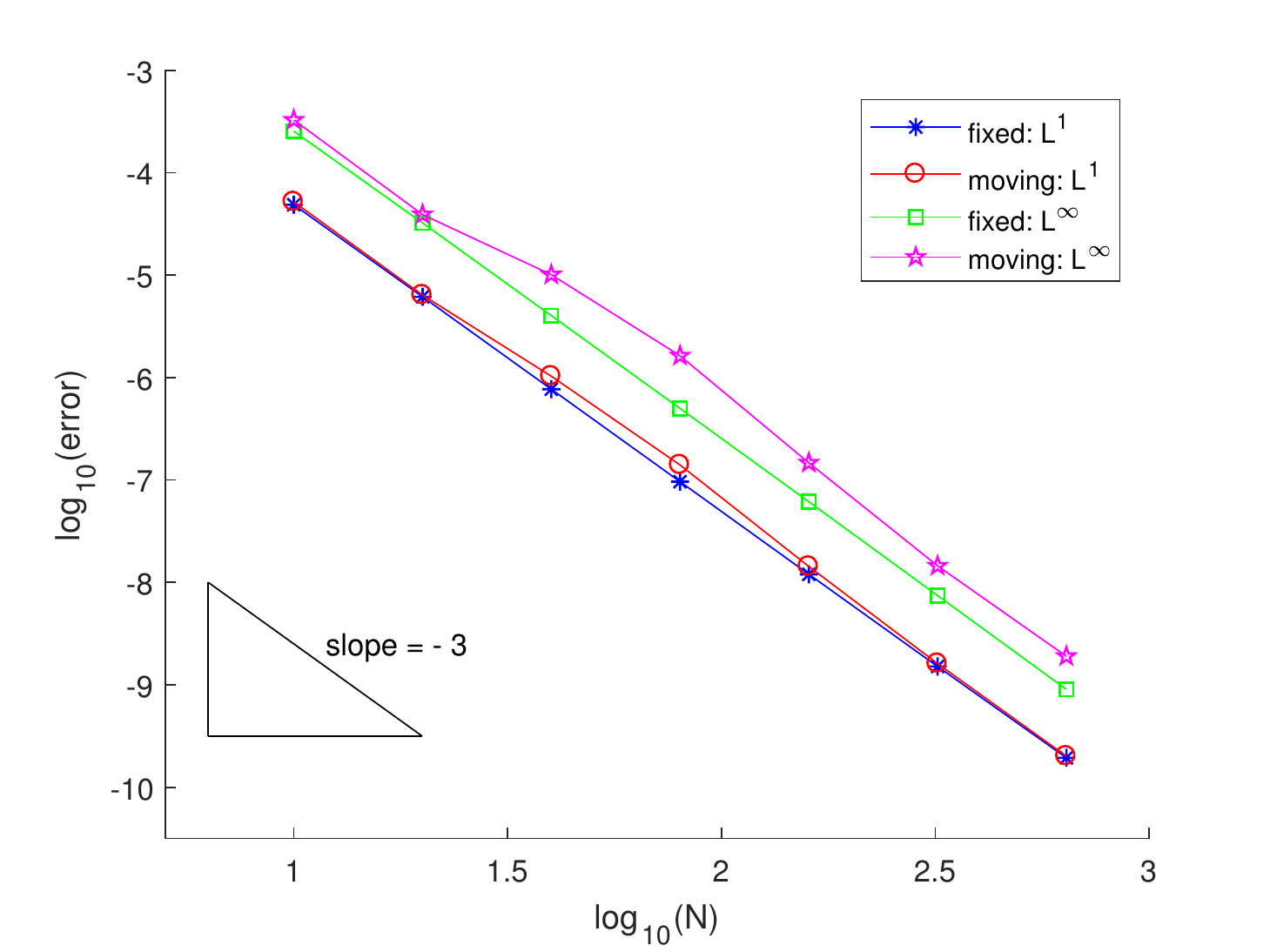}}
\caption{
\small{Example \ref{Ex1-1d}. The $L^1$ and $L^\infty$ norm of the error with a moving and fixed meshes.}
}\label{Fig:d1Ex1-order}
\end{figure}

\noindent
In this example we take $\sigma_t=22000$, $\sigma_s=1$, $c=3.0\times10^{8}$,  and
\begin{align*}
 q(x,\mu,t) = & -4\pi\mu^2\cos^3(\pi(x+t))\sin(\pi(x+t))(\frac{1}{c}+\mu)
 \\
 & +(\sigma_t\mu^2-\frac{\sigma_s}{3})\cos^4(\pi(x+t))
+\sigma_t-\sigma_s .
\end{align*}
The initial and boundary conditions are given by
\begin{equation*}
\begin{split}
&I(x,\mu,0)=\mu^2\cos^4(\pi x)+1,\quad\quad\quad \text{ for } -1 <\mu < 1,\;  0< x < 1 ,
\\
&I(0,\mu,t)=\mu^2\cos^4(\pi t)+1,\quad\quad\quad~ \text{ for }  0 <\mu \leq 1,\;  0< t\leq 0.1 ,
\\&I(1,\mu,t)=\mu^2\cos^4(\pi(1+t))+1,\quad \text{ for }  -1 \leq\mu < 0,\; 0< t\leq 0.1.
\end{split}
\end{equation*}
The problem has the exact solution as $I(x,\mu,t)=\mu^2\cos^4(\pi (x+t)) + 1$.
The $L^1$ and $L^\infty$ norm of the error in the numerical solutions obtained
with the $P^1$-DG and $P^2$-DG methods with fixed and
moving meshes is shown in Fig.~\ref{Fig:d1Ex1-order}.
It can be seen that both types of mesh lead to the same convergence order (2nd for $P^1$-DG and
3rd for $P^2$-DG) and comparable errors, which is consistent with the theoretical prediction.
One may also see that the solution error associated with a moving mesh is slightly larger than
that associated with a uniform mesh of the same size for this example. This can happen
for problems with smooth solutions
since the adaptation strategy used in the moving mesh method does not directly minimize
the error associated with the implicit Euler--DG discretization. Moreover, no special effort
has been made to optimize the parameters $\alpha_{h,m}$ in (\ref{mer}) and
$\tau$ in (\ref{xim}) for this specific example.

\begin{example}\label{Ex7-1d}
(A discontinuous test for the one-dimensional unsteady RTE simulating the absorbing-scattering model.)
\end{example}

\noindent
We take extinction coefficient $\sigma_t=1000$, $\sigma_s=1$, $c=3.0\times10^{8}$, and
\begin{align*}
q(x,y,\zeta,\eta,t)= & \sigma_t\sin(2\pi(\tanh(Rx)+5\mu t))
+\frac{10\pi\mu}{c}\cos(2\pi(\tanh(Rx)+5\mu t))
\\&+2\pi\mu R\cos(2\pi(\tanh(Rx)+5\mu t))(1-\tanh^2(Rx))
\\&+\frac{\sigma_s}{20\pi t}(\cos(2\pi(\tanh(Rx)+5t))-\cos(2\pi(\tanh(Rx)-5t)) ) +a(\sigma_t-\sigma_s),
\end{align*}
where $a=2$ and $R=200$.
The initial condition is $I(x,\mu,0)=\sin(2\pi\tanh(Rx))+a $ and the boundary conditions are
\begin{equation*}\label{e3}
\begin{split}
&I(-1,\mu,t)=\sin(2\pi(-\tanh(R)+5\mu t))+a,  \quad 0 <\mu \leq 1 ,~~0< t\leq T ,
\\
&I(1,\mu,t)=\sin(2\pi(\tanh(R)+5\mu t))+a,  \quad\quad -1 \leq\mu < 0,~~0< t\leq T.
\end{split}
\end{equation*}

\begin{figure}[H]
\centering
\includegraphics[width=0.45\textwidth]{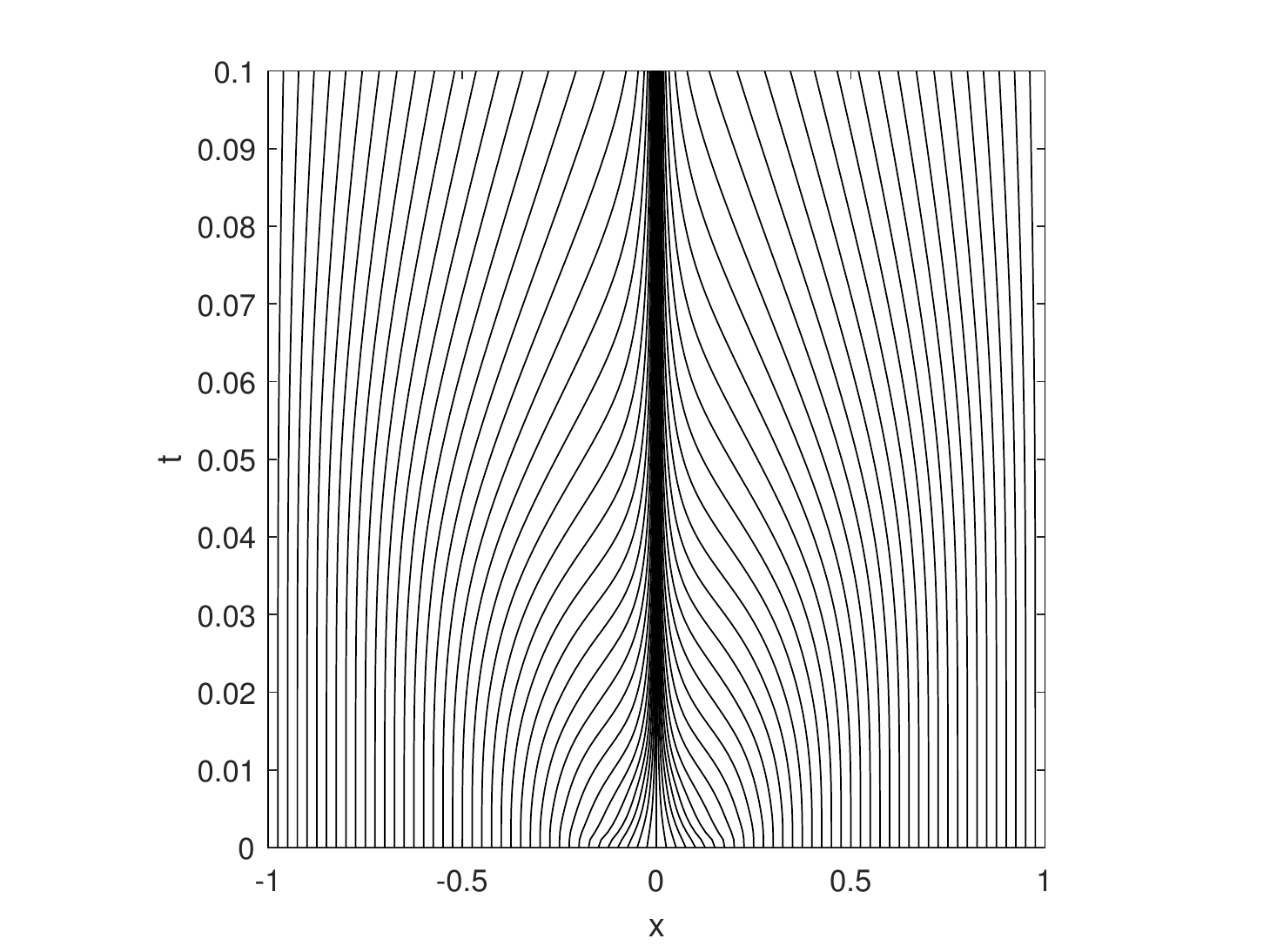}
\caption{
\small{Example \ref{Ex7-1d}. The mesh trajectories are obtained with the $P^2$-DG method with
a moving mesh of $N=80$. }
}\label{Fig:d1Ex7p2-mesh}
\end{figure}

\begin{figure}[H]
\centering
\subfigure[MM $N$=80, FM $N$=80]{
\includegraphics[width=0.45\textwidth]{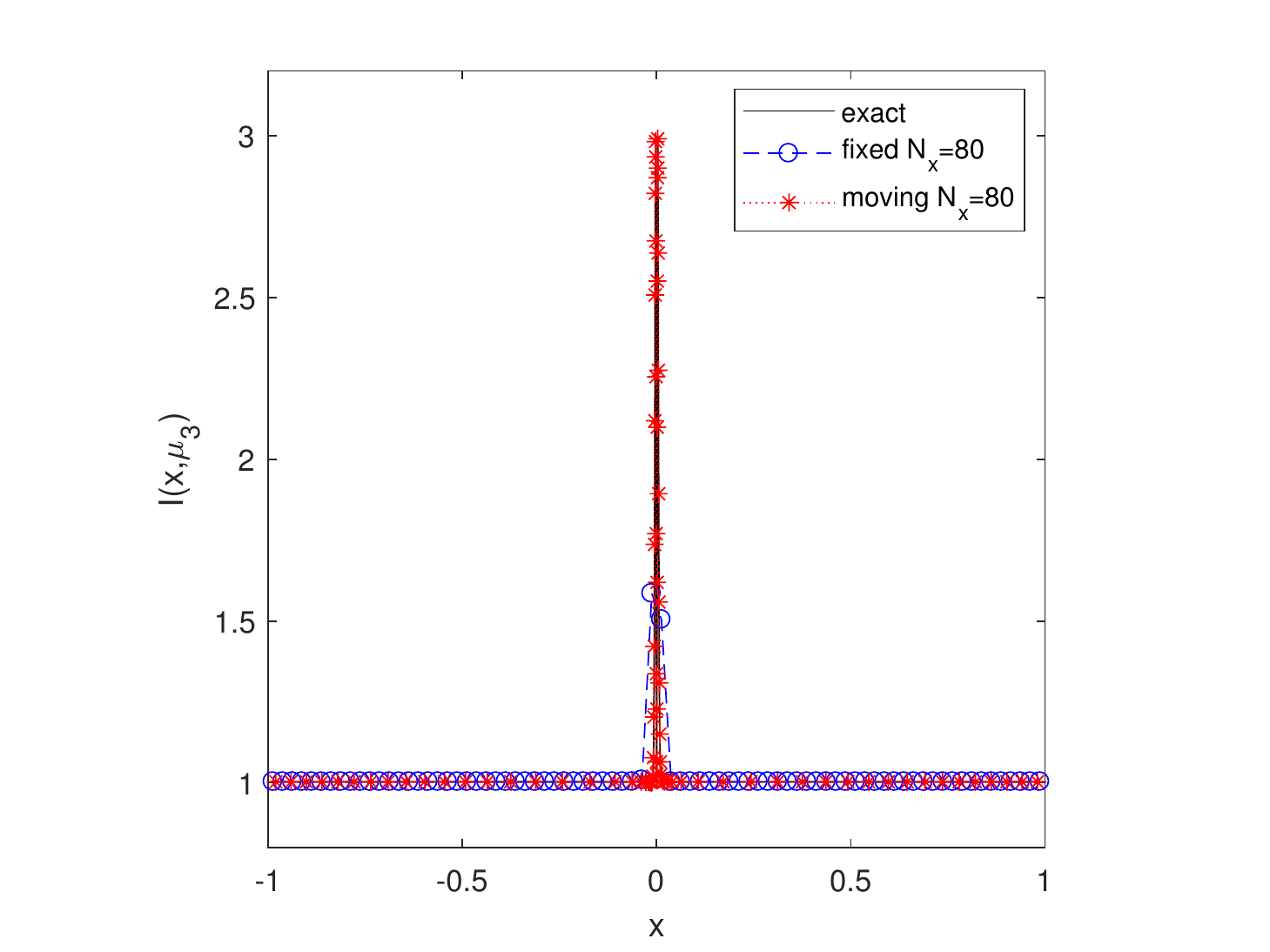}}
\subfigure[Close view of (a) near $x$=0]{
\includegraphics[width=0.45\textwidth]{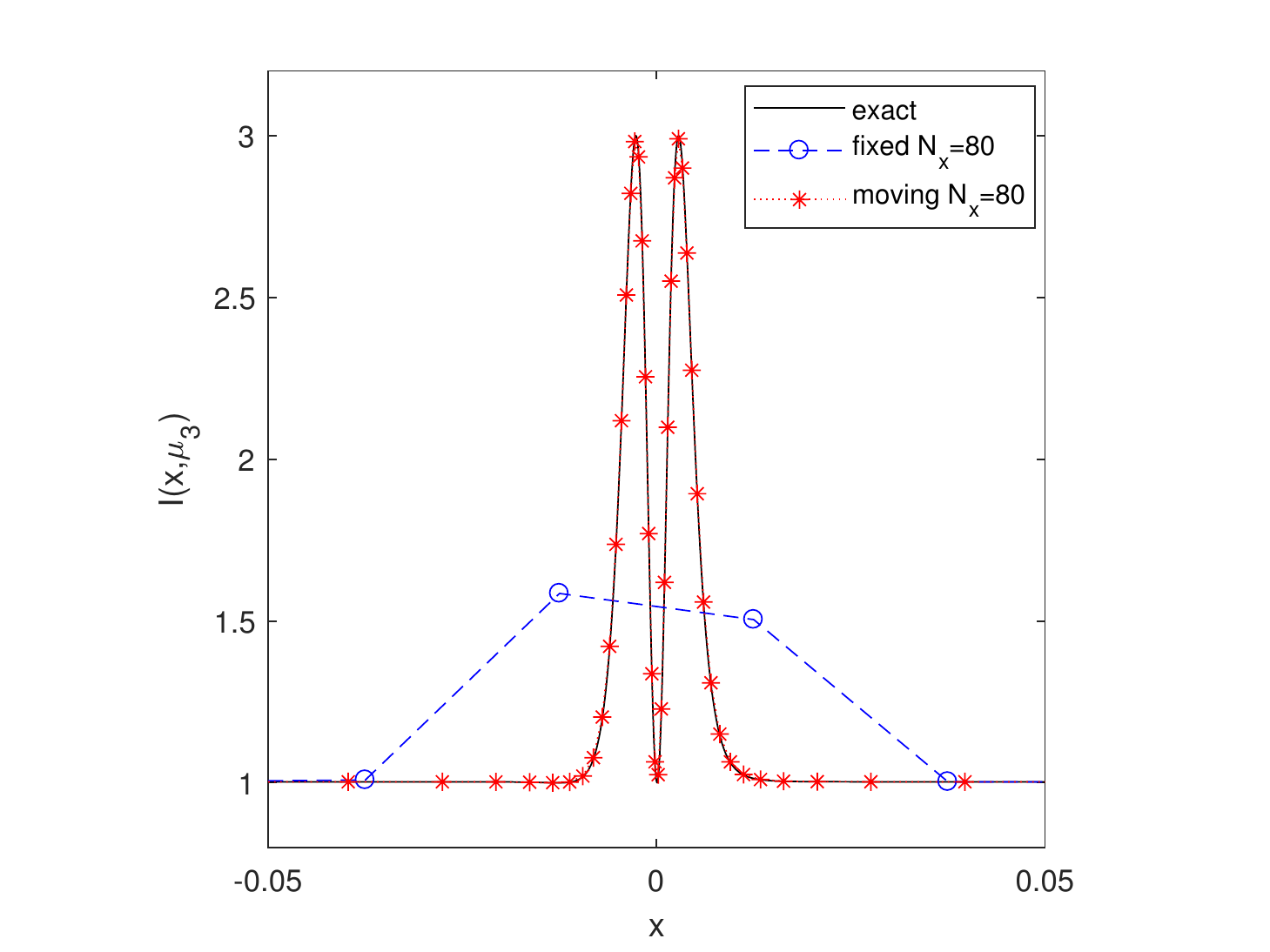}}
\subfigure[MM $N$=80, FM $N$=1280]{
\includegraphics[width=0.45\textwidth]{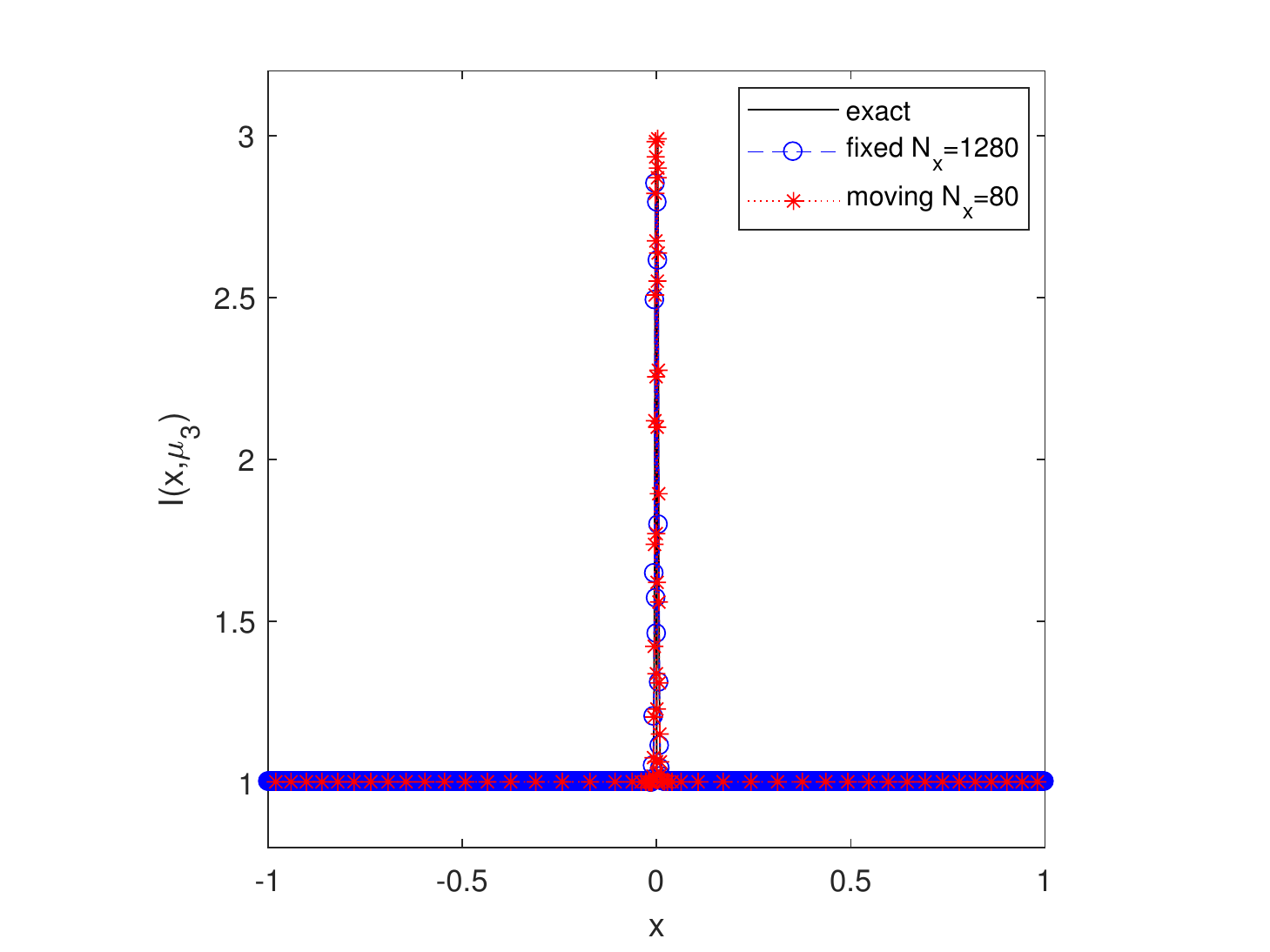}}
\subfigure[Close view of (c) near $x$=0]{
\includegraphics[width=0.45\textwidth]{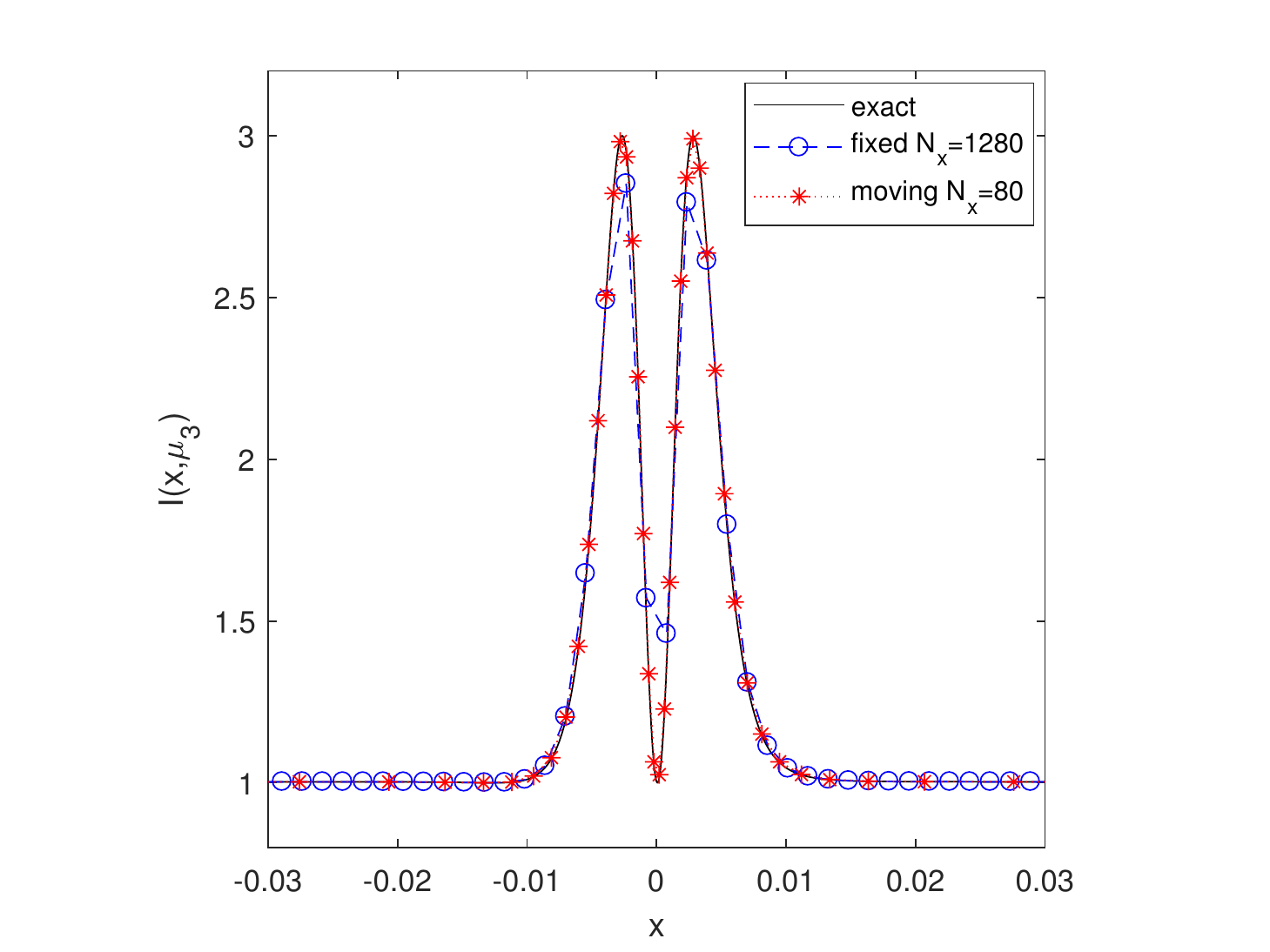}}
\caption{
\small{Example \ref{Ex7-1d}. The solution in the direction $\mu=-0.5255$ obtained with the $P^2$-DG method
with a moving mesh of $N$=80 is compared with those obtained with fixed meshes of $N$=80 and $N$=1280.
MM and FM stand for moving mesh and fixed mesh, respectively.}
}\label{Fig:d1Ex7p2-u3}
\end{figure}

\begin{figure}[H]
\centering
\subfigure[MM $N$=80, FM $N$=80]{
\includegraphics[width=0.45\textwidth]{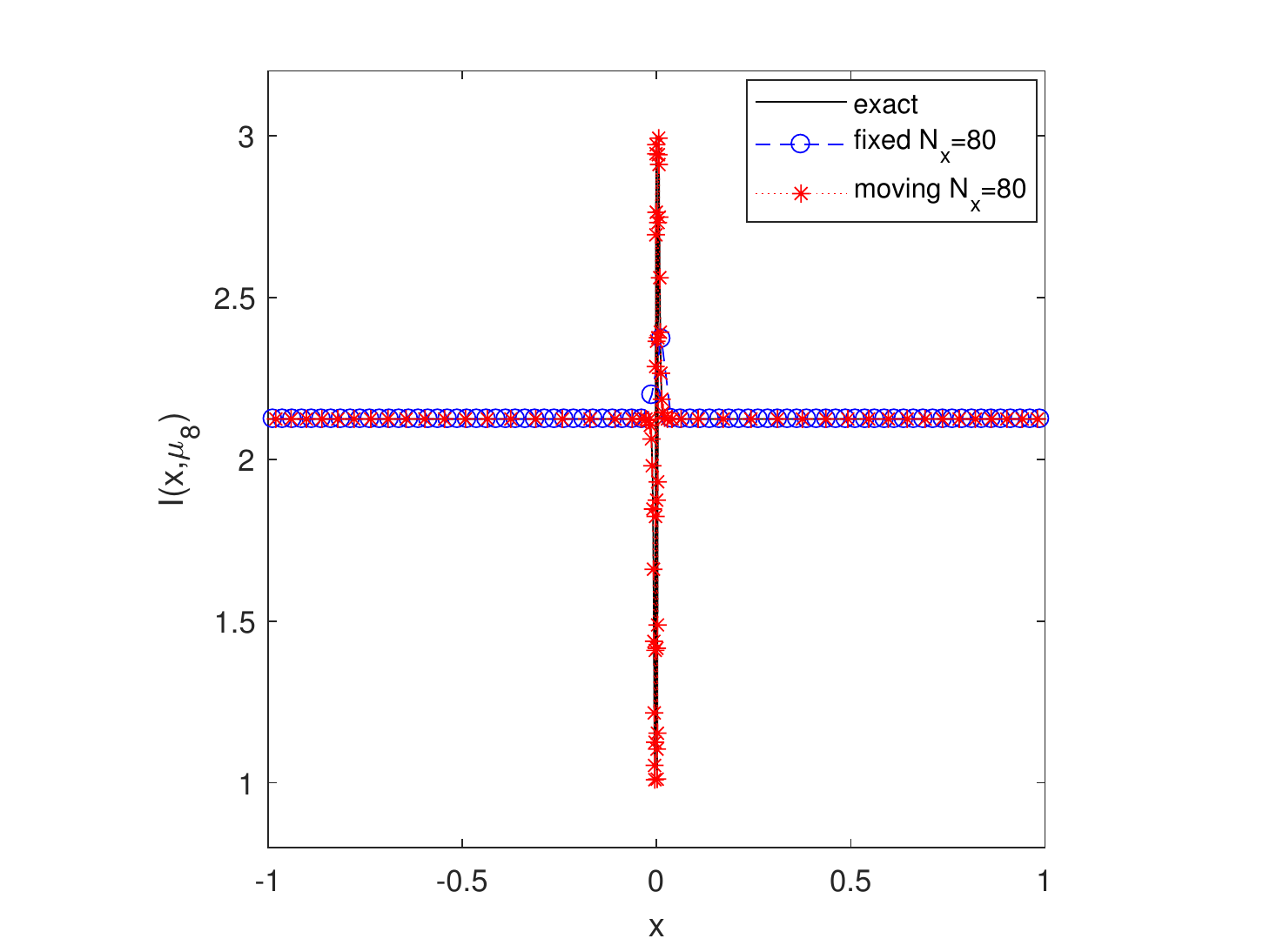}}
\subfigure[Close view of (a) near $x$=0]{
\includegraphics[width=0.45\textwidth]{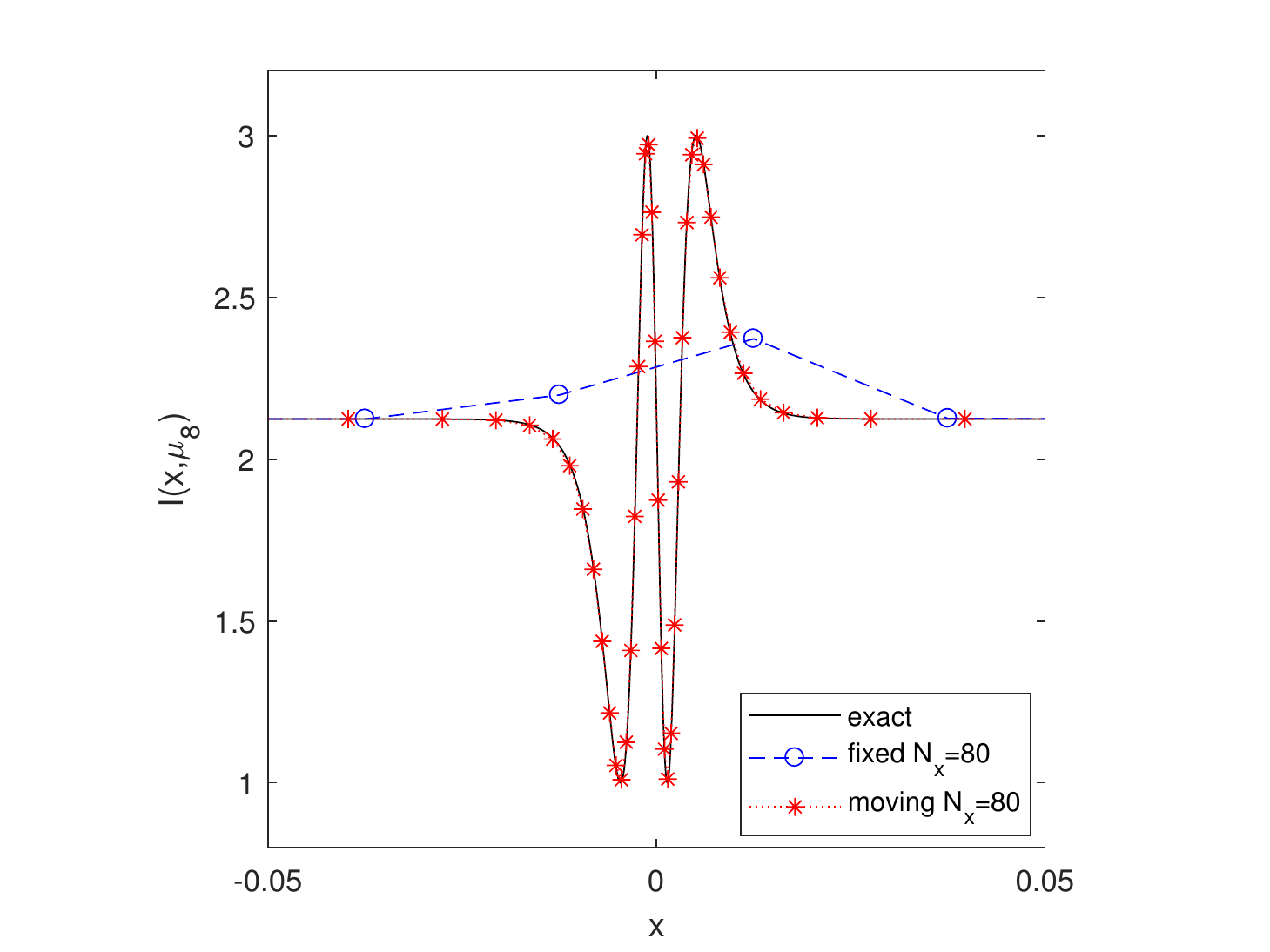}}
\subfigure[MM $N$=80, FM $N$=1280]{
\includegraphics[width=0.45\textwidth]{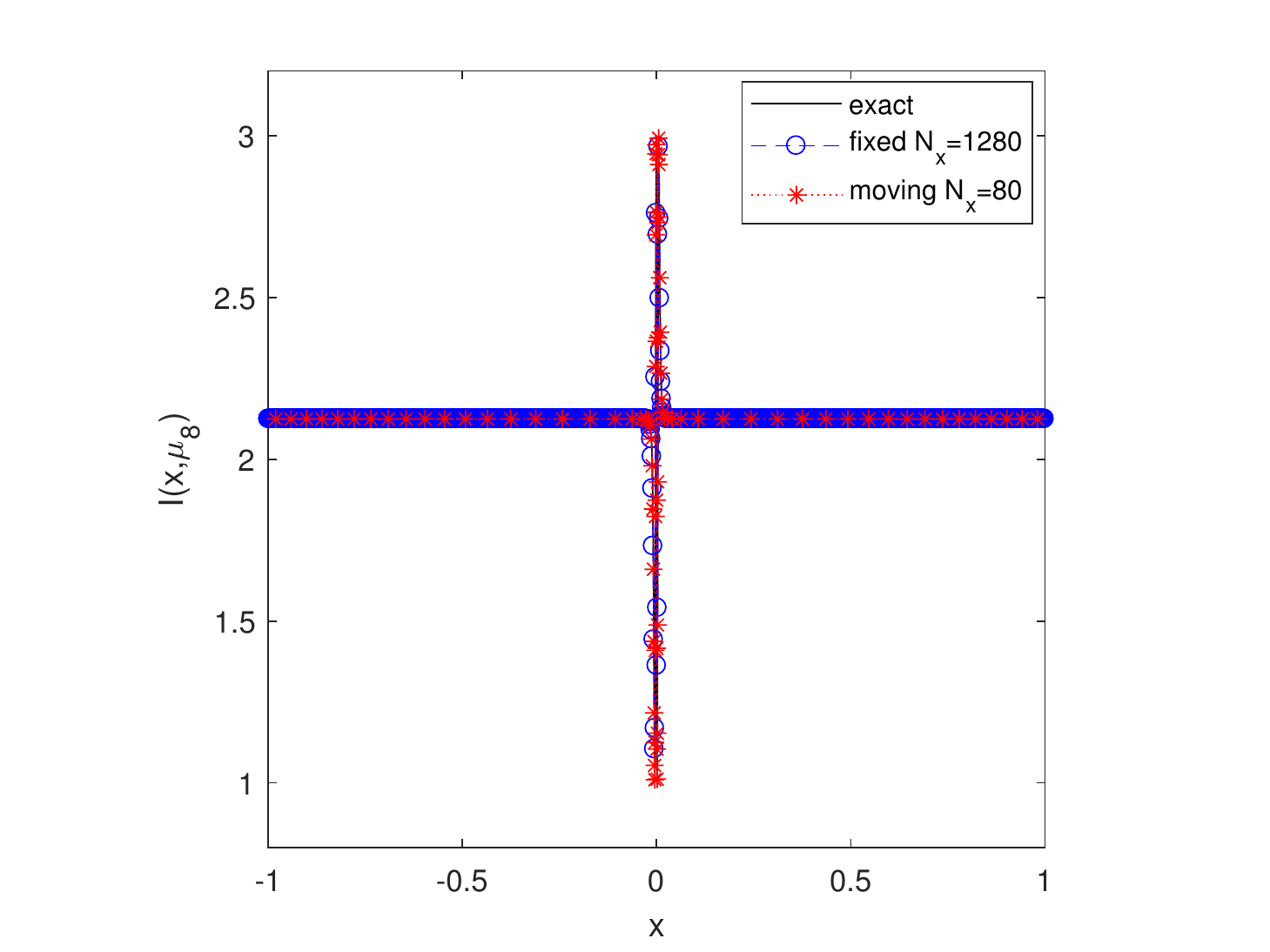}}
\subfigure[Close view of (c) near $x$=0 ]{
\includegraphics[width=0.45\textwidth]{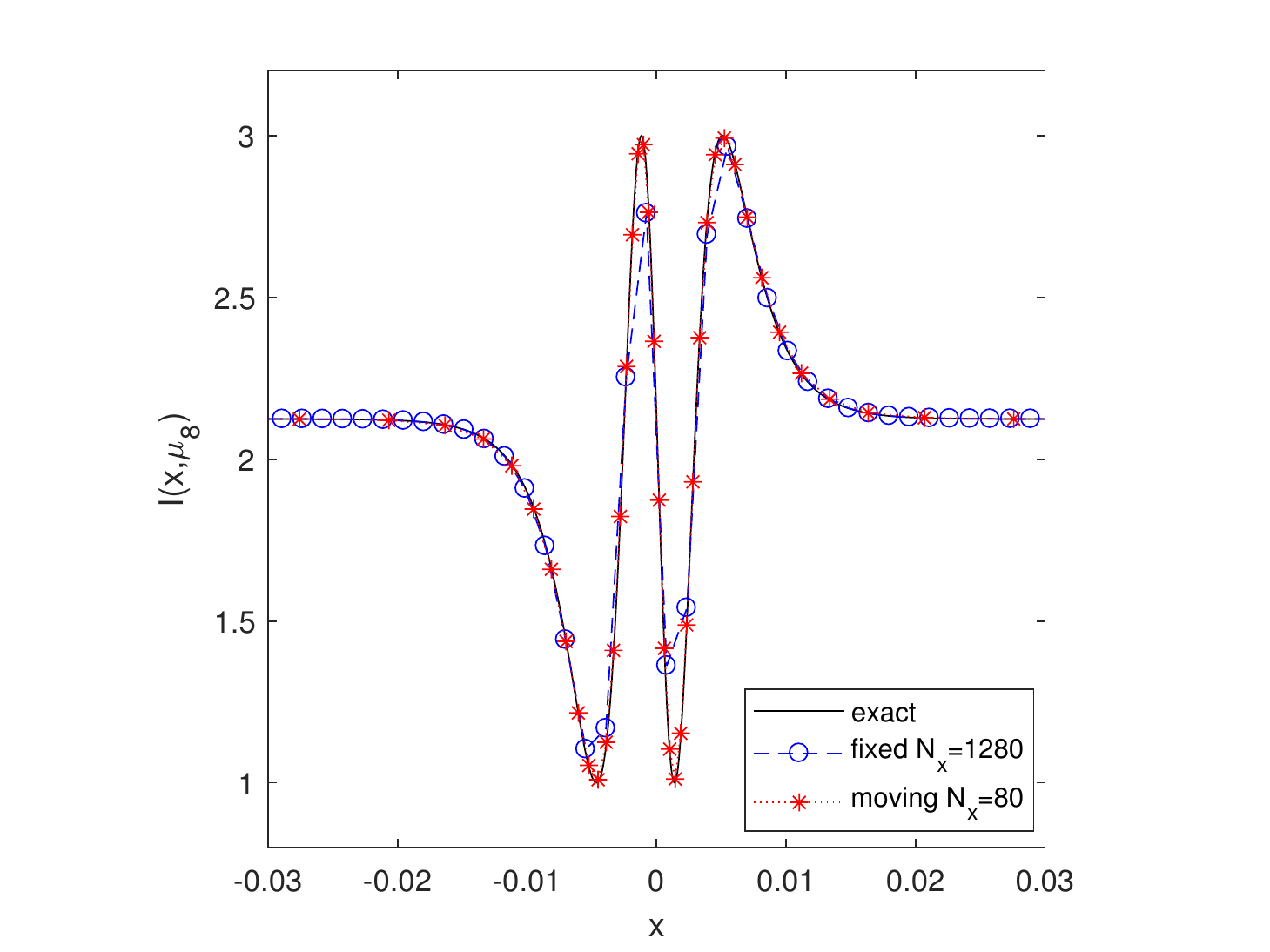}}
\caption{
\small{Example \ref{Ex7-1d}. The solution in the direction $\mu=0.9603$ obtained with the $P^2$-DG method
with a moving mesh of $N$=80 is compared with those obtained with fixed meshes of $N$=80 and $N$=1280.}
}\label{Fig:d1Ex7p2-u8}
\end{figure}

The exact solution of this example is $I(x,\mu,t)=\sin(2\pi(\tanh(Rx)+5\mu t))+a$.
The mesh trajectories for the $P^2$-DG method with a moving mesh of $N = 80$ are shown
in Fig.~\ref{Fig:d1Ex7p2-mesh}. The moving mesh solution ($N=80$) in the direction $\mu=-0.5255$
is compared with the fixed mesh solutions obtained with $N=80$ and $N=1280$ in
Fig.~\ref{Fig:d1Ex7p2-u3}. Similar results are shown in Fig.~\ref{Fig:d1Ex7p2-u8} for the angular direction
$\mu =0.9603 $. These results show that the moving mesh solution
($N=80$) is more accurate than those with fixed meshes of $N=80$ and $N=1280$.

\begin{figure}[H]
\centering
\subfigure[$P^1$-DG]{
\includegraphics[width=0.45\textwidth]{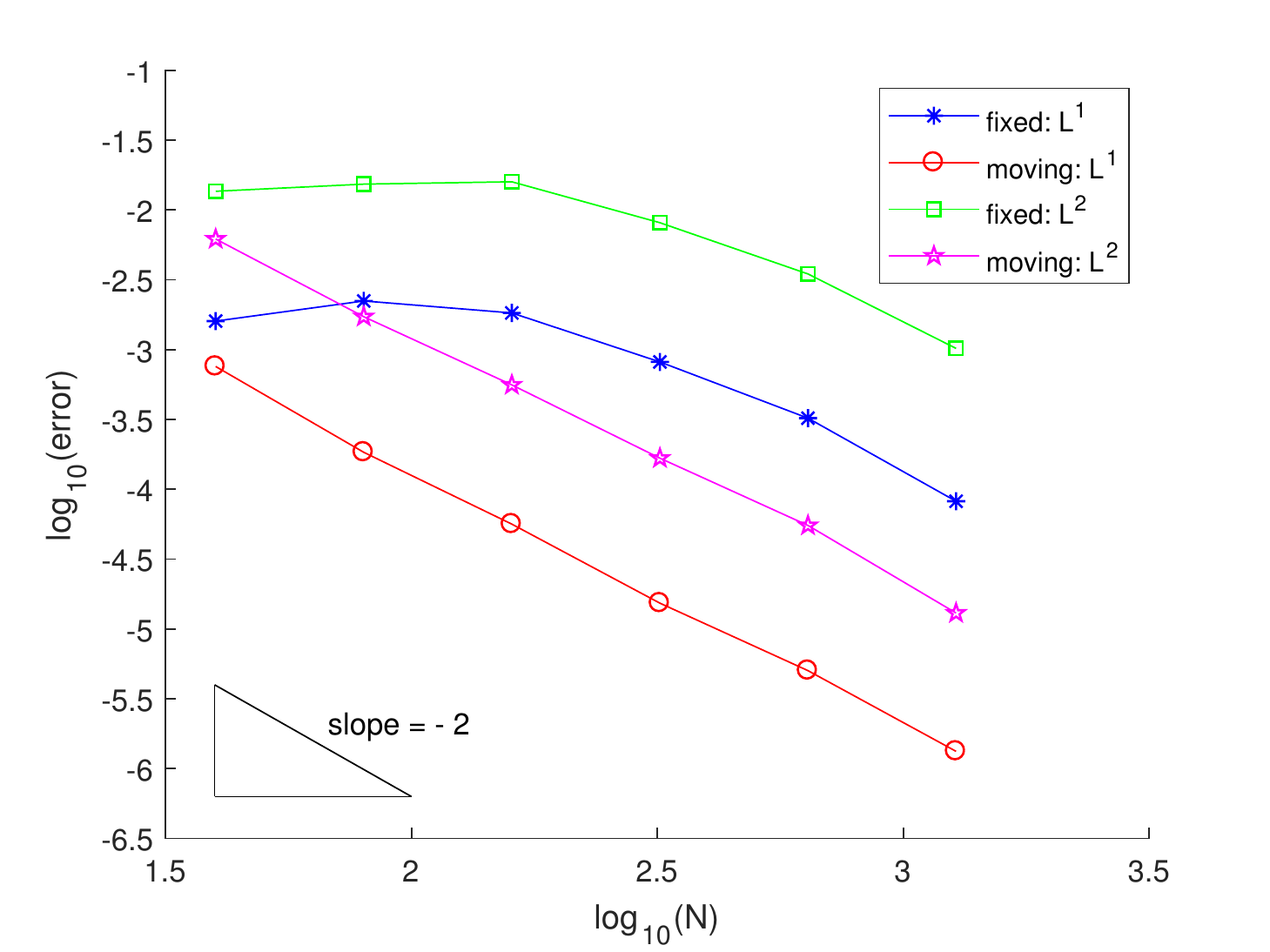}}
\subfigure[ $P^2$-DG]{
\includegraphics[width=0.45\textwidth]{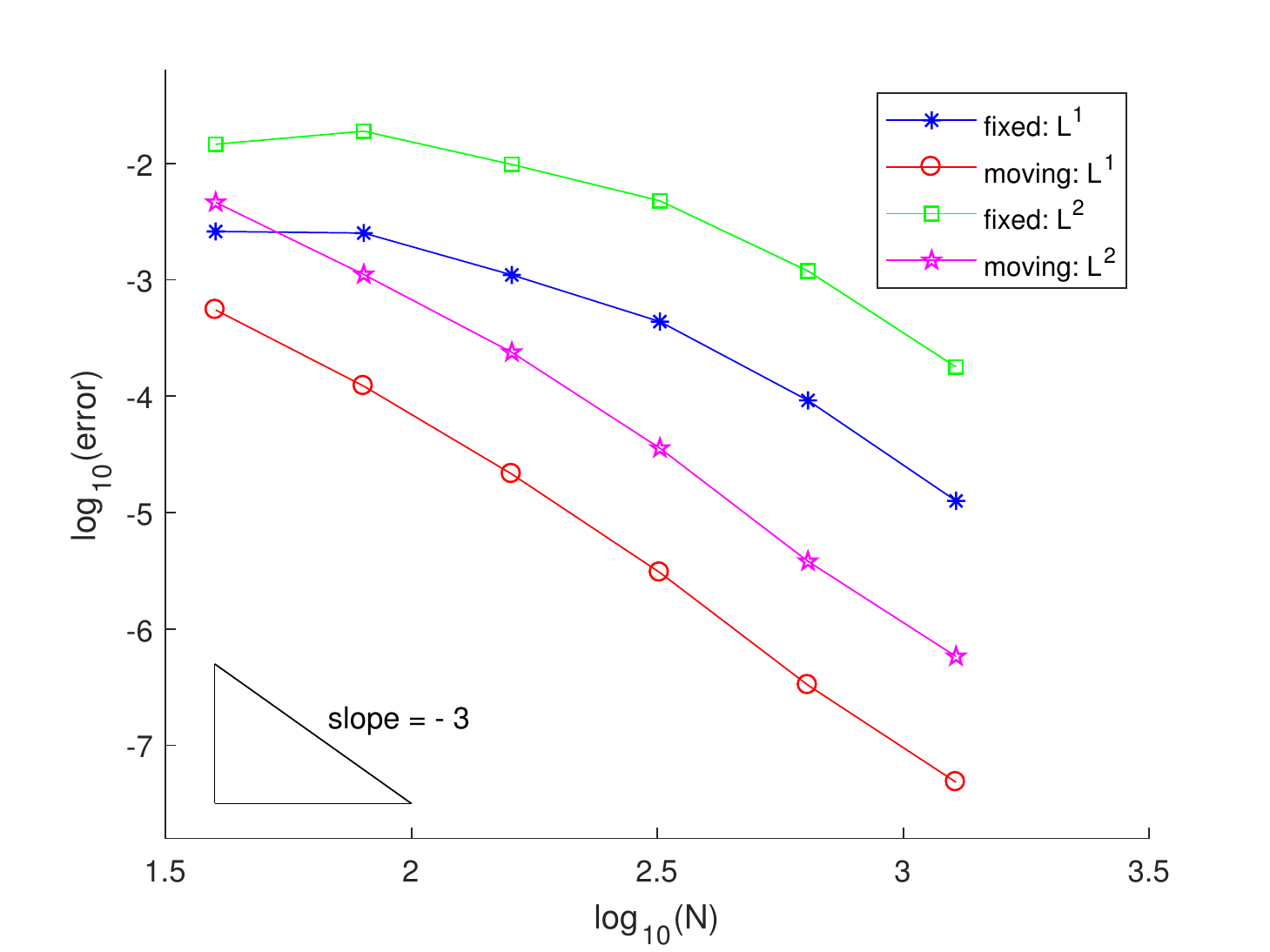}}
\caption{
\small{Example \ref{Ex7-1d}. The $L^1$ and $L^2$ norm of the error with moving and fixed meshes.}
}\label{Fig:d1Ex7-order}
\end{figure}

The error in the $L^1$ and $L^2$ norm is shown in Fig.~\ref{Fig:d1Ex7-order} for both $P^1$-DG and $P^2$-DG methods with fixed and moving meshes. It can be seen that both fixed and moving meshes lead to almost
the same order of convergence for relatively large $N$, i.e., 2nd order for $P^1$-DG and 3rd-order for $P^2$-DG.
However, a moving mesh always produces more accurate solutions than a fixed mesh of the same size
for this example.

\begin{figure}[H]
\centering
\includegraphics[width=0.45\textwidth]{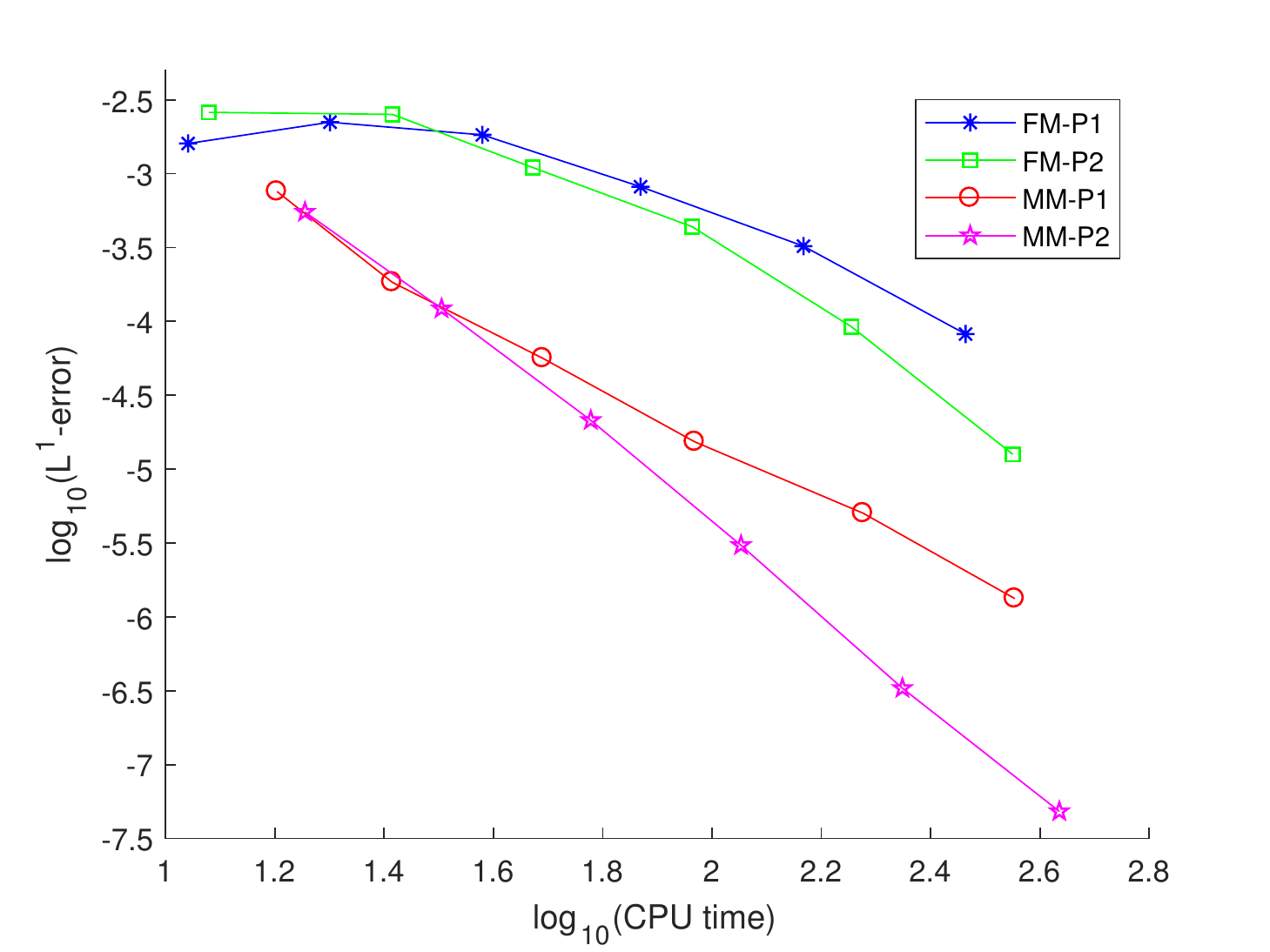}
\caption{
\small{Example \ref{Ex7-1d}. The $L^1$ norm of the error is plotted against the CPU time.}
}\label{Fig:d1Ex7-cpuL1}
\end{figure}

To show the efficiency of the methods, we plot in Fig.~\ref{Fig:d1Ex7-cpuL1}
the $L^1$ norm of the error against the CPU time measured in seconds on a Thinkpad T440 with Matlab 2017a.
One can see that moving mesh $P^1$-DG (resp. $P^2$-DG)
is more efficient than fixed mesh $P^1$-DG (resp. $P^2$-DG) in the sense that
the former leads to a smaller error than the latter for a fixed amount of the CPU time.
Moreover, when $N > 50$, MM (resp. FM) $P^2$-DG is more efficient than MM (resp. FM) $P^1$-DG.
Thus, a moving mesh improves the computational efficiency and
the quadratic DG method has better efficiency than the linear one on both fixed and moving meshes
when $N$ is sufficiently large.


\begin{example}\label{Ex6-1d}
(A discontinuous example of the one-dimensional unsteady RTE for the absorbing-scattering model.)
\end{example}

\begin{figure}[H]
\centering
\includegraphics[width=0.45\textwidth]{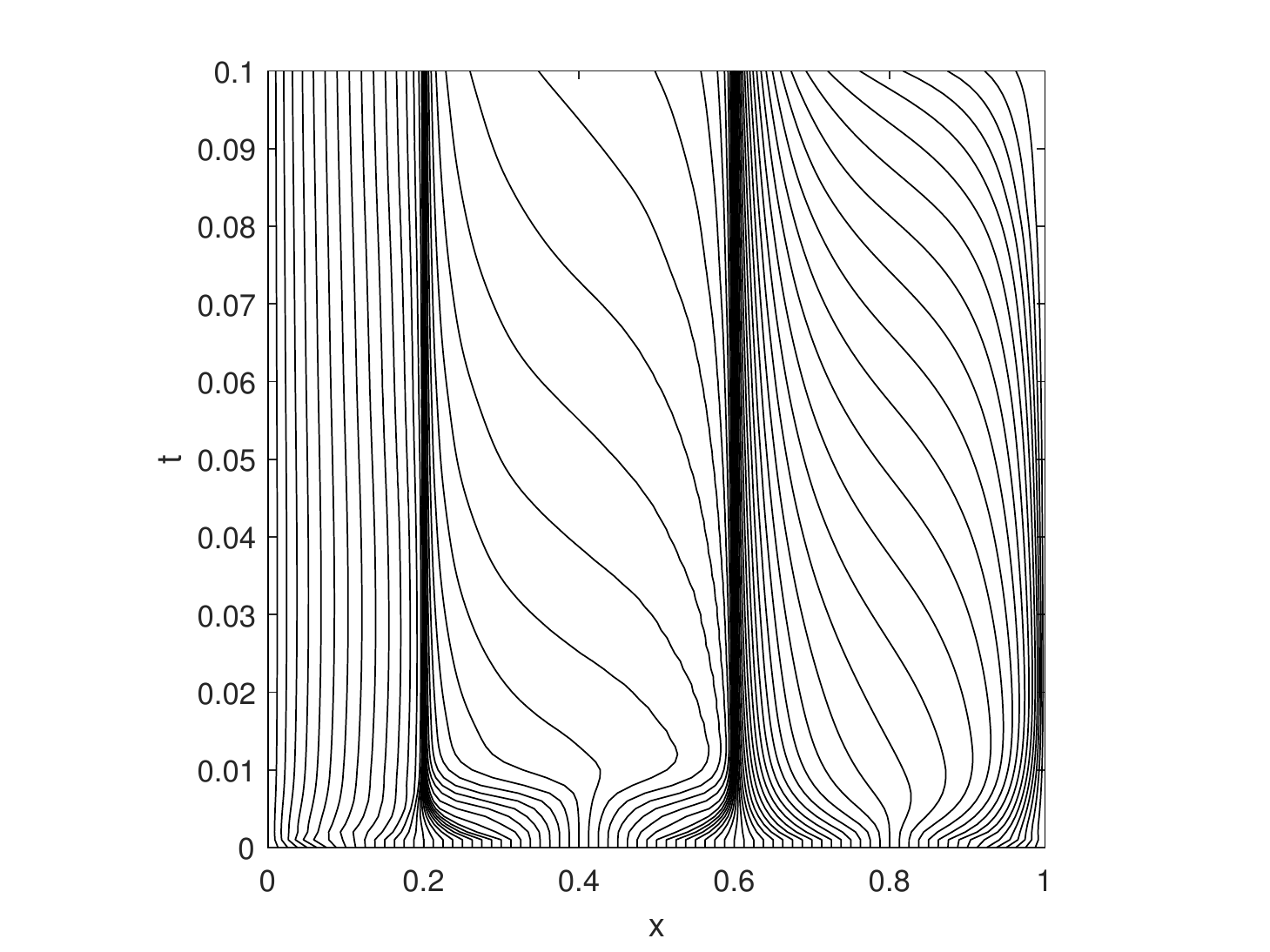}
\caption{
\small{Example \ref{Ex6-1d}. The mesh trajectories are obtained with the $P^2$-DG method with
a moving mesh of $N=80$. }
}\label{Fig:d1Ex6p2-mesh}
\end{figure}

\begin{figure}[H]
\centering
\subfigure[MM $N$=80, FM $N$=80]{
\includegraphics[width=0.45\textwidth]{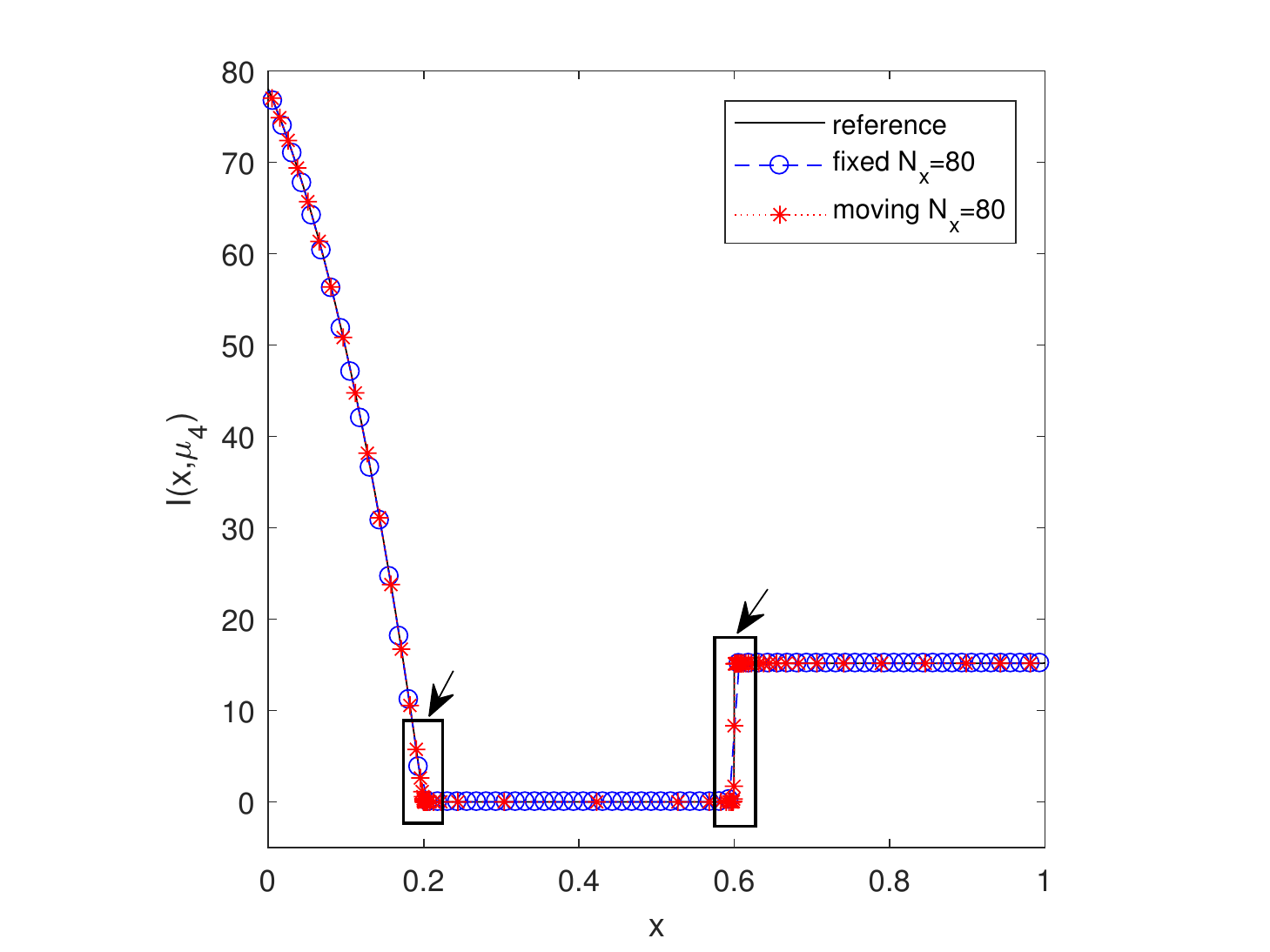}}
\subfigure[Close view of (a) near $x$=0.2 and $x$=0.6]{
\includegraphics[width=0.45\textwidth]{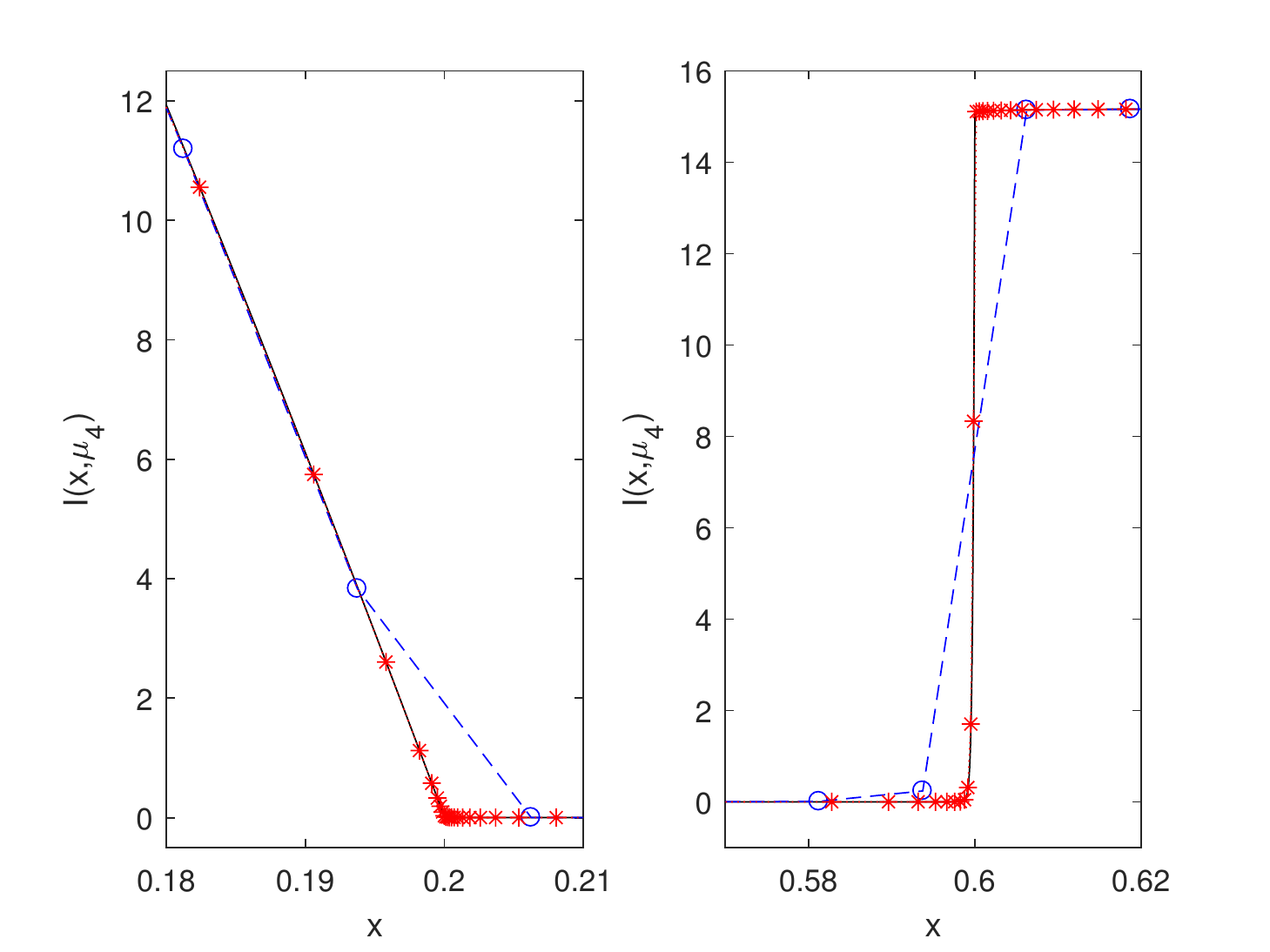}}
\subfigure[MM $N$=80, FM $N$=1280]{
\includegraphics[width=0.45\textwidth]{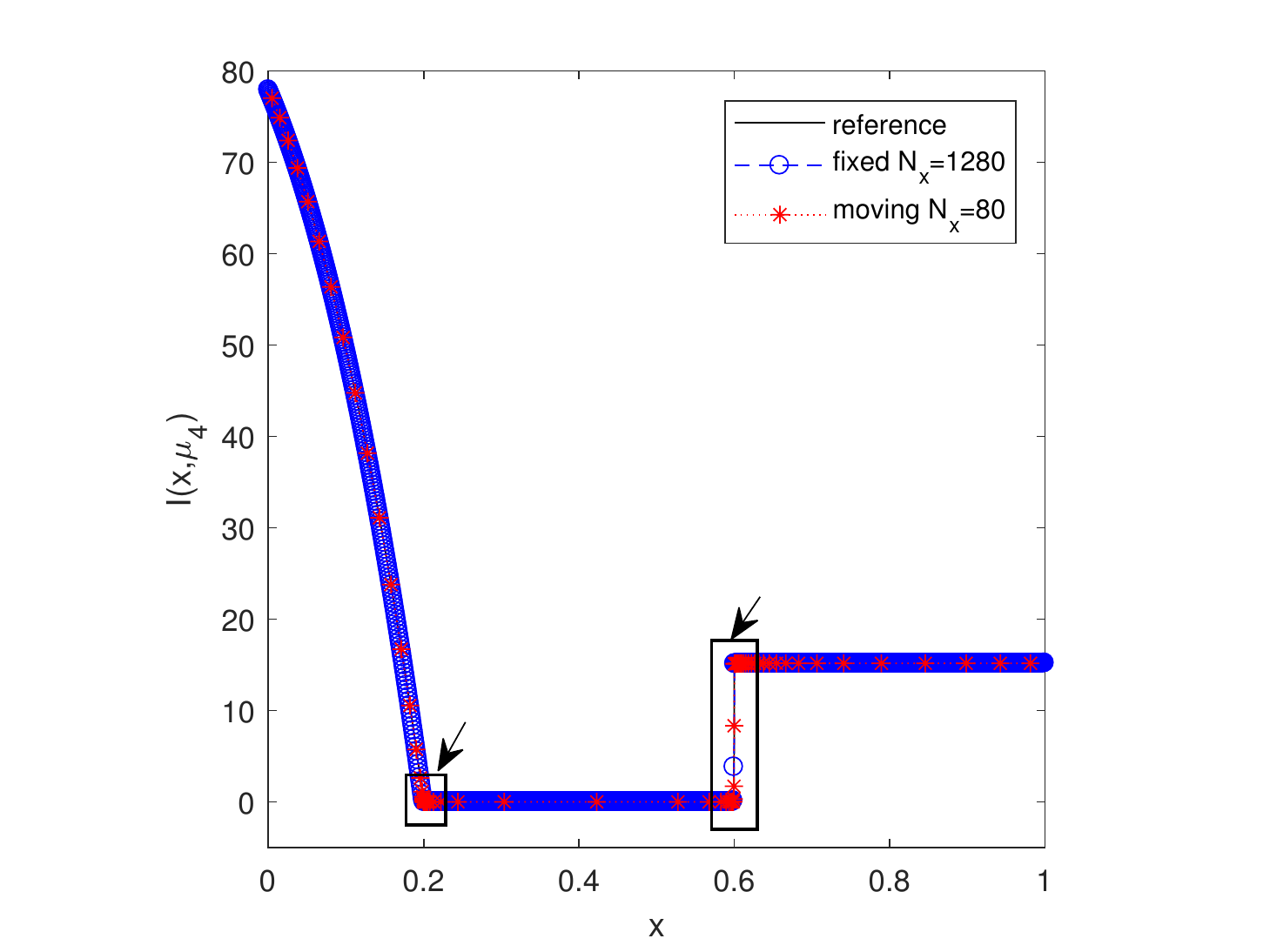}}
\subfigure[Close view of (c) near $x$=0.2 and $x$=0.6]{
\includegraphics[width=0.45\textwidth]{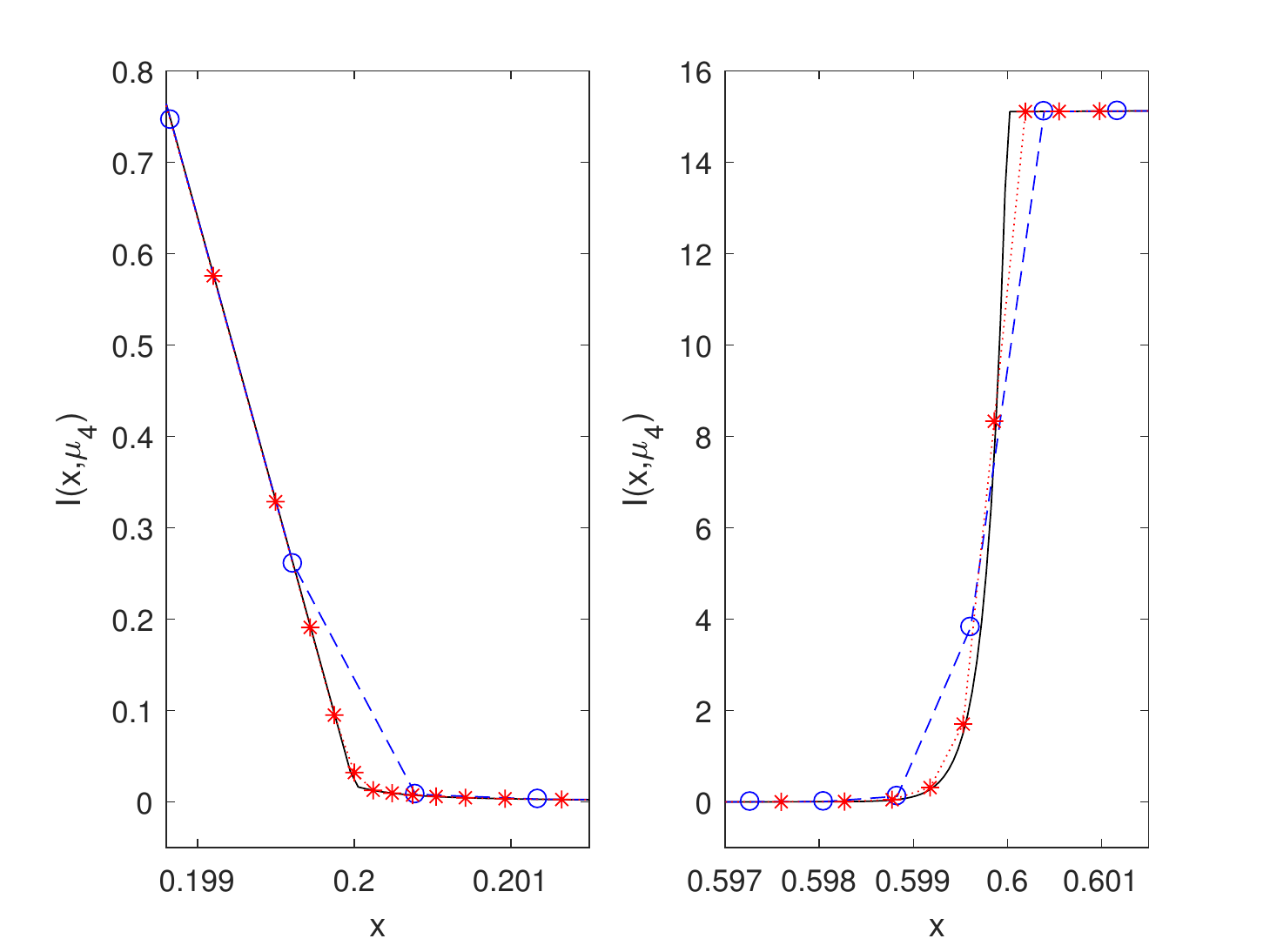}}
\caption{
\small{Example \ref{Ex6-1d}. The solution in the direction $\mu=-0.1834$ obtained with the $P^2$-DG method
with a moving mesh of $N$=80 is compared with those obtained with fixed meshes of $N$=80 and $N$=1280.}
}\label{Fig:d1Ex6p2-u4}
\end{figure}

\begin{figure}[H]
\centering
\subfigure[MM $N$=80, FM $N$=80]{
\includegraphics[width=0.45\textwidth]{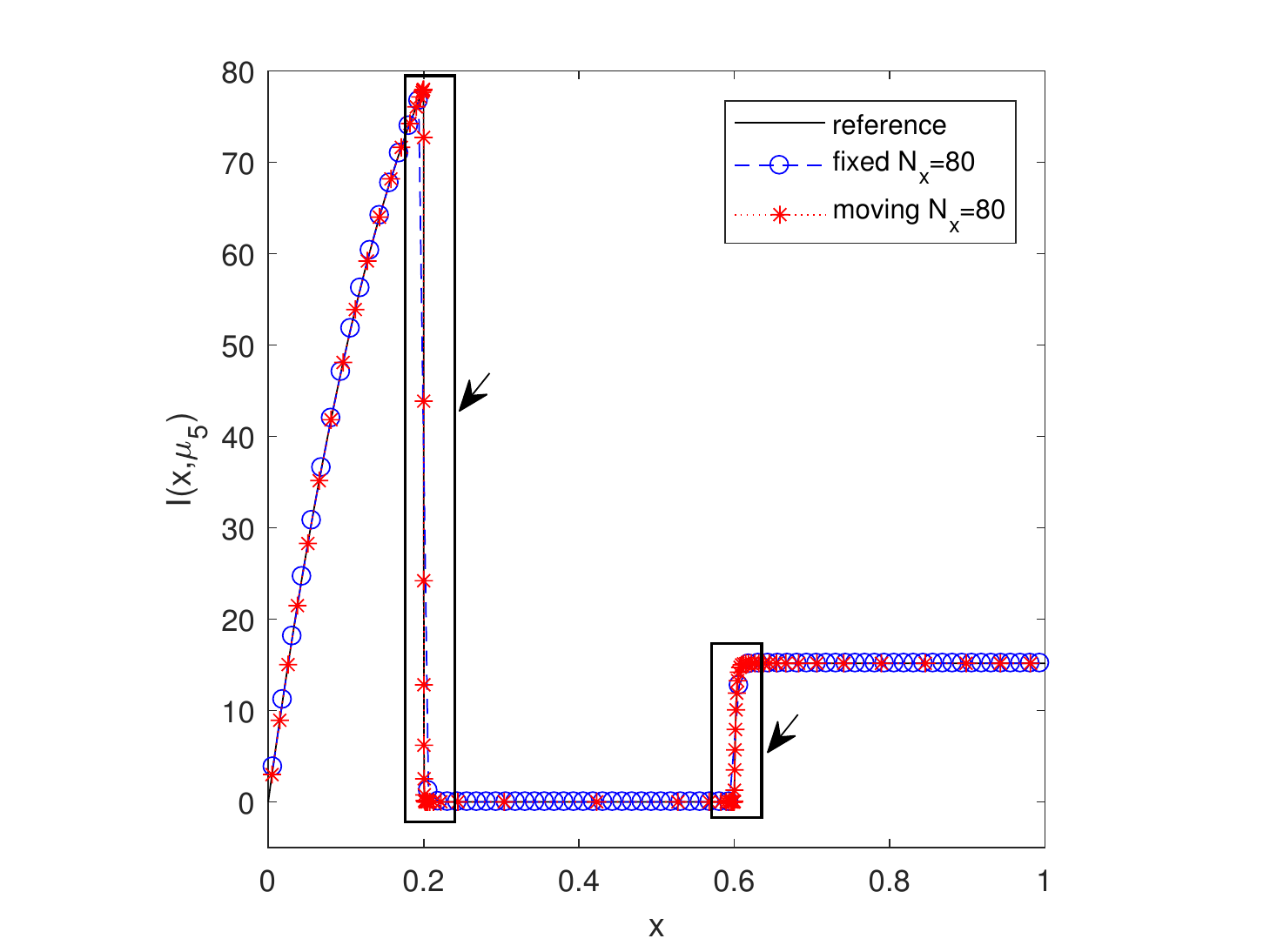}}
\subfigure[Close view of (a) near $x$=0.2 and $x$=0.6]{
\includegraphics[width=0.45\textwidth]{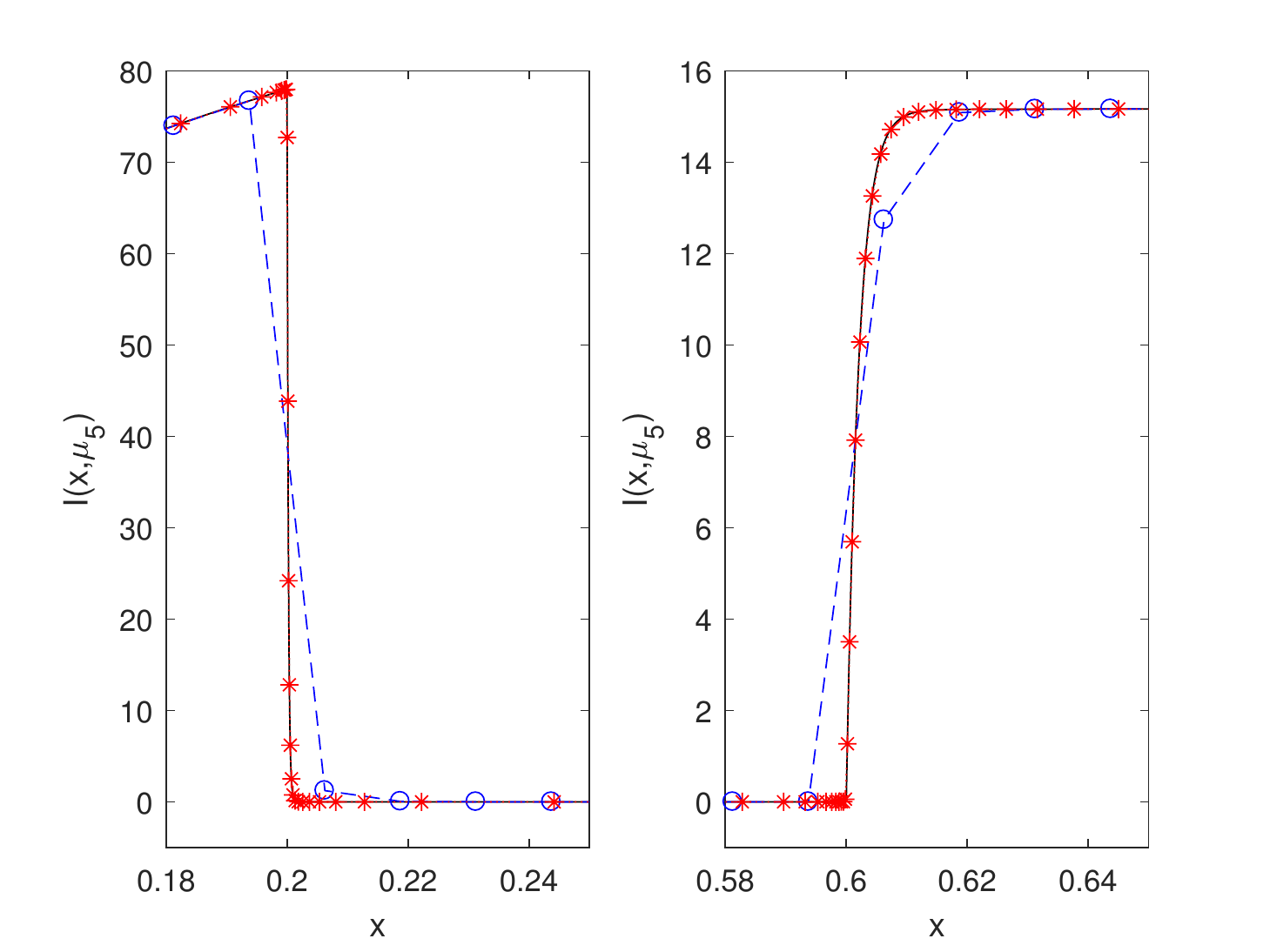}}
\subfigure[MM $N$=80, FM $N$=1280]{
\includegraphics[width=0.45\textwidth]{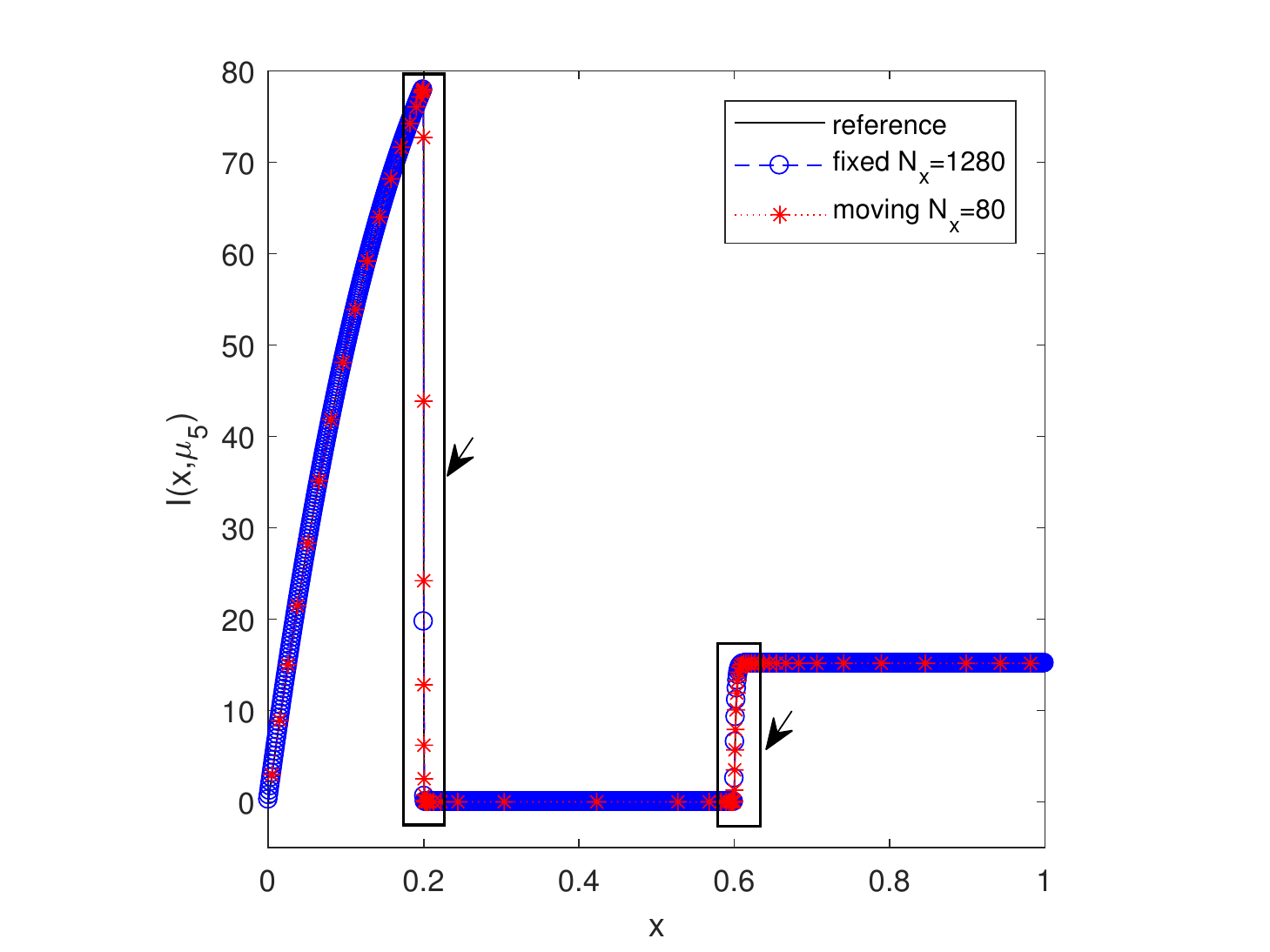}}
\subfigure[Close view of (c) near $x$=0.2 and $x$=0.6]{
\includegraphics[width=0.45\textwidth]{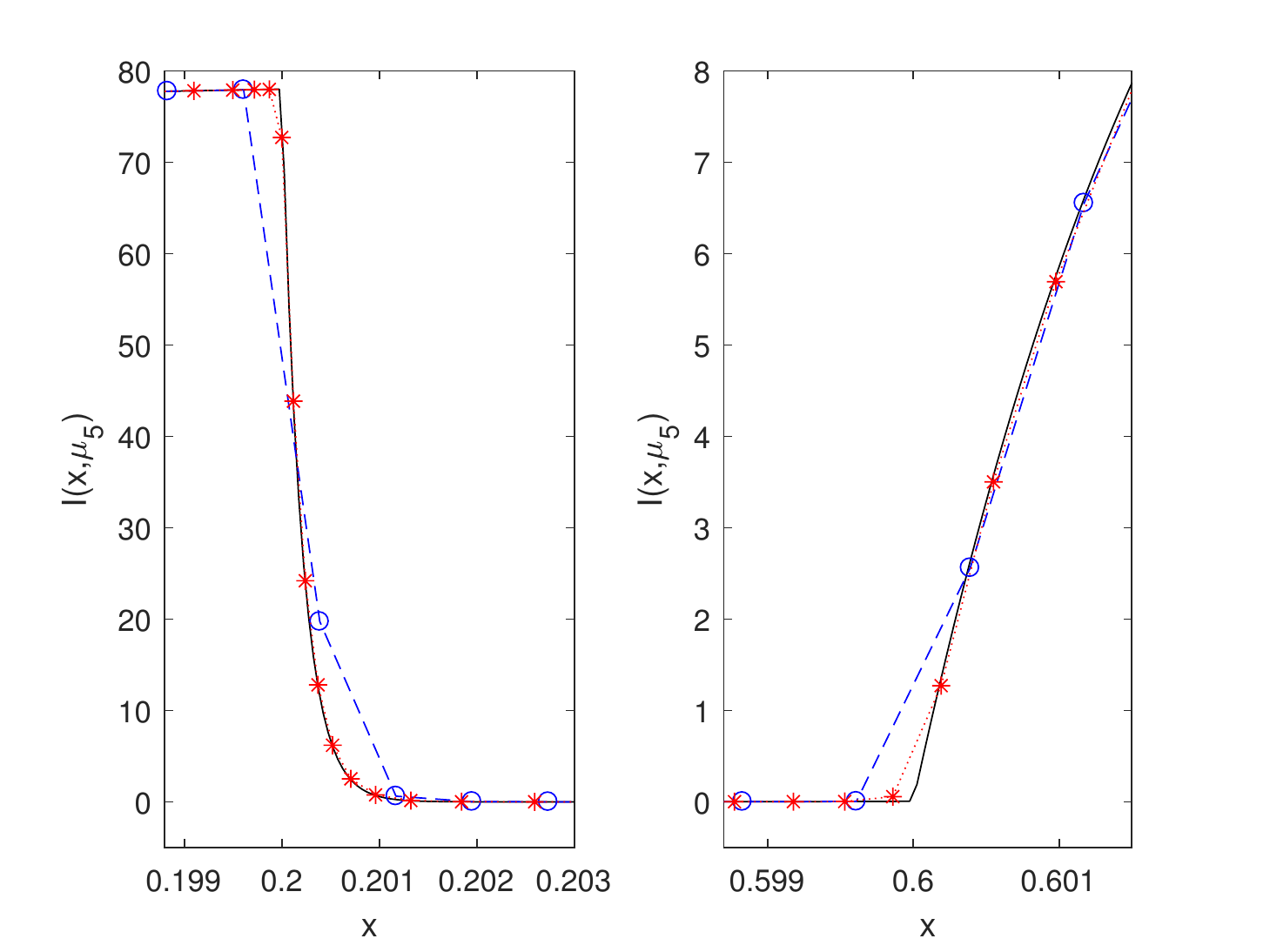}}
\caption{
\small{Example \ref{Ex6-1d}. The solution in the direction $\mu=0.1834$ obtained with the $P^2$-DG method
with a moving mesh of $N$=80 is compared with those obtained with fixed meshes of $N$=80 and $N$=1280.}
}\label{Fig:d1Ex6p2-u5}
\end{figure}

\noindent
In this example, we take $\sigma_s=1$, $c=3.0\times10^8$, and
\begin{equation*}
\sigma_t=
\begin{cases}
1,\quad   &\text{for }0\leq x< 0.2,\\
900,\quad &\text{for }0.2\leq x< 0.6,\\
90,\quad  &\text{for }0.6\leq x\leq1,\\
\end{cases}
\end{equation*}
 \begin{equation*}
q(x,\mu,t)=
\begin{cases}
100e^{-t},\quad&\text{for }0\leq x< 0.2,\\
1,\quad &\text{for }0.2\leq x<0.6,\\
1000e^{3t},\quad &\text{for }0.6\leq x\leq1.\\
\end{cases}
\end{equation*}
The initial condition is $I(x,\mu,0)=15x$
and the boundary conditions are given by
\begin{equation*}
\begin{split}
&I(0,\mu,t)=0, \quad \quad\quad~~ \text{for }0<\mu \leq 1,~~0< t\leq 0.1 ,
\\&I(1,\mu,t)=15+2t,\quad \text{for }-1 \leq\mu < 0,~~0< t\leq 0.1.
\end{split}
\end{equation*}
The solution of this problem has two sharp layers. Since its analytical form is unavailable
for comparison purpose we take the numerical solution obtained with the $P^2$-DG method with
a fixed mesh of $N=20000$ as the reference solution.
The mesh trajectories for the $P^2$-DG method with a moving mesh of $N = 80$ are shown
in Fig.~\ref{Fig:d1Ex6p2-mesh}. The moving mesh solution ($N=80$) in the direction $\mu=-0.1834$
is compared with the fixed mesh solutions obtained with $N=80$ and $N=1280$ in
Fig.~\ref{Fig:d1Ex6p2-u4}. Similar results are shown in Fig.~\ref{Fig:d1Ex6p2-u5} for the angular direction
$\mu = 0.1834$. The results show that the moving mesh solution
($N=80$) is more accurate than those with fixed meshes of $N=80$ and $N=1280$.

%

\begin{example}\label{Ex1-2d}
(An accuracy test for the two-dimensional unsteady RTE for the absorbing-scattering model.)
\end{example}

\begin{figure}[H]
\centering
\subfigure[$P^1$-DG]{
\includegraphics[width=0.45\textwidth]{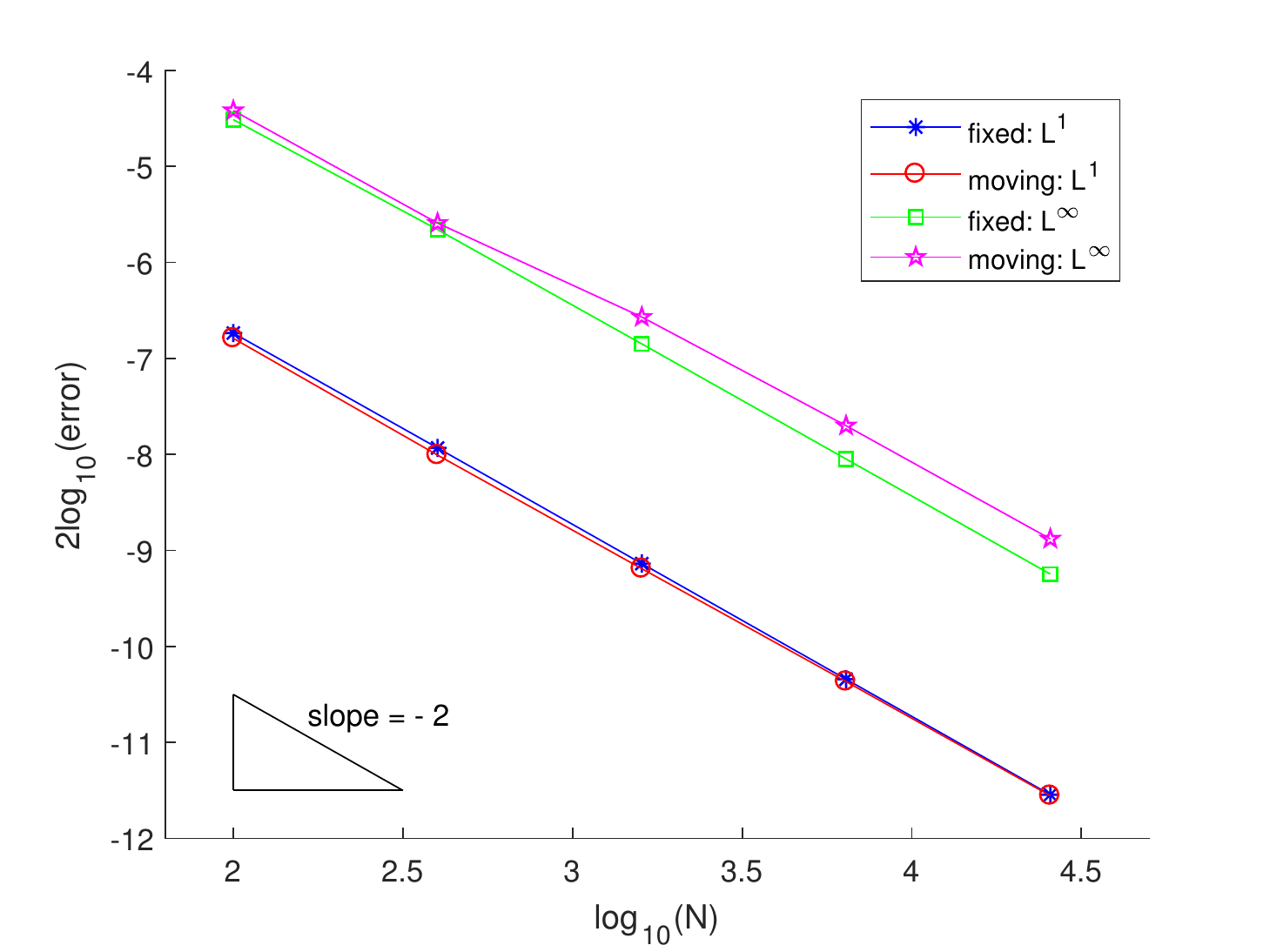}}
\subfigure[ $P^2$-DG]{
\includegraphics[width=0.45\textwidth]{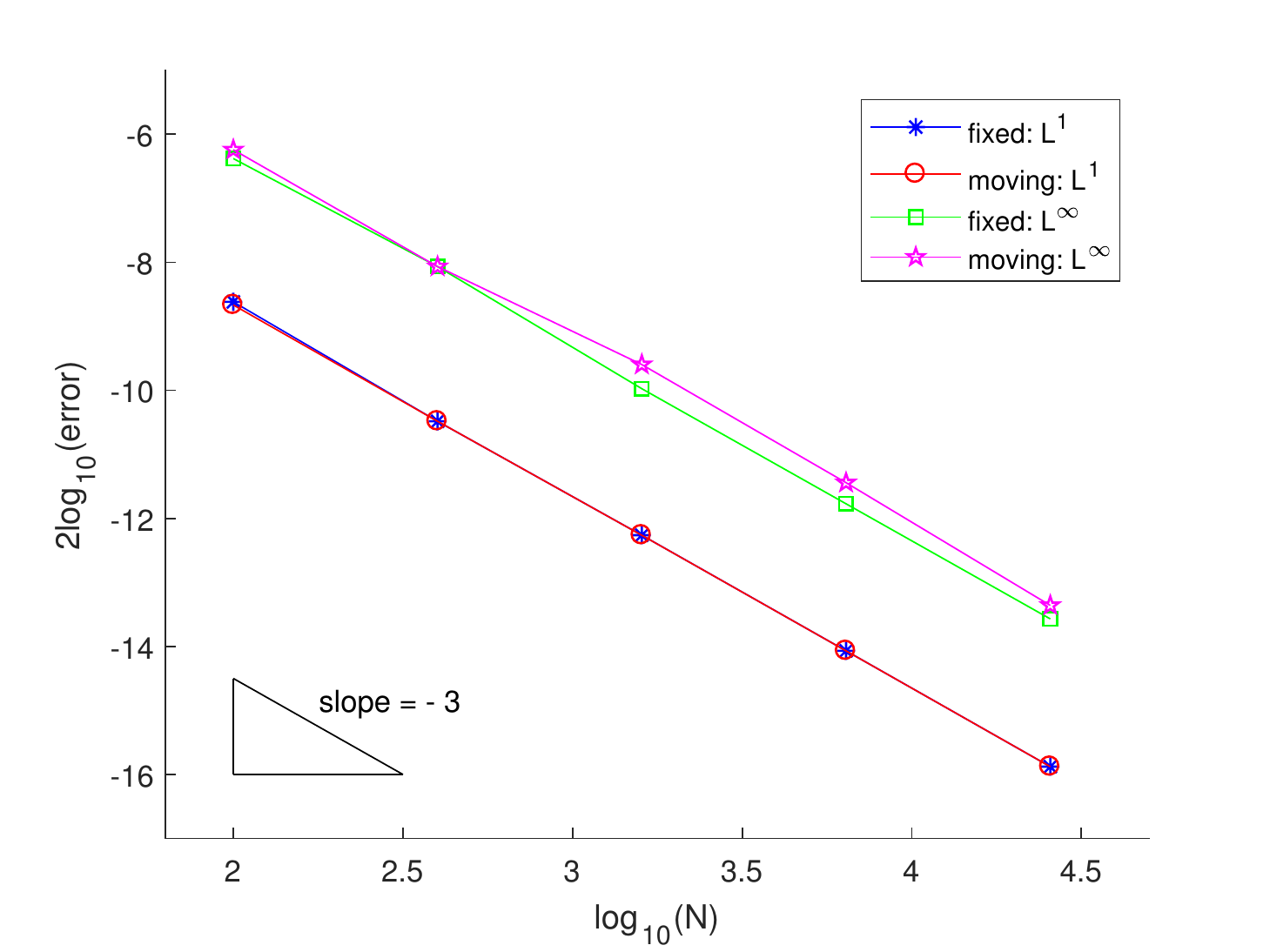}}
\caption{
\small{Example \ref{Ex1-2d}. The $L^1$ and $L^\infty$ norm of the error with a moving and fixed meshes.}
}\label{Fig:d2Ex1-order}
\end{figure}

\noindent
In this example, we take $\sigma_t=22000$, $\sigma_s=1$, $c=3.0\times10^{8}$, and
\begin{align*}
q(x,y,\zeta,\eta,t)= & e^{t}\Big{(}
-2\pi(\zeta+\eta)(\zeta^2+\eta^2)\cos^3\big{(}\frac{\pi}{2}(x+y)\big{)}
\sin\big{(}\frac{\pi}{2}(x+y)\big{)}
\\
& +\big{(}(\frac{1}{c}+\sigma_t)(\zeta^2+\eta^2)-\frac{2}{3}\sigma_s\big{)}
\cos^4\big{(}\frac{\pi}{2}(x+y)\big{)}
+(\frac{1}{c}+\sigma_t-\sigma_s)\Big{)}.
\end{align*}
The initial condition is $I(x,y,\zeta,\eta,0)=(\zeta^2+\eta^2)\cos^4(\frac{\pi}{2} (x+y))+1$
and the boundary conditions are given by
\begin{equation*}
\begin{split}
&I(x,0,\zeta,\eta,t)=e^{t}\big{(}(\zeta^2+\eta^2)\cos^4(\frac{\pi}{2}x)+1\big{)},
\quad\quad\quad~\eta>0 ,
\\
&I(x,1,\zeta,\eta,t)=e^{t}\big{(}(\zeta^2+\eta^2)\cos^4(\frac{\pi}{2}(x+1))+1\big{)},
\quad \eta<0 ,
\\
&I(0,y,\zeta,\eta,t)=e^{t}\big{(}(\zeta^2+\eta^2)\cos^4(\frac{\pi}{2}y)+1\big{)},
\quad\quad\quad~\zeta>0,
\\
&I(1,y,\zeta,\eta,t)=e^{t}\big{(}(\zeta^2+\eta^2)\cos^4(\frac{\pi}{2}(1+y))+1\big{)},
\quad \zeta<0.
\end{split}
\end{equation*}
The exact solution of this problem is
$I(x,y,\zeta,\eta,t)=e^{t}\big{(}(\zeta^2+\eta^2)\cos^4( \frac{\pi}{2}(x+y)) + 1\big{)}$.
The error in the $L^1$ and $L^\infty$ norm is plotted in Fig.~\ref{Fig:d2Ex1-order} for
the $P^1$-DG and $P^2$-DG methods with fixed and moving meshes. Once again, it can be seen that
both fixed and moving meshes lead to comparable results and the same convergence order (2nd for $P^1$-DG and
3rd for $P^2$-DG).

\begin{example}\label{Ex3-2d}
(A discontinuous example of the two-dimensional unsteady RTE for the transparent model.)
\end{example}

\noindent
In this test, we solve the two-dimensional unsteady RTE \eqref{s3.1} with $\sigma_t=0$,
$\sigma_s=0$, $c=3.0\times10^{8}$, $q = 0$, $\zeta = 0.3$, and $\eta = 0.5$.
The computational domain is $(0,1)\times (0,1)$. The initial condition is
\begin{equation*}
I(x,y,\zeta,\eta,0)=\left\{\begin{array}{ll}
0,\quad &\text{for }y < \frac{\eta}{\zeta}x,
\\
\cos^6\big{(}\frac{\pi}{2}y\big{)}, \quad & \text{otherwise}.\\
\end{array}
\right.
\end{equation*}
The boundary conditions are
\begin{equation*}
I(0,y,\zeta,\eta,t)=\cos^6\Big{(}\frac{\pi}{2}y\Big{)}\cos^{10}(t),\qquad
I(x,0,\zeta,\eta,t)=0.
\end{equation*}
The exact solution of this example is
\begin{equation*}
I(x,y,\zeta,\eta,t)=\left\{\begin{array}{ll}
0,\quad &\text{for }y < \frac{\eta}{\zeta}x,
\\
\cos^6\big{(}\frac{\pi}{2}(y-\frac{\eta}{\zeta}x )\big{)}\cos^{10}(t-\frac{x}{c\zeta} ), \quad &\text{otherwise}.
\end{array}
\right.
\end{equation*}
Notice that only a single angular direction is chosen in this example and the integral term in \eqref{s3.1}
is not involved.

\begin{figure}[H]
\centering
\subfigure[Radiative intensity on MM $N$=1600]{
\includegraphics[width=0.45\textwidth]{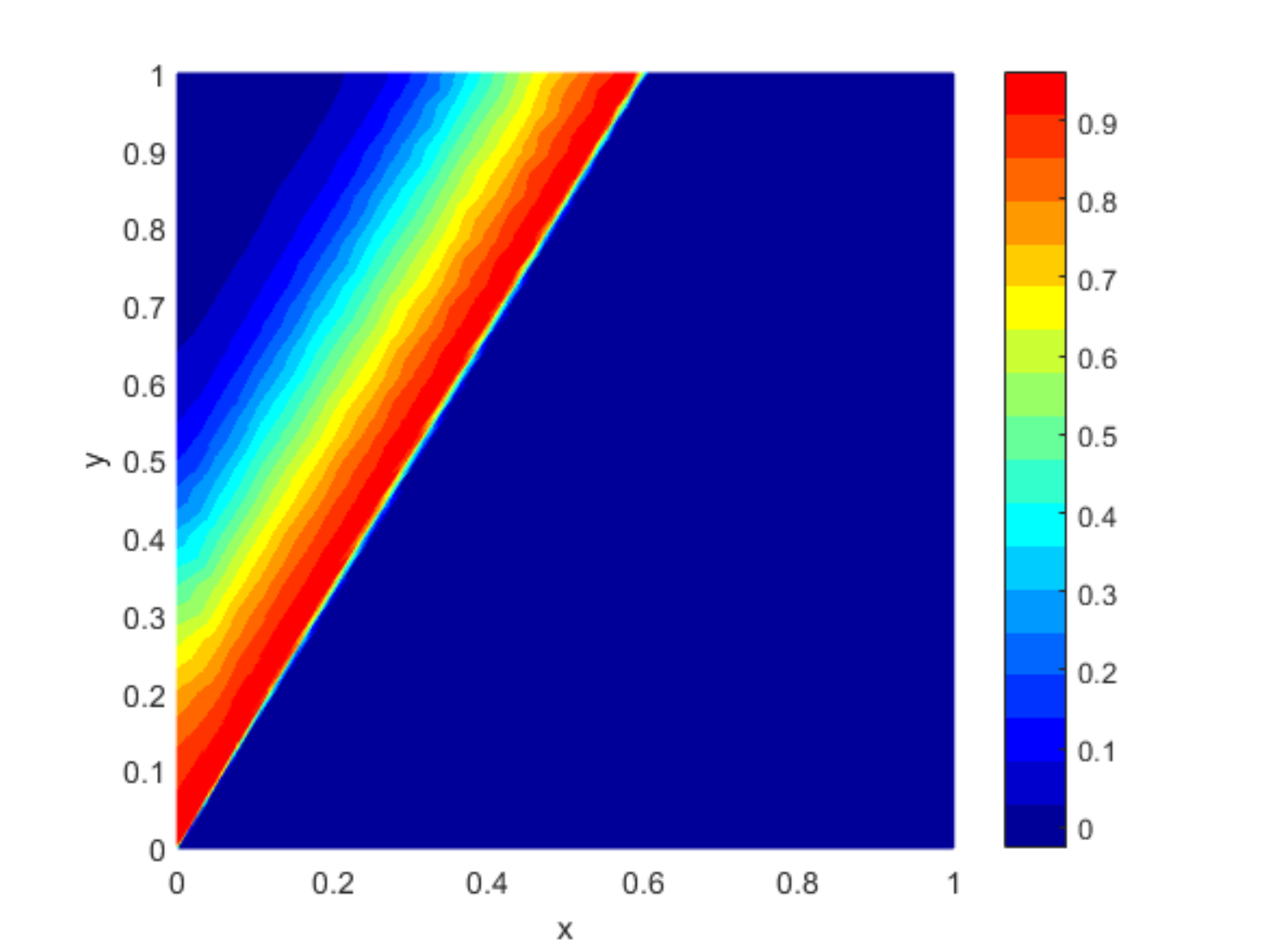}}
\subfigure[MM $N$=1600 at $t=0.1$]{
\includegraphics[width=0.45\textwidth]{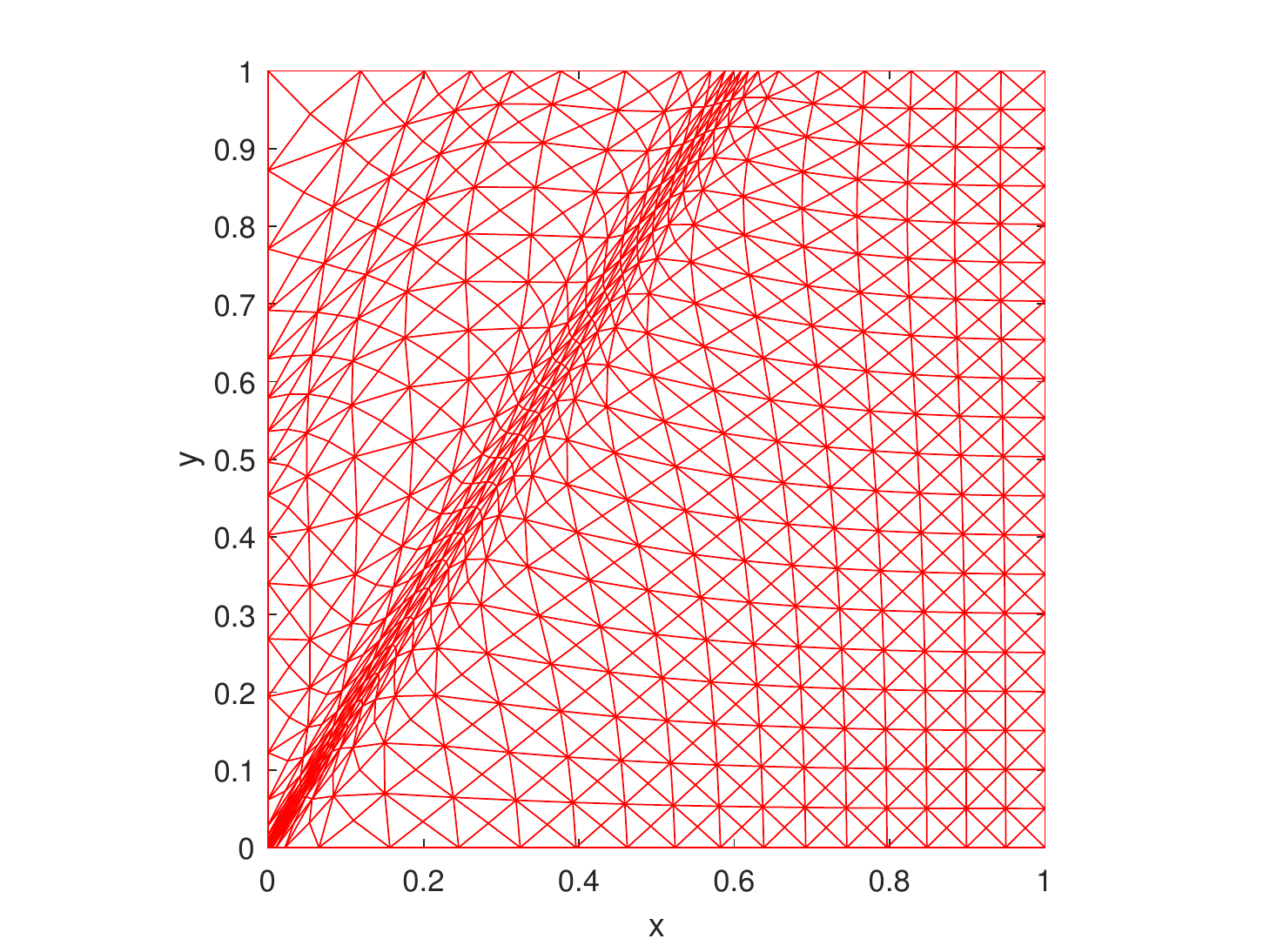}}
\subfigure[Radiative intensity on FM $N$=1600]{
\includegraphics[width=0.45\textwidth]{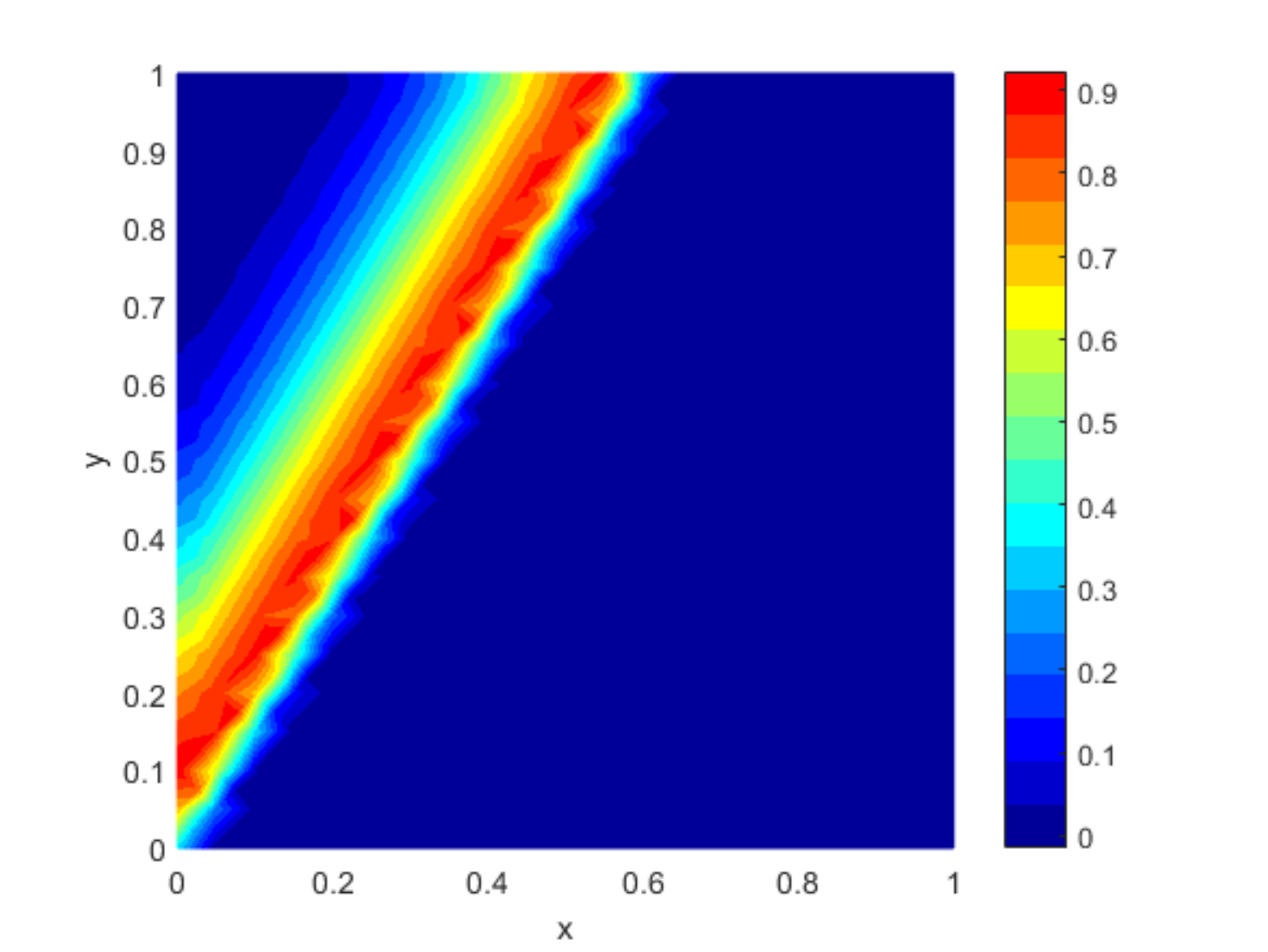}}
\subfigure[Radiative intensity on FM $N$=57600]{
\includegraphics[width=0.45\textwidth]{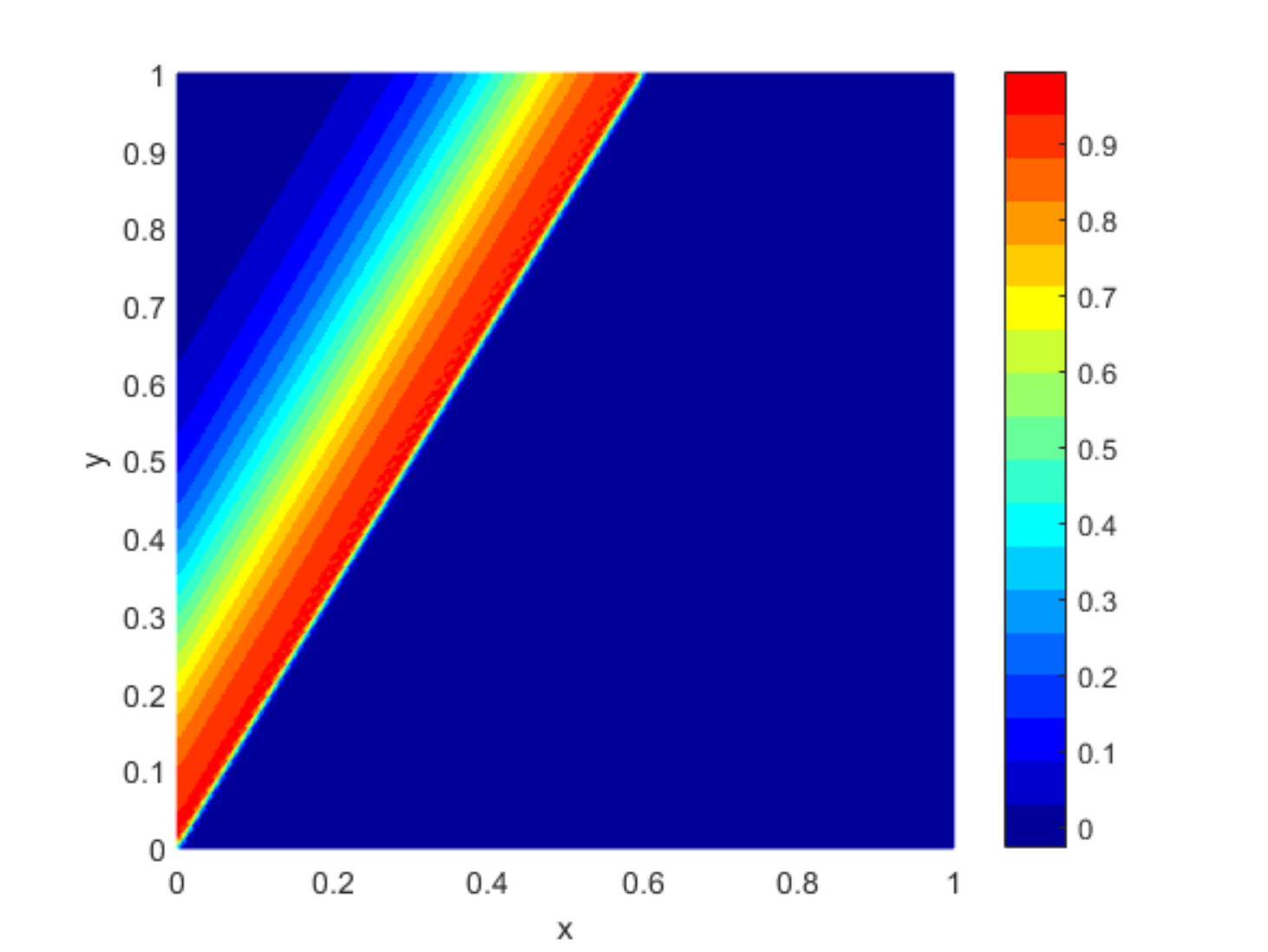}}
\caption{\small{Example \ref{Ex3-2d}. The radiative intensity contours (and mesh) at $t=0.1$
obtained by the $P^2$-DG method
with fixed and moving meshes.}}
\label{Fig:d2Ex3p2-1}
\end{figure}

\begin{figure}[H]
\centering
\subfigure[MM $N$=1600, FM $N$=1600]{
\includegraphics[width=0.45\textwidth]{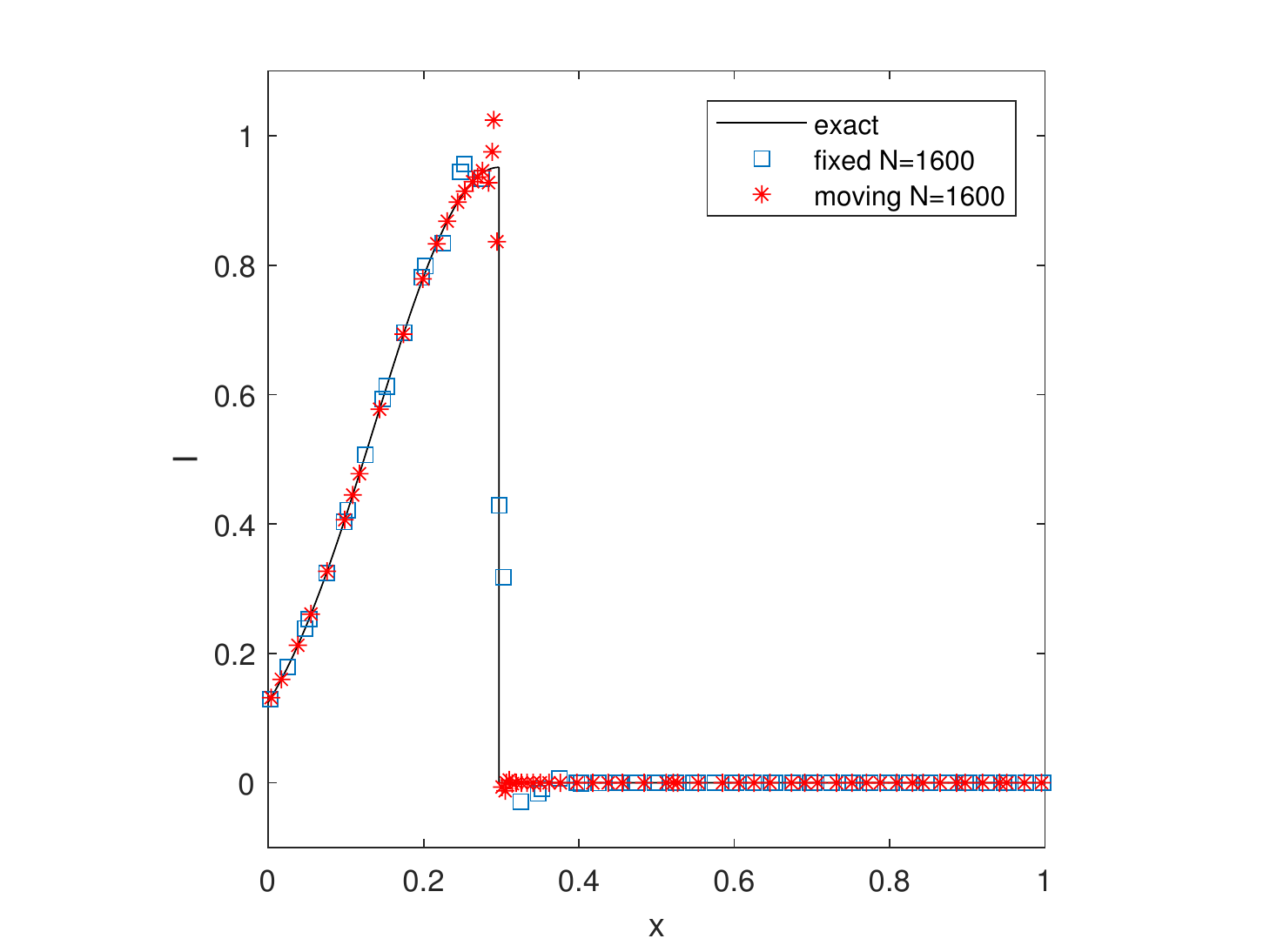}}
\subfigure[MM $N$=1600, FM $N$=57600]{
\includegraphics[width=0.45\textwidth]{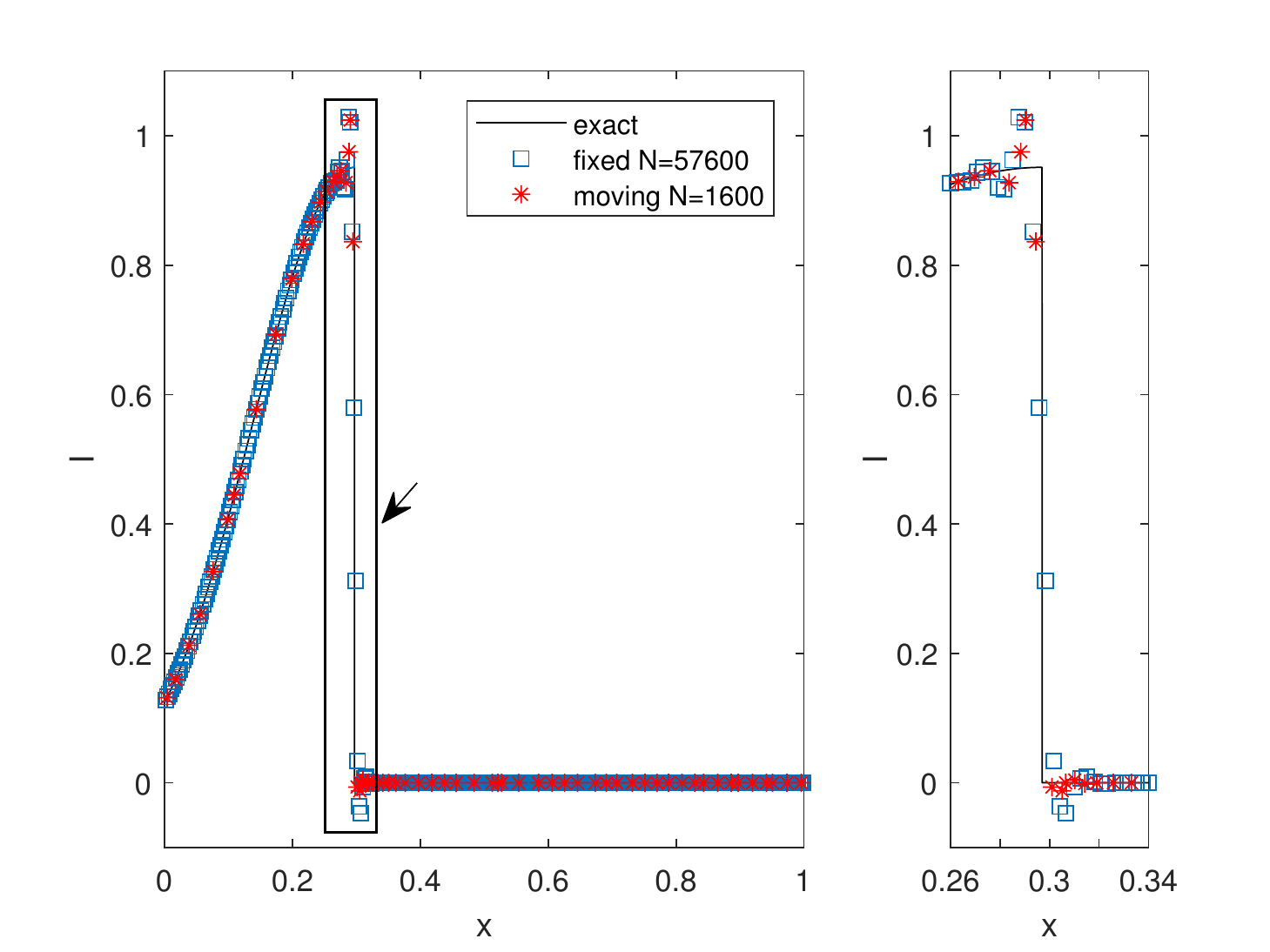}}
\caption{\small{Example \ref{Ex3-2d}. The comparison of the radiative intensity cut along the line $y=0.495$ obtained
by the $P^2$-DG method with a moving mesh of $N=1600$ and fixed meshes of $N=1600$ and $N=57600$.}}
\label{Fig:d2Ex3p2-2}
\end{figure}

The radiative intensity contours obtained with a moving mesh of $N=1600$ and fixed meshes of
$N=1600$ and $N=57600$ are shown in Fig.~\ref{Fig:d2Ex3p2-1}.
In Fig.~\ref{Fig:d2Ex3p2-2}, we compare the radiative intensity cut along the line $y=0.495$
for the moving mesh of $N=1600$ and the fixed meshes of $N=1600$ and $N=57600$. The results show that the moving mesh solution ($N=1600$) is better than that with the fixed mesh of $N=1600$ and is comparable with that with the fixed mesh of $N=57600$.

\begin{figure}[H]
\centering
\subfigure[$P^1$-DG]{
\includegraphics[width=0.45\textwidth]{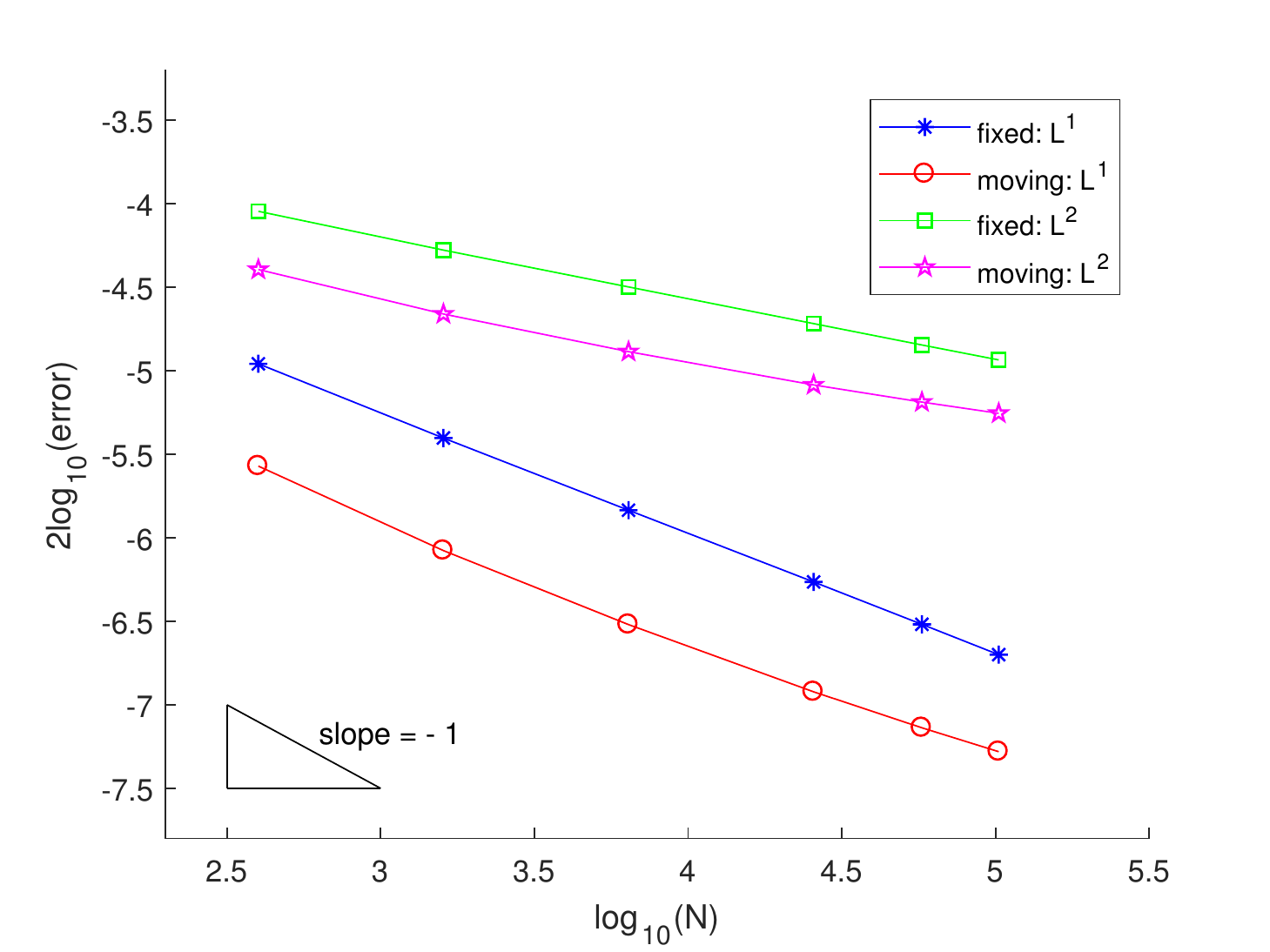}}
\subfigure[ $P^2$-DG]{
\includegraphics[width=0.45\textwidth]{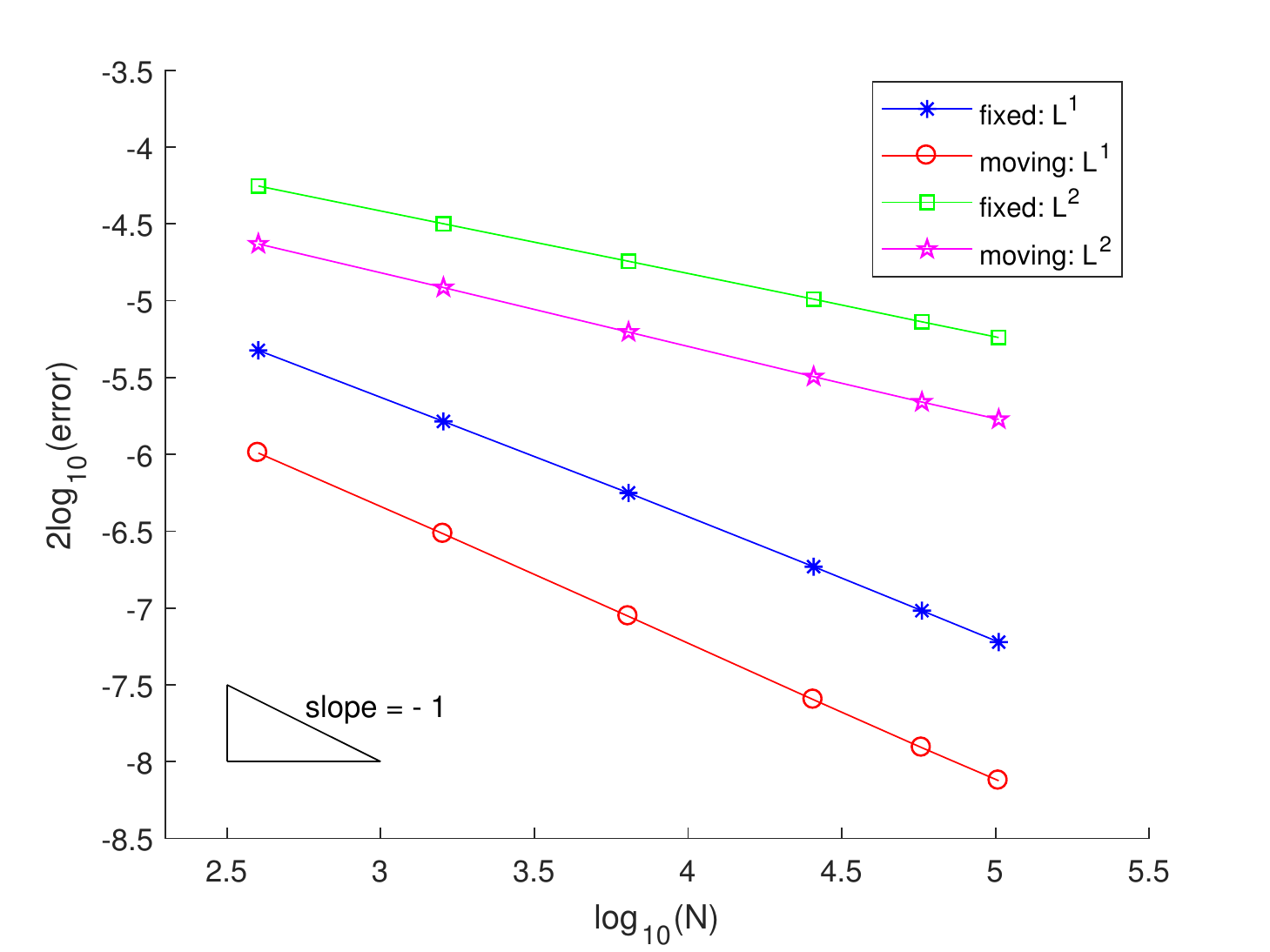}}
\caption{
\small{Example \ref{Ex3-2d}. The $L^1$ and $L^2$ norm of the error with a moving and fixed meshes.}
}\label{Fig:d2Ex3-order}
\end{figure}

The error in the $L^1$ and $L^2$ norm is shown in Fig.~\ref{Fig:d2Ex3-order} for the $P^1$-DG and $P^2$-DG
methods with fixed and moving meshes. It is worth pointing out that we cannot expect $P^1$-DG and $P^2$-DG
can achieve their optimal order for this problem since the solution is discontinuous.
One can see from the figure that both fixed and moving meshes lead to almost the same
convergence order. The order of $P^1$-DG is about 0.72 in $L^1$ norm and 0.36 in $L^2$ norm while
that of $P^2$-DG is about 0.81 in $L^1$ norm and 0.41 in $L^2$ norm.
Moreover, the figure shows that a moving mesh produces more accurate solutions than a fixed mesh
of the same number of elements for this example.

\begin{example}\label{Ex2-2d}
(A discontinuous example of the two-dimensional unsteady RTE for the purely absorbing model.)
\end{example}

\noindent
In this example, we choose $\sigma_t=1$, $\sigma_s=0$, $c=3.0\times10^{8}$, $q = 0$, $\zeta = 0.4$, and $\eta = 0.9$. The computational domain is $(0,1)\times (0,1)$. The initial condition is
\begin{equation*}
I(x,y,\zeta,\eta,0)=\begin{cases}
\tanh(500(x-0.5))+1,\quad & \text{for }y < \frac{\eta}{\zeta}x,\\
1, \quad &\text{otherwise}.\\
\end{cases}
\end{equation*}
The boundary conditions are
\begin{equation*}
I(0,y,\zeta,\eta,t)= e^{(y^2t)},\qquad
I(x,0,\zeta,\eta,t)=\tanh(500(x-0.5))+1.
\end{equation*}
The exact solution of this example is
\begin{equation*}
I(x,y,\zeta,\eta,t)=\begin{cases}
(\tanh(500(x-\frac{\zeta}{\eta}y-0.5))+1)e^{-\frac{\sigma_t}{\eta}y},\quad &\text{for } y < \frac{\eta}{\zeta}x,\\
e^{\big{(}(y-\frac{\eta}{\zeta}x)^2(t-\frac{x}{c\zeta})-\frac{\sigma_t}{\zeta}x\big{)}}, \quad &\text{otherwise},\\
\end{cases}
\end{equation*}
which exhibits a discontinuity along with $y=\frac{\eta}{\zeta}x$ and a sharp layer along with
$x=\frac{\zeta}{\eta}y+0.5$. Like the previous example, only a single angular direction is chosen in this example
and the integral term in \eqref{s3.1} is not involved.

\begin{figure}[H]
\centering
\subfigure[Radiative intensity on MM $N$=6400]{
\includegraphics[width=0.45\textwidth]{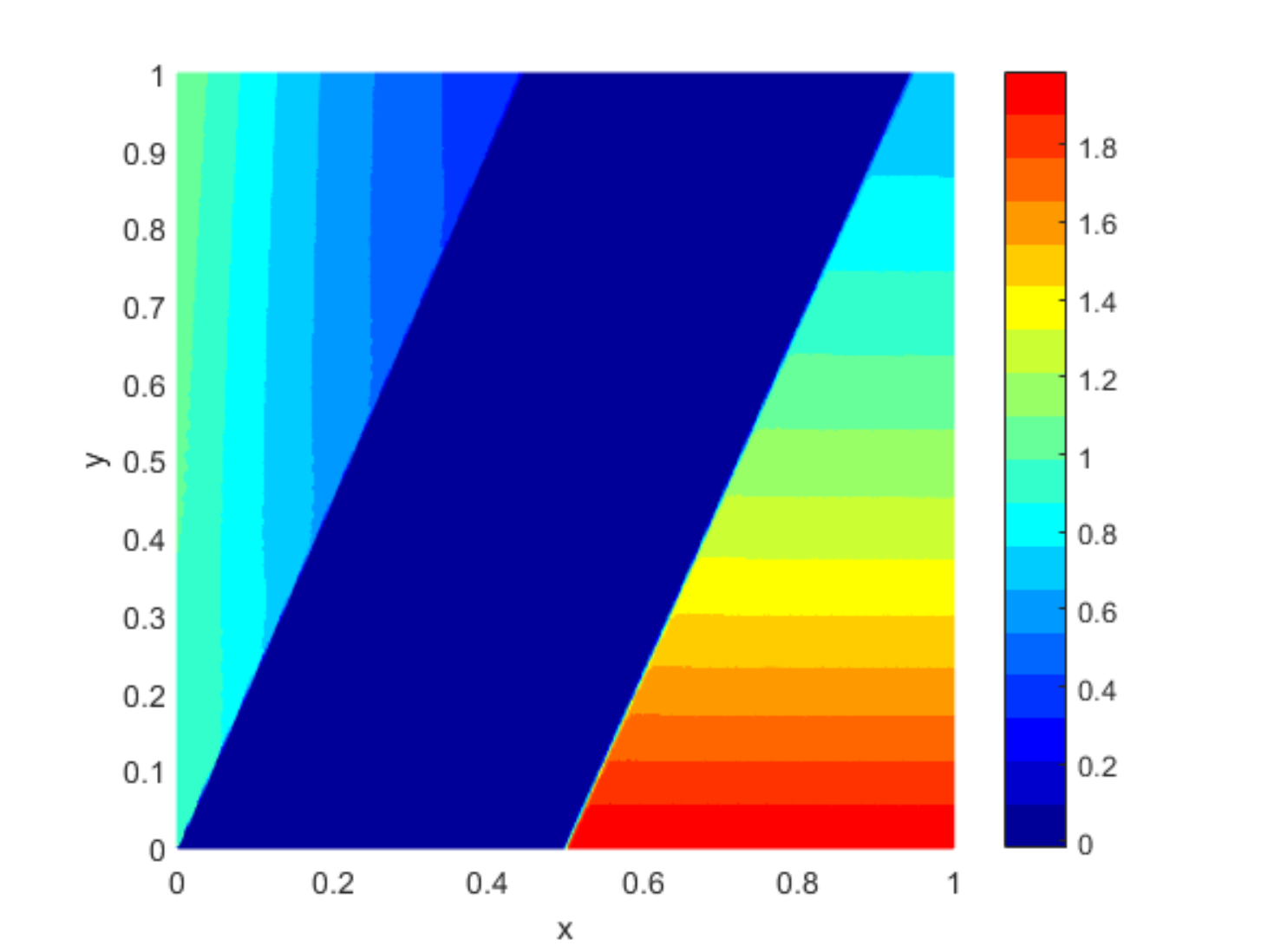}}
\subfigure[MM $N$=6400 at $t=0.1$]{
\includegraphics[width=0.45\textwidth]{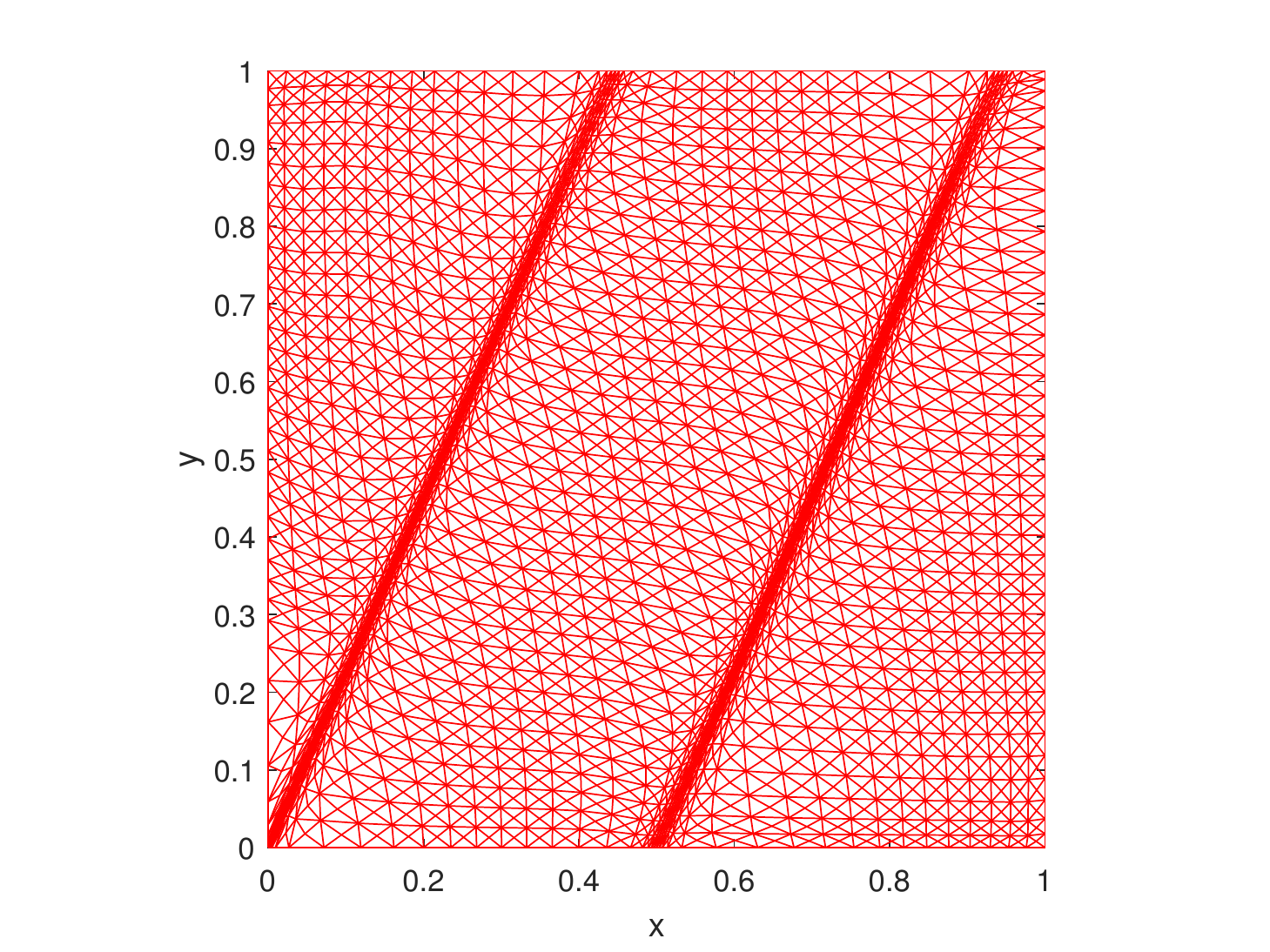}}
\subfigure[Radiative intensity on FM $N$=6400]{
\includegraphics[width=0.45\textwidth]{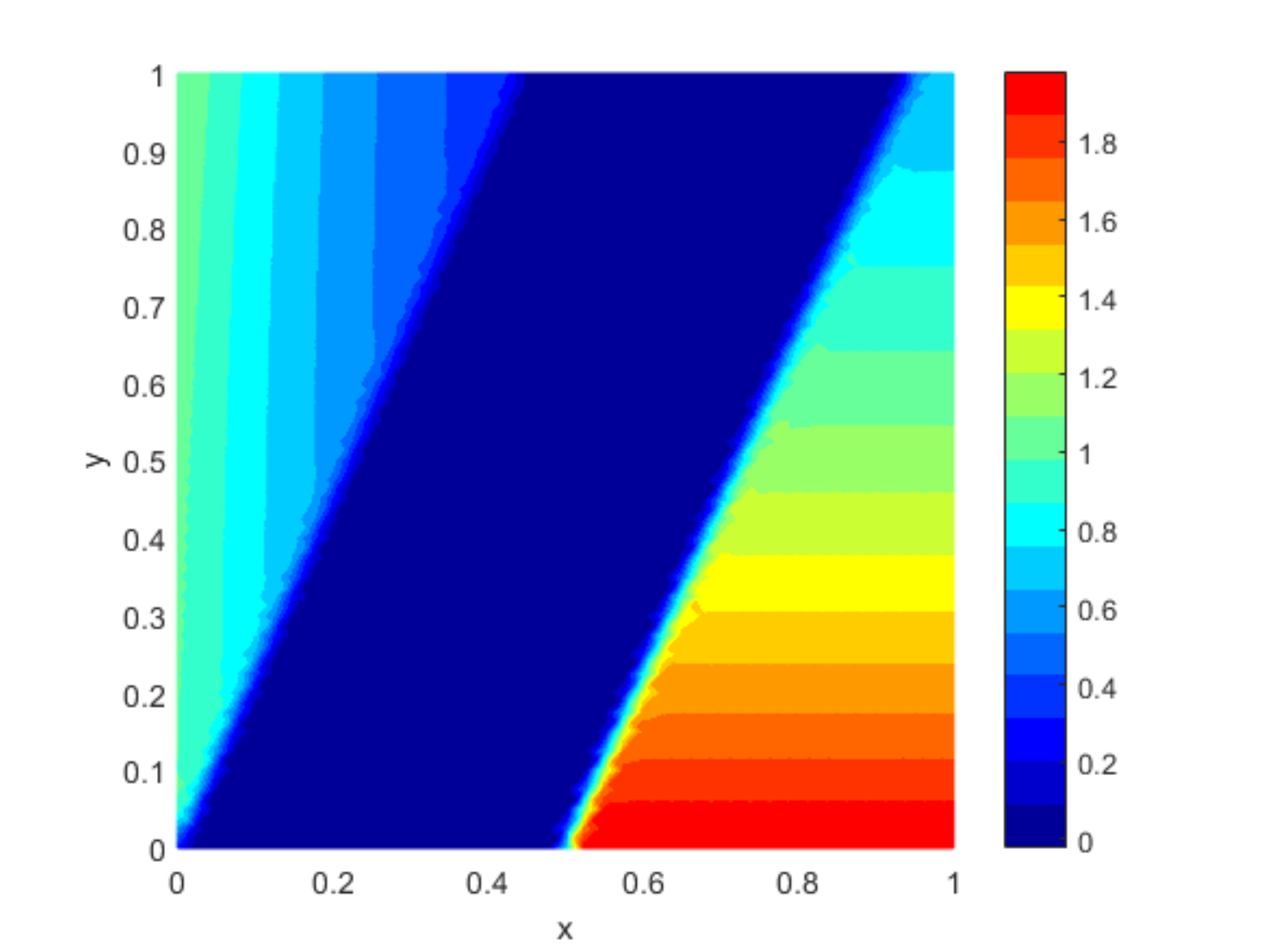}}
\subfigure[Radiative intensity on FM $N$=102400]{
\includegraphics[width=0.45\textwidth]{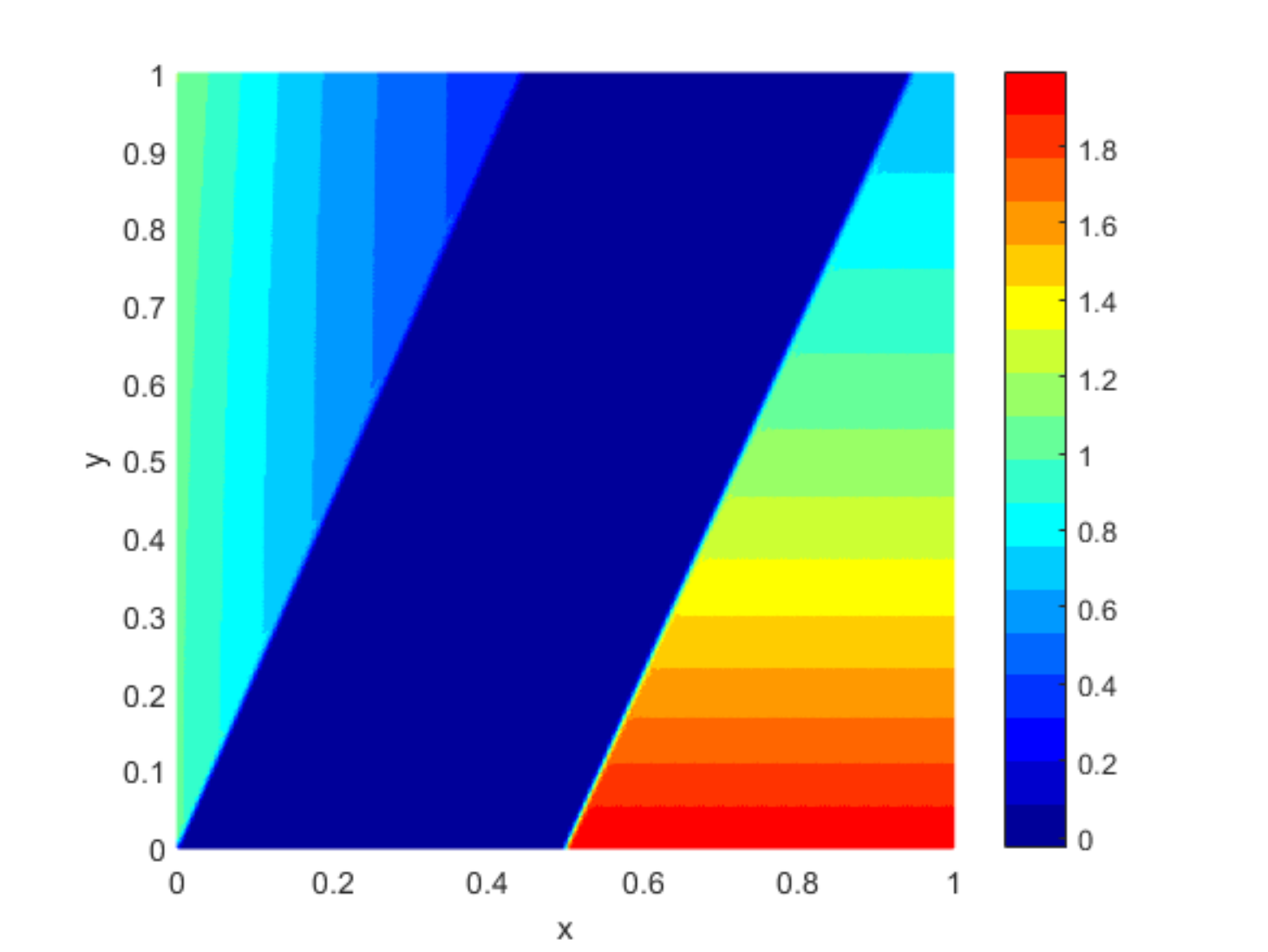}}
\caption{
\small{Example \ref{Ex2-2d}. The radiative intensity contours (and mesh) at $t = 0.1$ obtained by
the $P^2$-DG method with moving and fixed meshes.}}
\label{Fig:d2Ex2p2-1}
\end{figure}

\begin{figure}[H]
\centering
\subfigure[MM $N$=6400, FM $N$=6400]{
\includegraphics[width=0.45\textwidth]{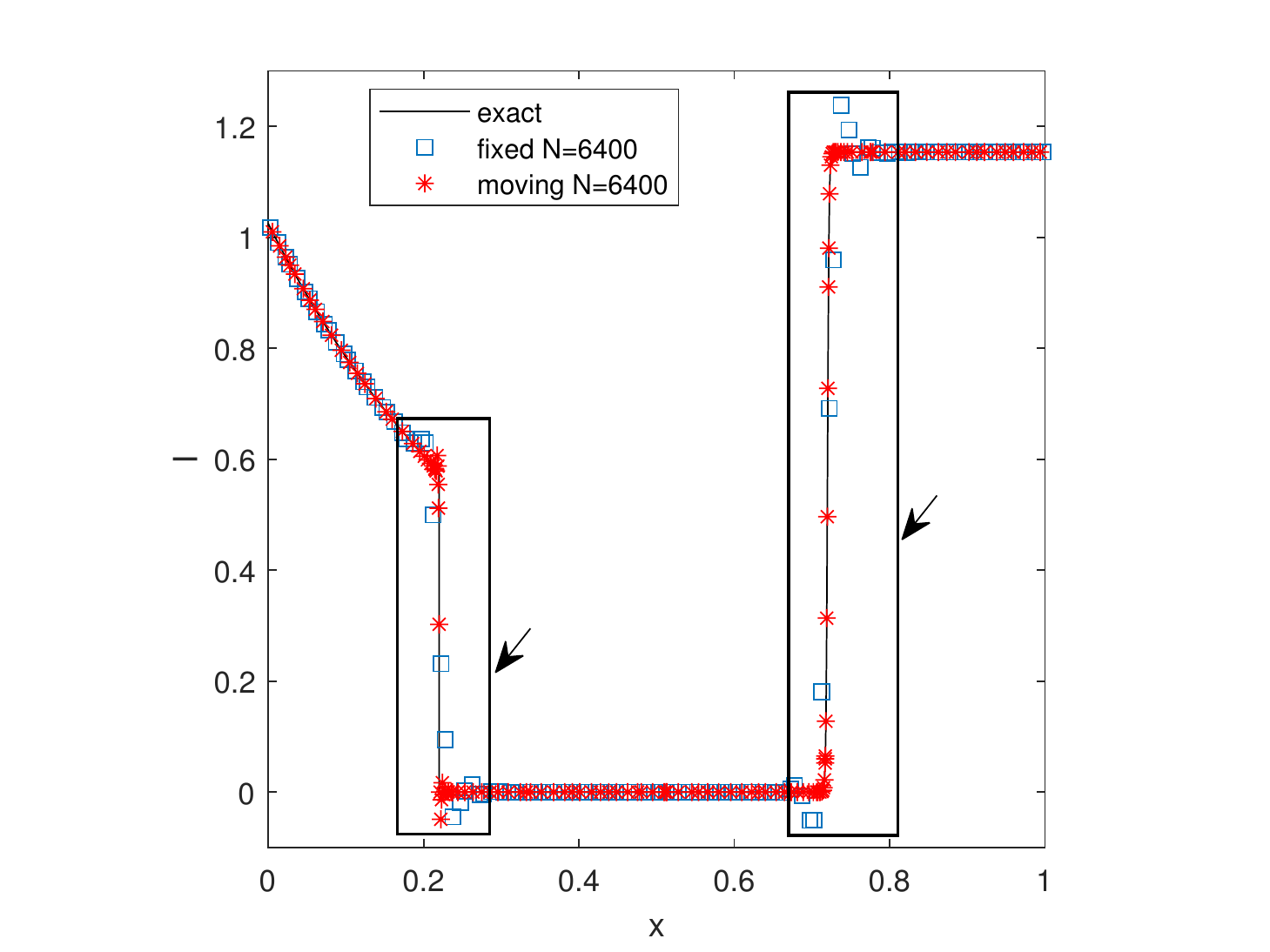}}
\subfigure[Close view of (a) near $x$=0.22 and $x$=0.72]{
\includegraphics[width=0.45\textwidth]{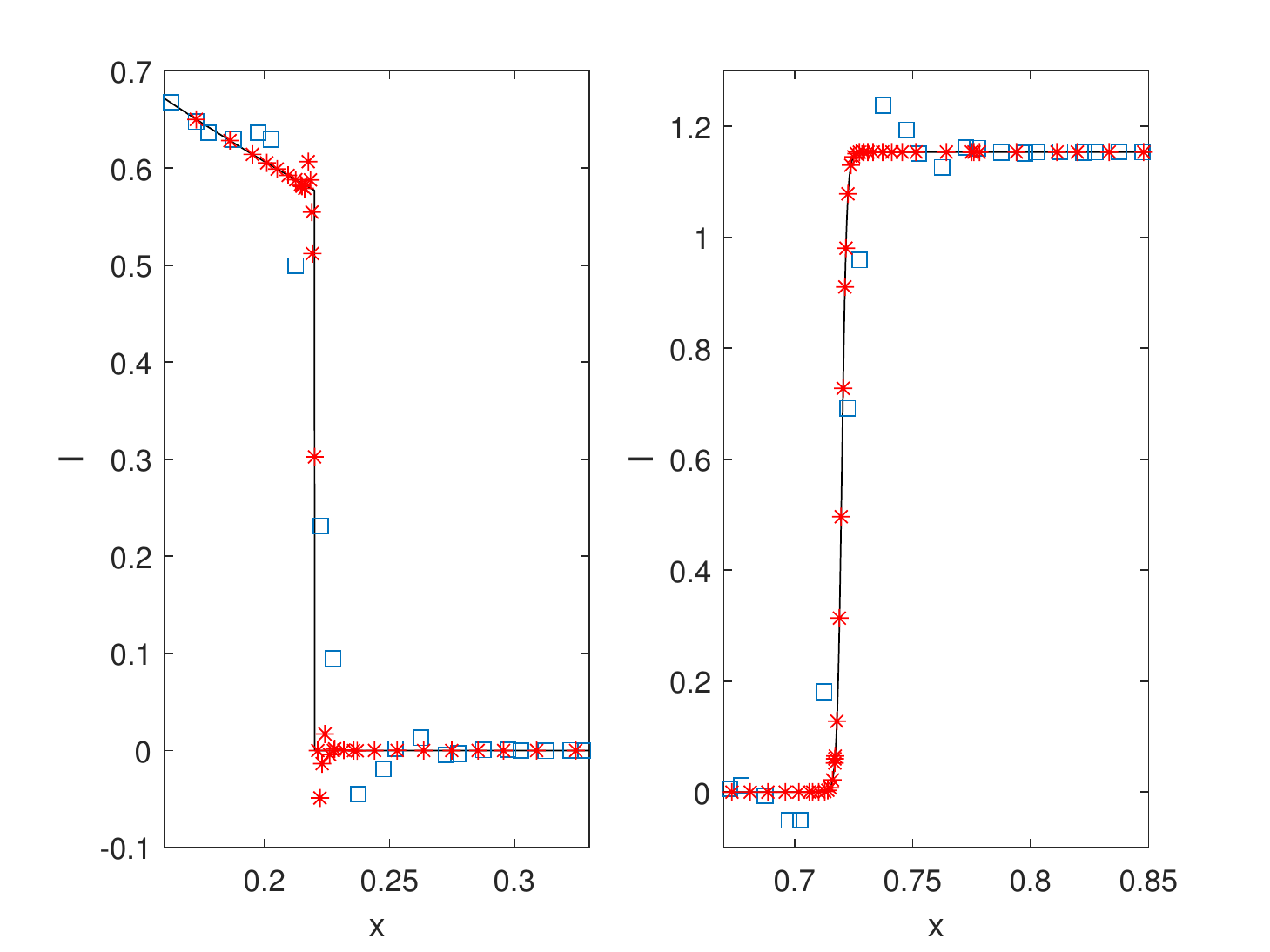}}
\subfigure[MM $N$=6400, FM $N$=102400]{
\includegraphics[width=0.45\textwidth]{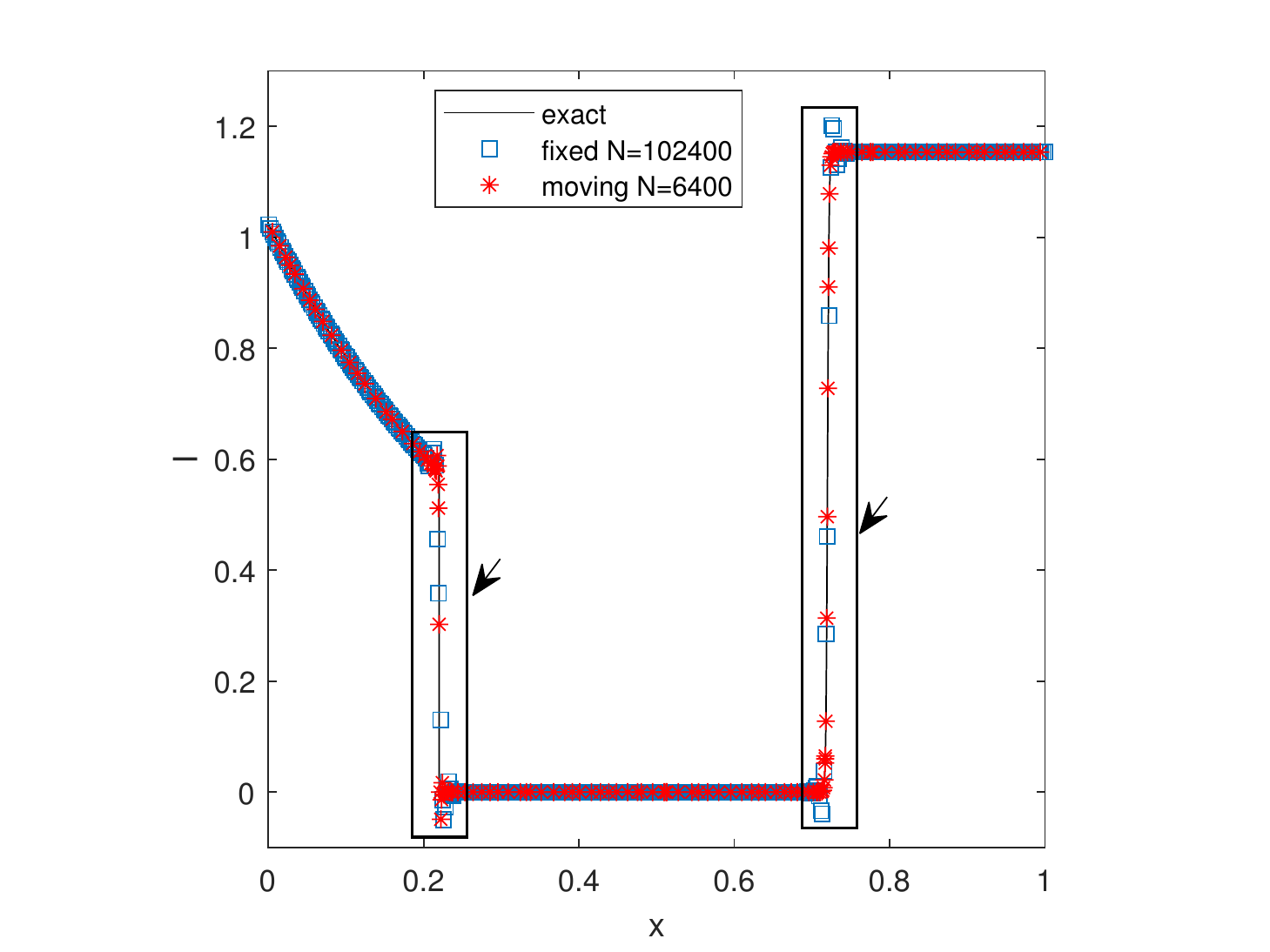}}
\subfigure[Close view of (c) near $x$=0.22 and $x$=0.72]{
\includegraphics[width=0.45\textwidth]{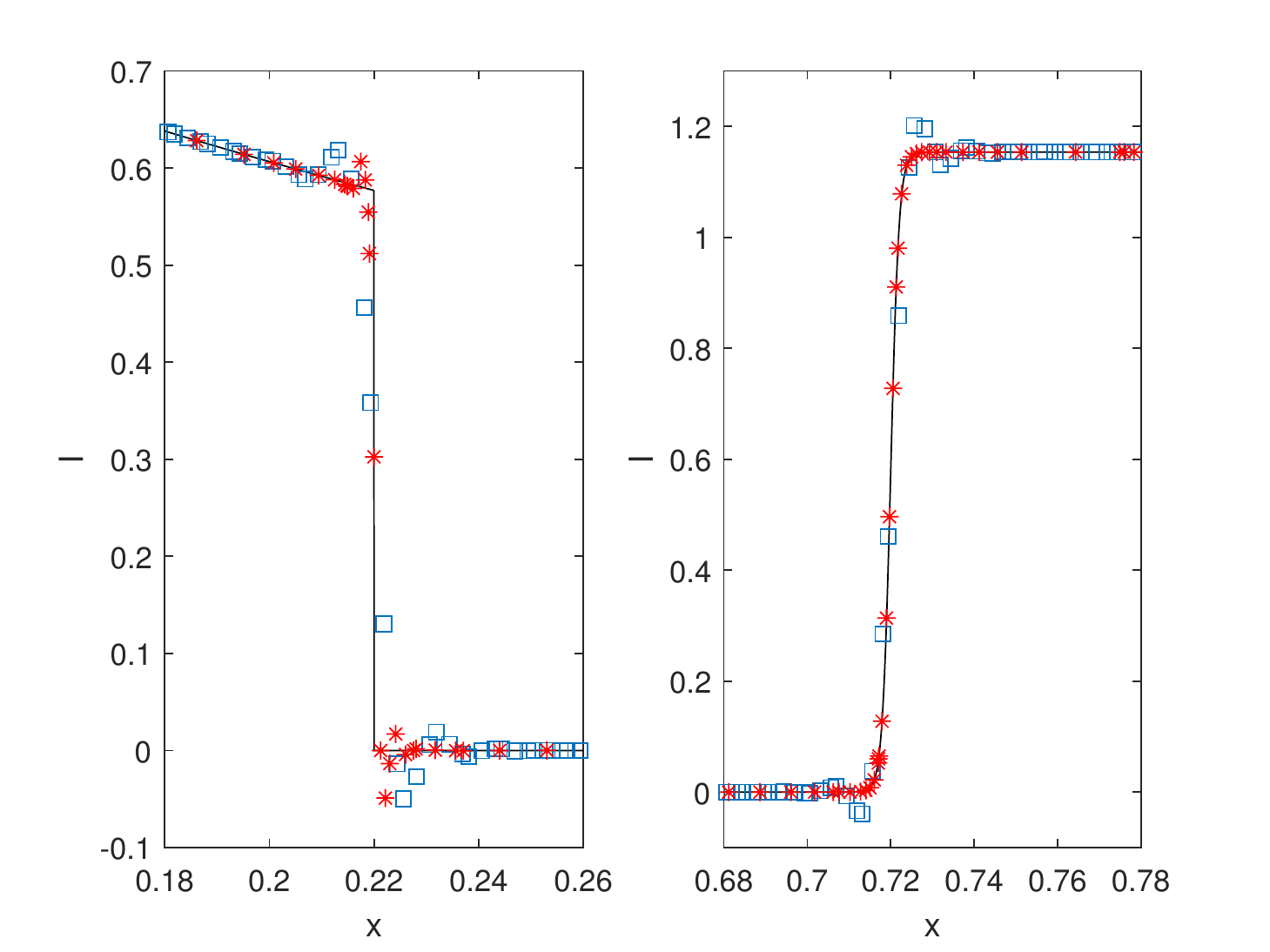}}
\caption{\small{Example \ref{Ex2-2d}. The radiative intensity cut along the line $y=0.495$ obtained a moving mesh
of $N=6400$ is compared with those obtained with fixed meshes of $N=6400$ and $N=102400$.}}
\label{Fig:d2Ex2p2-2}
\end{figure}

\begin{figure}[H]
\centering
\subfigure[$P^1$-DG]{
\includegraphics[width=0.45\textwidth]{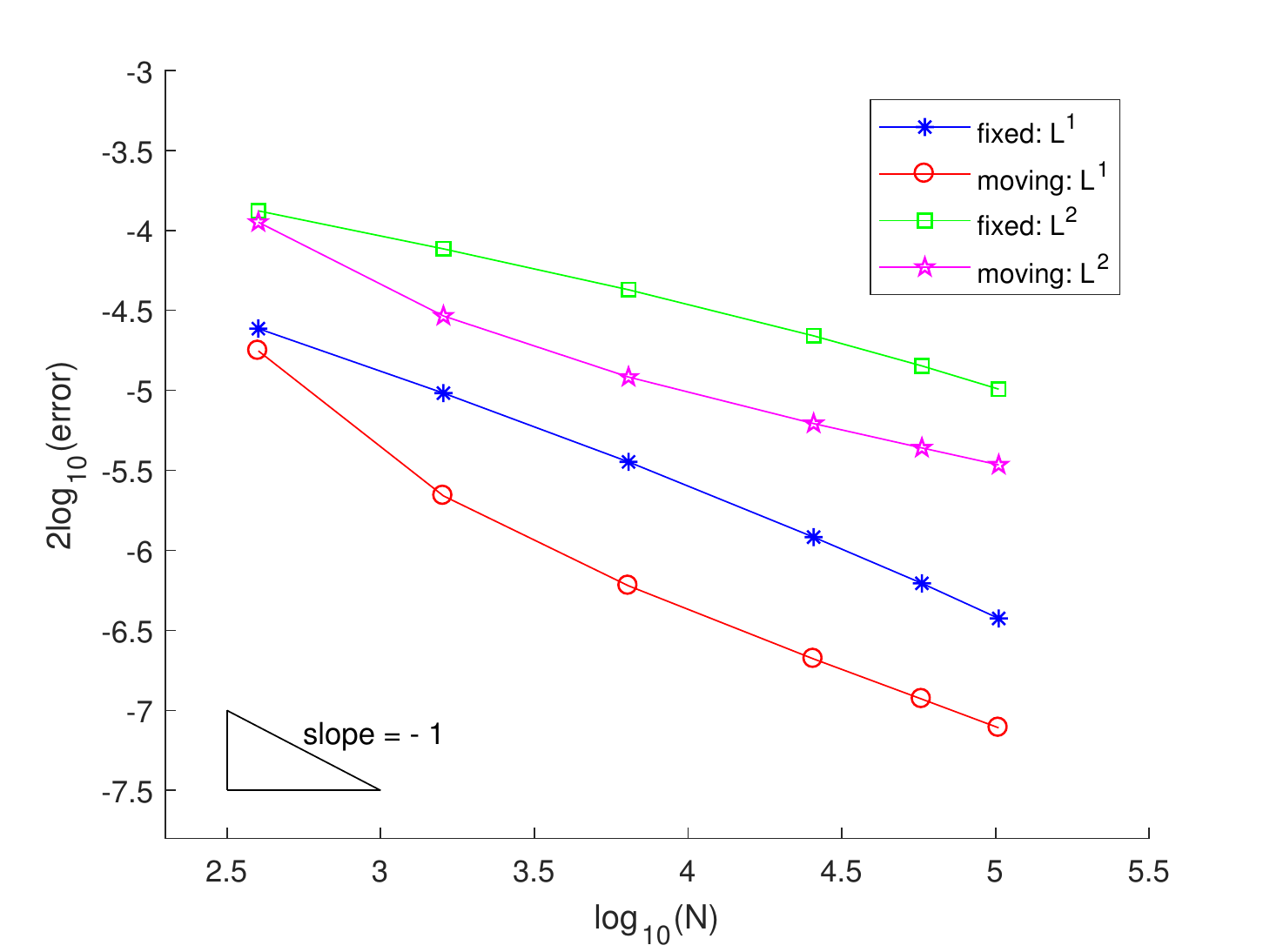}}
\subfigure[ $P^2$-DG]{
\includegraphics[width=0.45\textwidth]{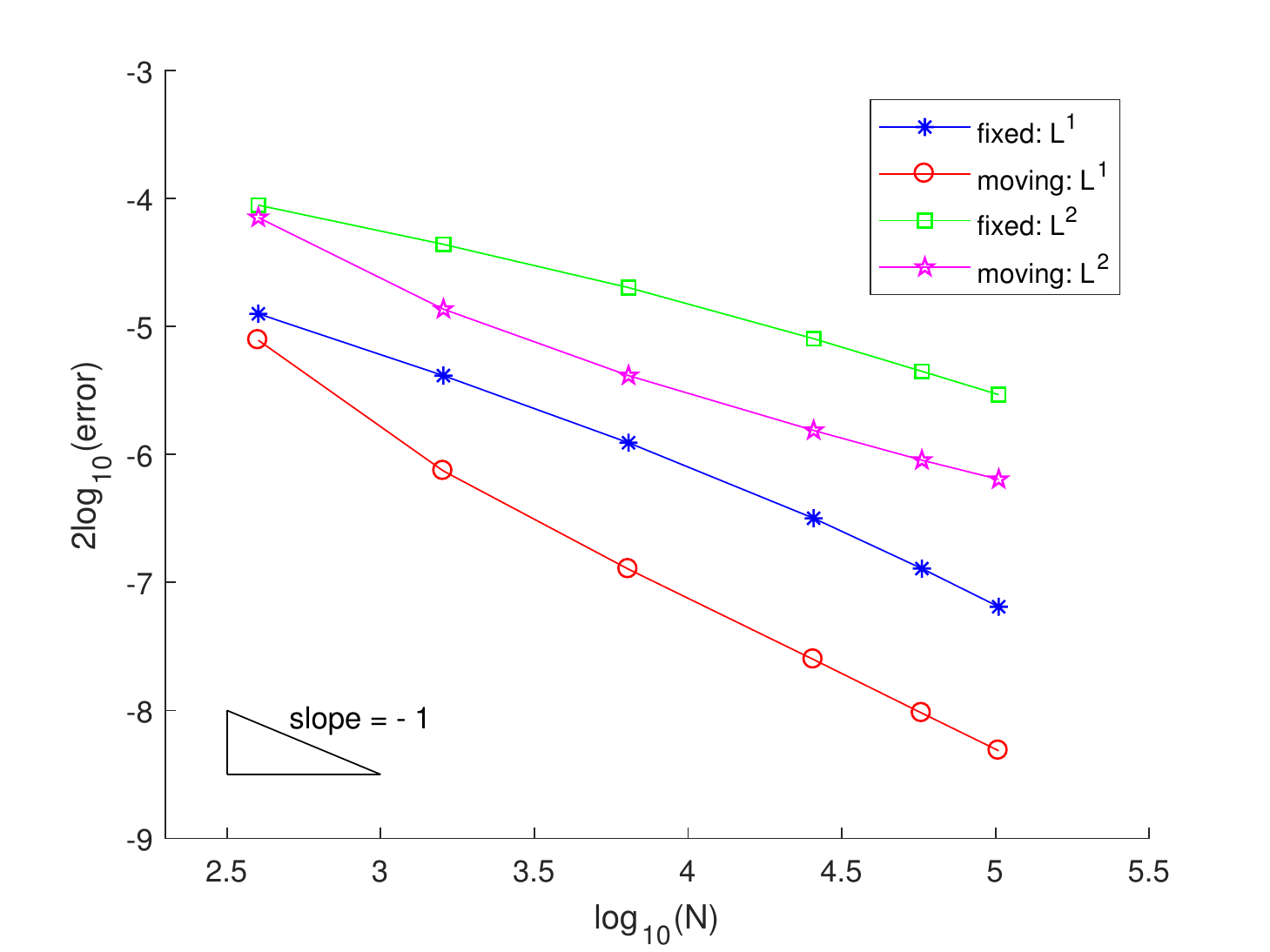}}
\caption{
\small{Example \ref{Ex2-2d}. The $L^1$ and $L^2$ norm of the error with a moving and fixed meshes.}
}\label{Fig:d2Ex2-order}
\end{figure}

The radiative intensity contours obtained with the $P^2$-DG method with a moving mesh of $N = 6400$
and fixed meshes of $N = 6400$ and $N = 102400$ are shown in Fig.~\ref{Fig:d2Ex2p2-1}.
In Fig.~\ref{Fig:d2Ex2p2-2}, the radiative intensity cut along the line $y=0.495$ is compared for
moving and fixed meshes. The advantage of using a moving mesh is clear.
The error in the $L^1$ abd $L^2$ norm is plotted as a function of $N$ in Fig.~\ref{Fig:d2Ex2-order}
for fixed and moving meshes. The $P^1$-DG method shows an order of about 0.87 in $L^1$ norm
and 0.57 in $L^2$ norm while $P^2$-DG has an order of 1.2 in $L^1$ norm
and 0.73 in $L^2$ norm for both fixed and moving meshes. A moving mesh produces more accurate
solutions than a fixed mesh of the same size for this example.

\begin{figure}[H]
\centering
\includegraphics[width=0.45\textwidth]{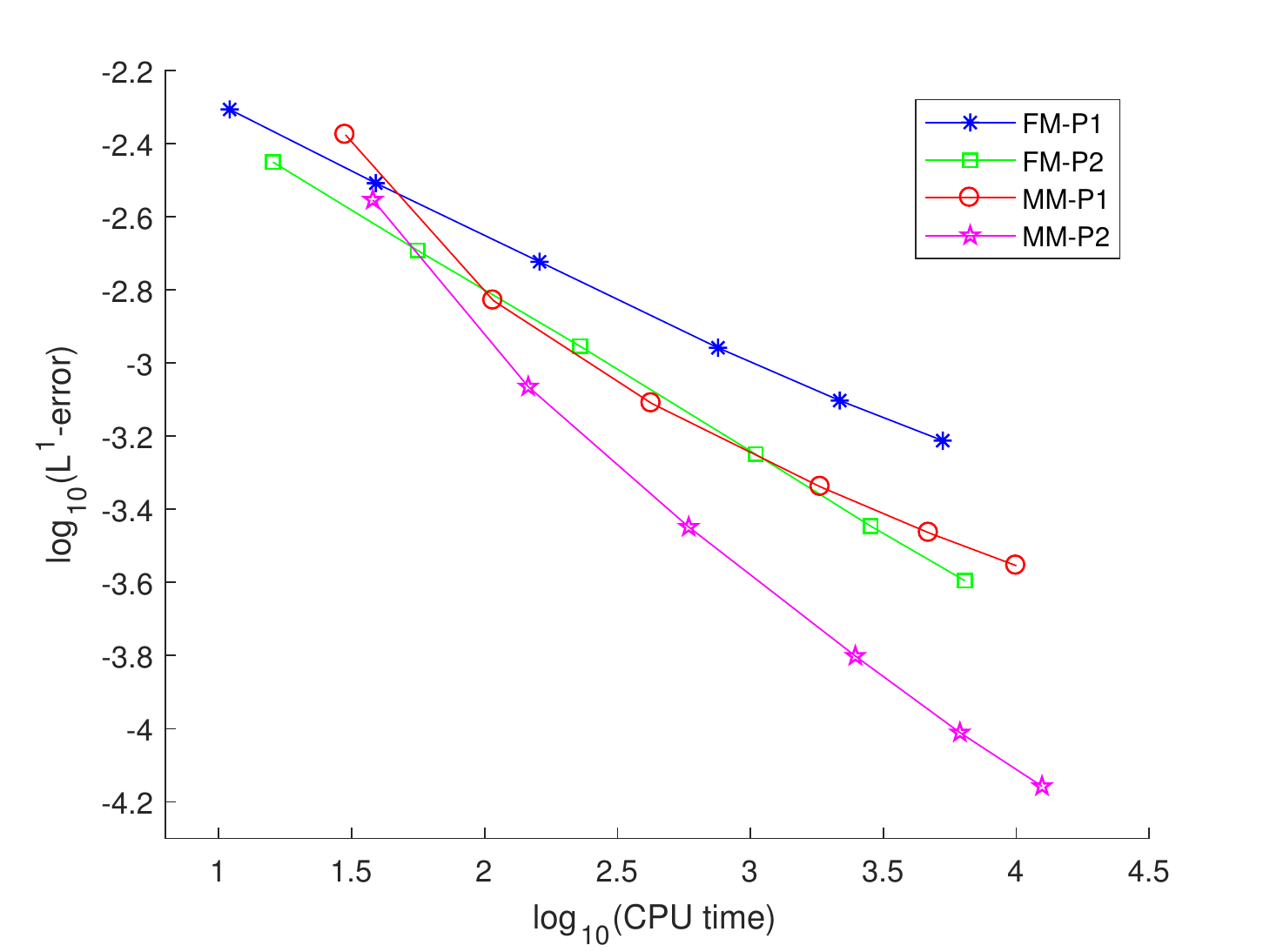}
\caption{
\small{Example \ref{Ex2-2d}. The $L^1$ norm of the error is plotted against the CPU time.}
}\label{Fig:d2Ex2-cpuL1}
\end{figure}

To show the efficiency of the methods for this two-dimensional example, we plot
the $L^1$ norm of the error against the CPU time in Fig.~\ref{Fig:d2Ex2-cpuL1}.
We can see that the error is smaller for moving mesh $P^1$-DG (resp. $P^2$-DG)
than fixed mesh $P^1$-DG (resp. $P^2$-DG) for a fixed amount of the CPU time.
Moreover, the better efficiency of a higher-order method is more obvious in this example
than the one-dimensional example~\ref{Ex7-1d}:
Fixed mesh $P^2$-DG is nearly equally or more efficient than moving mesh $P^1$-DG.


\begin{example}\label{Ex4-2d}
(A steep transition layer example of the two-dimensional unsteady RTE for the absorbing-scattering model.)
\end{example}

\begin{figure}[H]
\centering
\subfigure[Radiative intensity on MM $N$=6400]{
\includegraphics[width=0.45\textwidth]{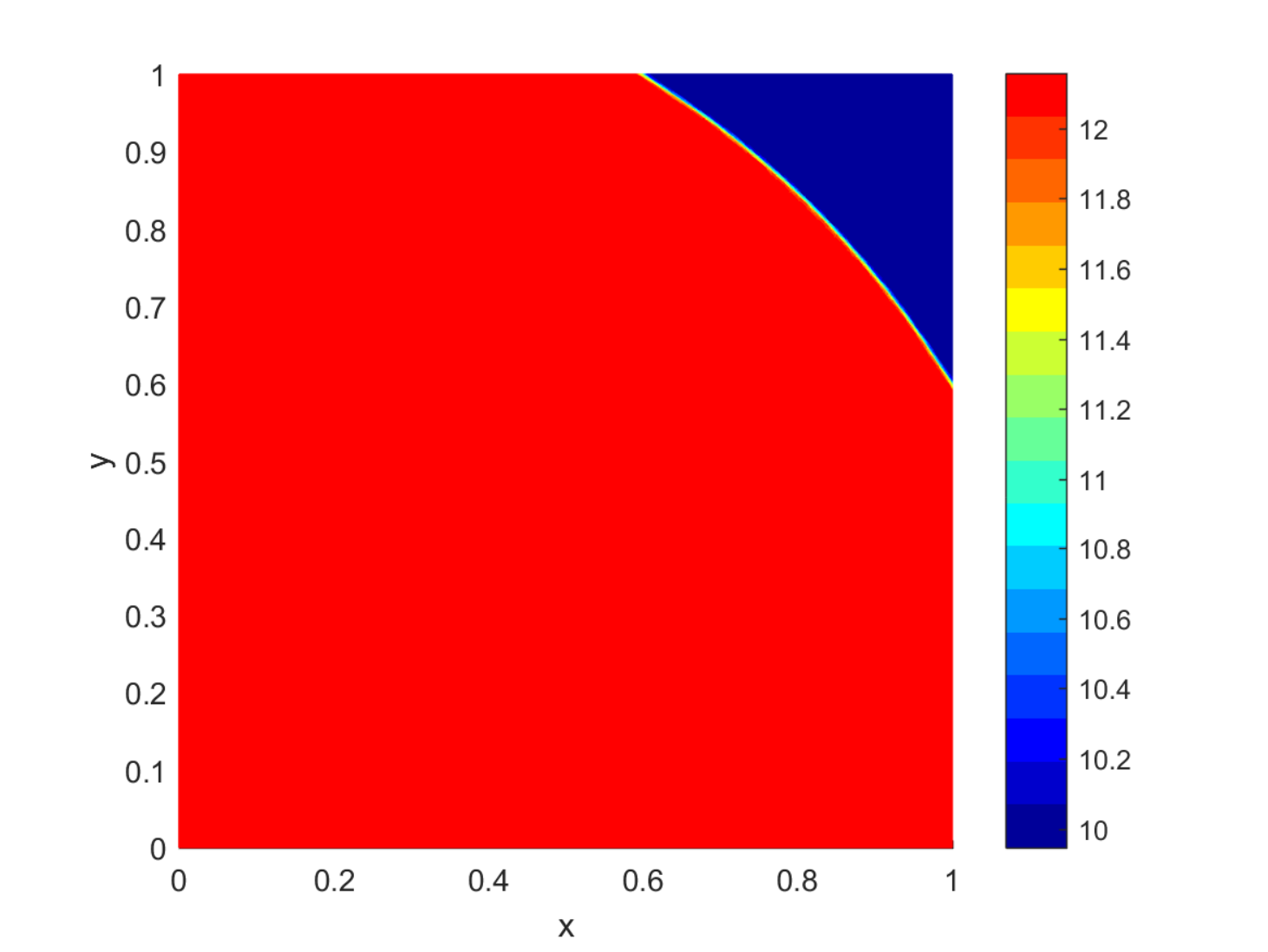}}
\subfigure[MM $N$ = 6400 at $t=0.1$]{
\includegraphics[width=0.45\textwidth]{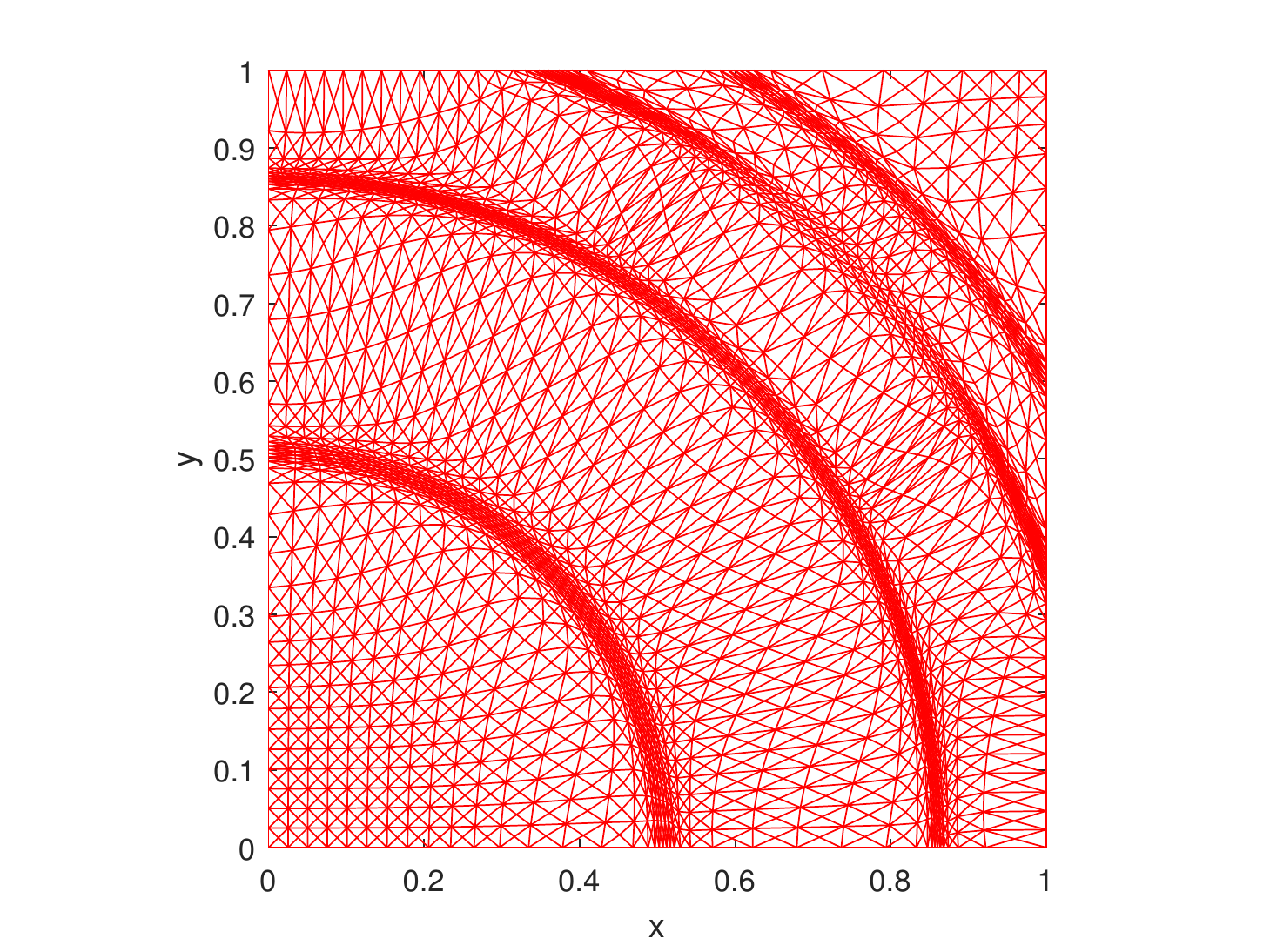}}
\subfigure[ Radiative intensity on FM $N$=6400]{
\includegraphics[width=0.45\textwidth]{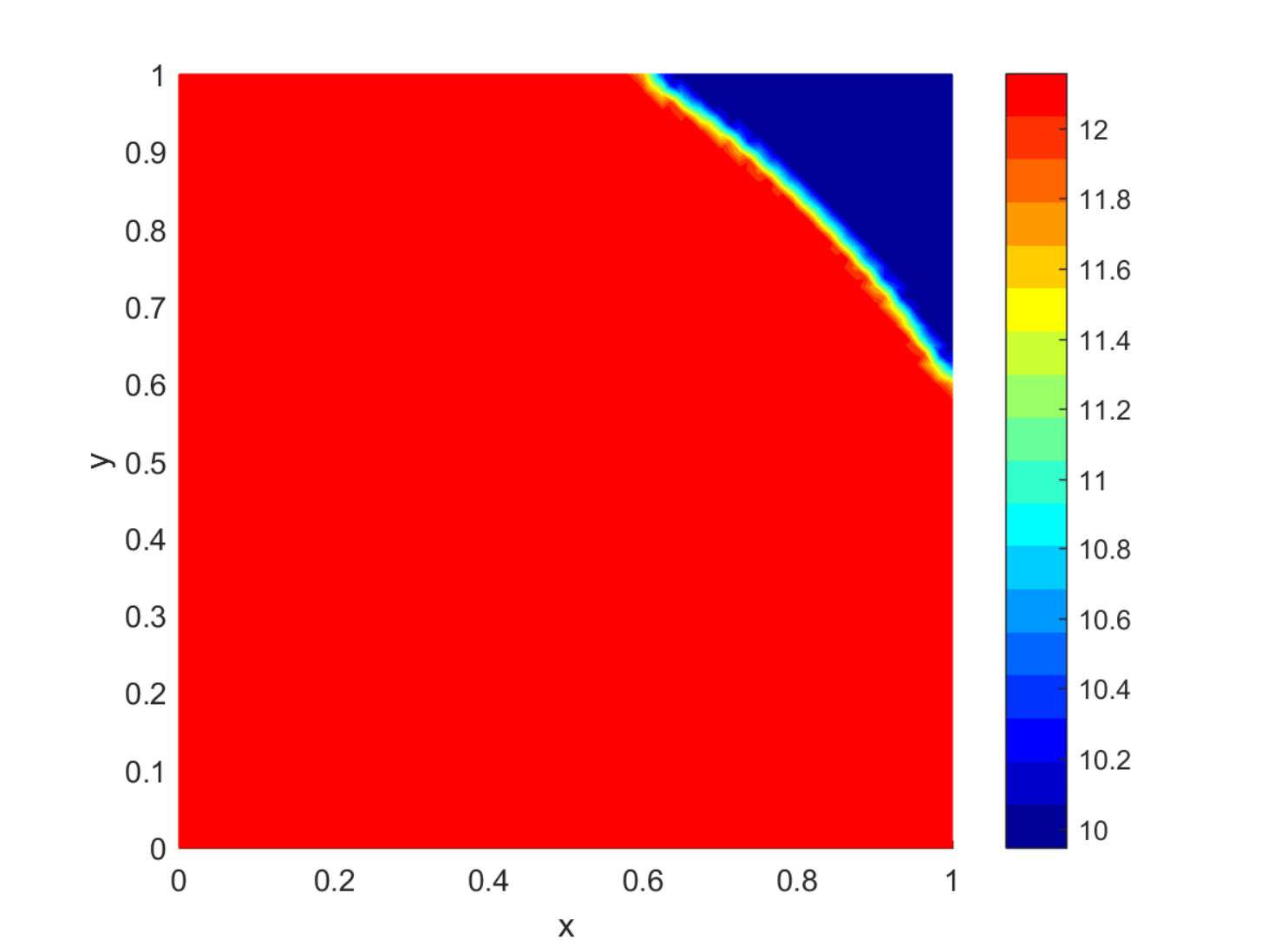}}
\subfigure[Radiative intensity on FM $N$=57600]{
\includegraphics[width=0.45\textwidth]{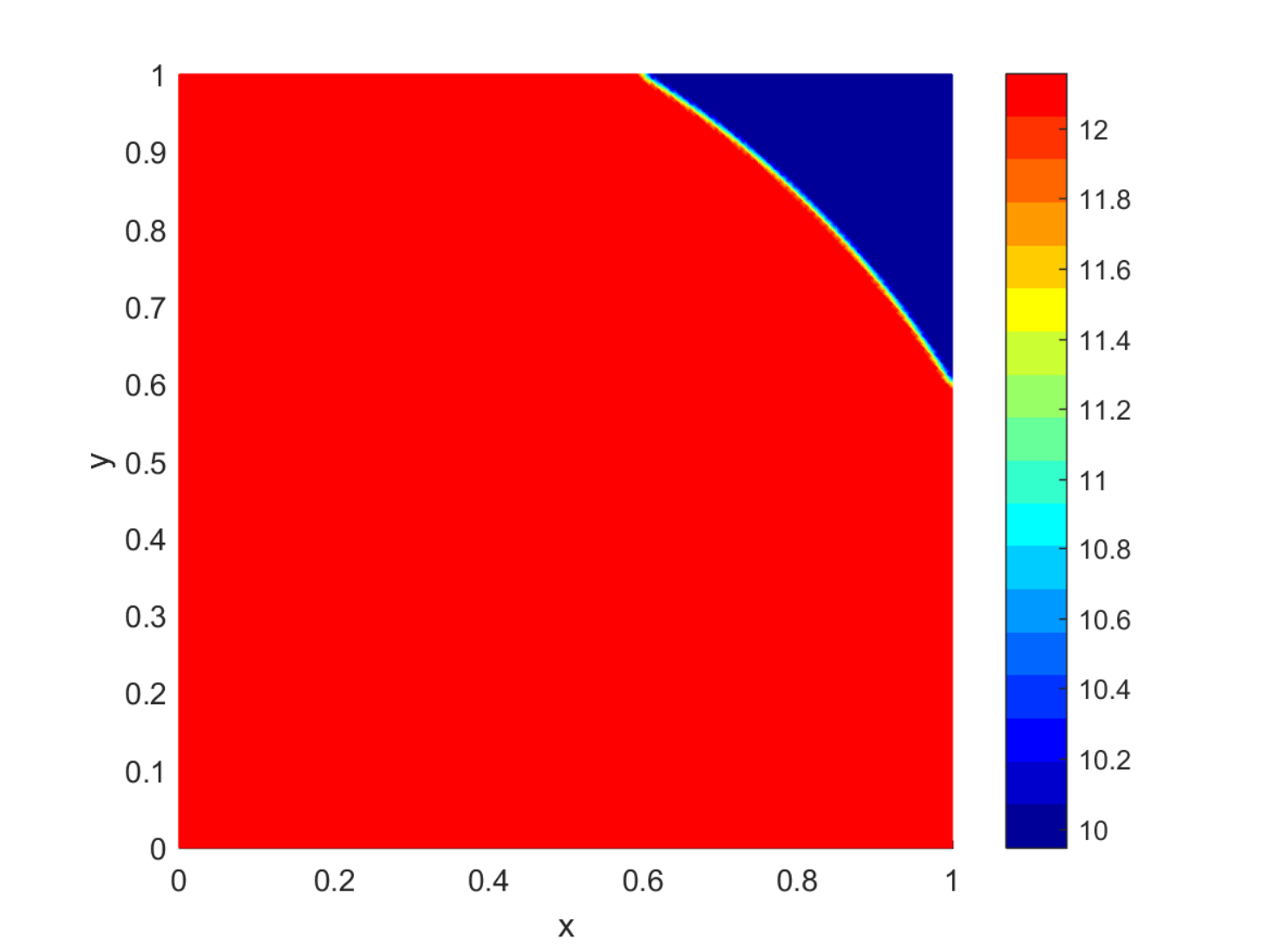}}
\caption{\small{Example \ref{Ex4-2d}. The radiative intensity contours in the direction $\Omega = (-0.2578,-0.1068)$
(and mesh) at $t = 0.1$ obtained by $P^2$-DG method with the moving mesh of $N=6400$ and fixed meshes
of $N=6400$ and $N=57600$.}}
\label{Fig:d2Ex4p2-1}
\end{figure}

\begin{figure}[H]
\centering
\subfigure[Radiative intensity on MM $N$=6400]{
\includegraphics[width=0.45\textwidth]{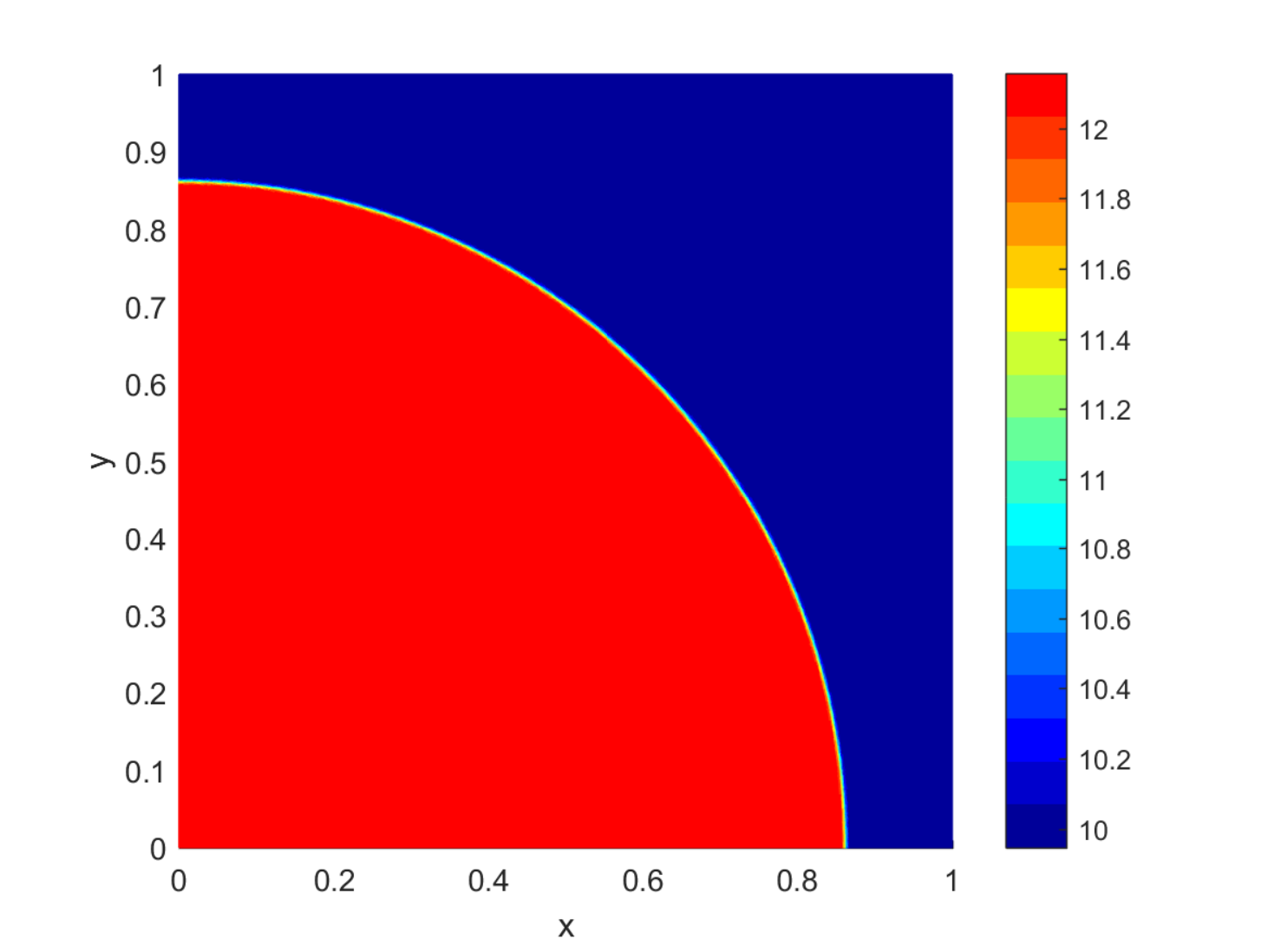}}
\subfigure[MM $N$ = 6400 at $t=0.1$]{
\includegraphics[width=0.45\textwidth]{Ex4_2d_figure_P2/Ex4_2d_P2Mesh40-eps-converted-to.pdf}}
\subfigure[ Radiative intensity on FM $N$=6400]{
\includegraphics[width=0.45\textwidth]{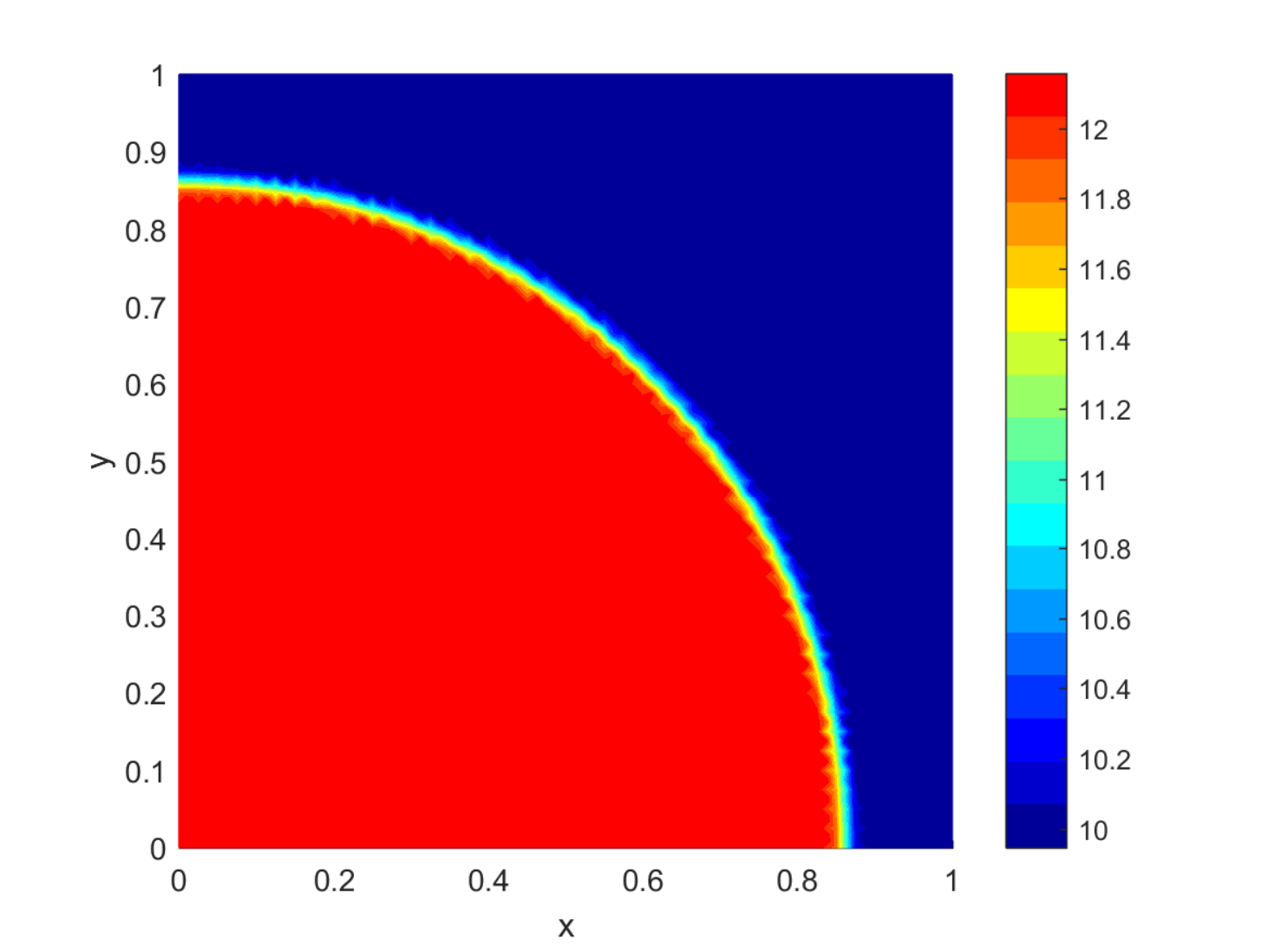}}
\subfigure[Radiative intensity on FM $N$=57600]{
\includegraphics[width=0.45\textwidth]{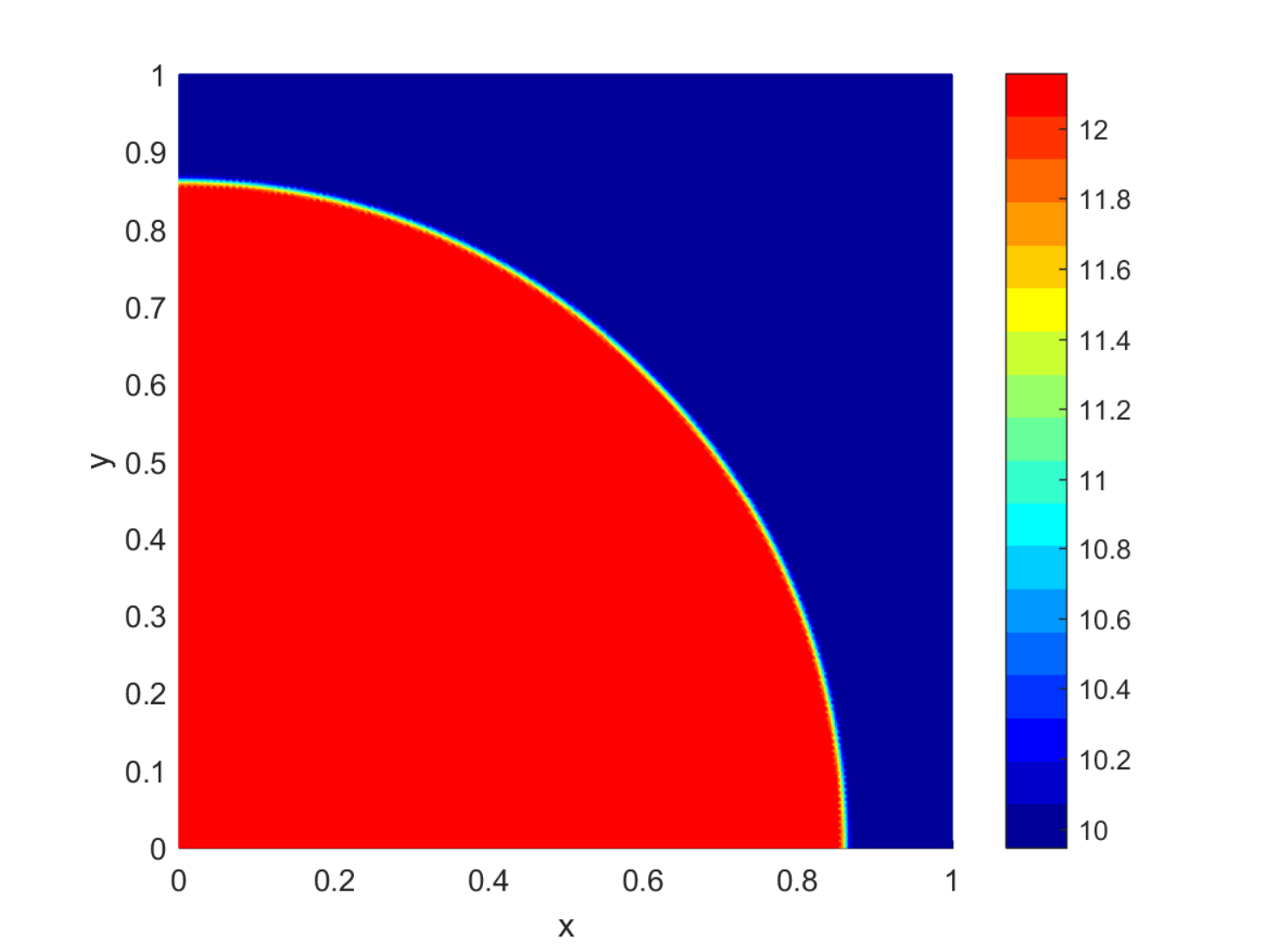}}
\caption{\small{Example \ref{Ex4-2d}. The radiative intensity contours in the direction $\Omega = (0.7860,0.3256)$
(and mesh) at $t = 0.1$ obtained by $P^2$-DG method with the moving mesh of $N=6400$ and fixed meshes of $N=6400$ and $N=57600$.}}
\label{Fig:d2Ex4p2-2}
\end{figure}

\noindent
In this example, the parameters are taken as $\sigma_t=10000$, $\sigma_s=1$, $c=3.0\times10^{8}$, and
\begin{align*}
q(x,y,\zeta,\eta,t)= & e^t\Big{(}
(\frac{1}{c}+\sigma_t)
\big{(}a-\tanh\big{(}R(x^2+y^2-\sqrt{2(1-\zeta^2-\eta^2)}\big{)}\big{)}
\\ &
-2R(\zeta x+\eta y)\big{(}1-\tanh^2\big{(}R(x^2+y^2-\sqrt{2(1-\zeta^2-\eta^2)}\big{)}\big{)}
\\
&+\frac{\sqrt{2}\sigma_s}{2R}\big{(}\ln\big{(}\cosh(R(\sqrt{2}-x^2-y^2))\big{)}
-\ln\big{(}\cosh(R(x^2+y^2))\big{)}\big{)}-\sigma_sa
\Big{)} ,
\end{align*}
where $R=200$ and $a=10$.
The initial condition is
\[
I(x,y,\zeta,\eta,0)=\Big{(}a-\tanh\big{(}R(x^2+y^2-\sqrt{2(1-\zeta^2-\eta^2)}\big{)}\Big{)}
\]
and the boundary conditions are
\begin{equation*}
\begin{split}
&I(x,0,\zeta,\eta,t)=
e^t\Big{(}a-\tanh\big{(}R(x^2-\sqrt{2(1-\zeta^2-\eta^2)}\big{)}\Big{)},
\quad\quad\quad \eta>0  ,
\\&I(x,1,\zeta,\eta,t)=
e^t\Big{(}a-\tanh\big{(}R(x^2+1-\sqrt{2(1-\zeta^2-\eta^2)}\big{)}\Big{)},
\quad~\eta<0  ,
\\&I(0,y,\zeta,\eta,t)=
e^t\Big{(}a-\tanh\big{(}R(y^2-\sqrt{2(1-\zeta^2-\eta^2)}\big{)}\Big{)},
\quad\quad\quad \zeta>0  ,
\\&I(1,y,\zeta,\eta,t)=
e^t\Big{(}a-\tanh\big{(}R(1+y^2-\sqrt{2(1-\zeta^2-\eta^2)}\big{)}\Big{)},
\quad~ \zeta<0 .
\end{split}
\end{equation*}
This problem has the exact solution
\[
I(x,y,\zeta,\eta,t)=e^t\Big{(}a-\tanh\big{(}R(x^2+y^2-\sqrt{2(1-\zeta^2-\eta^2)}\big{)}\Big{)},
\]
for which the location of the steep transition layers changes with the angular variable $\Omega=(\zeta,\eta)$.
The radiative intensity contours for the directions $\Omega=(-0.2578,-0.1068)$ and $\Omega=(0.7860,0.3256)$
are shown in Figs.~\ref{Fig:d2Ex4p2-1} and \ref{Fig:d2Ex4p2-2}, respectively, for the $P^2$-DG method with
moving and fixed meshes. It can be seen that the elements of the moving mesh are concentrated
in the regions of the sharp transition layers in the radiative intensities for all angular directions.
(The mesh shows four layers while only two layers are shown in the intensity contours
in Figs.~\ref{Fig:d2Ex4p2-1} and \ref{Fig:d2Ex4p2-2}. The figures for the intensities for other directions
are omitted to save space.)
It can be seen that the moving mesh ($N=6400$) provides a better resolution of the layers
than the fixed mesh of $N=6400$ and is comparable with the fixed mesh of $N=57600$.

\begin{figure}[H]
\centering
\subfigure[$P^1$-DG]{
\includegraphics[width=0.45\textwidth]{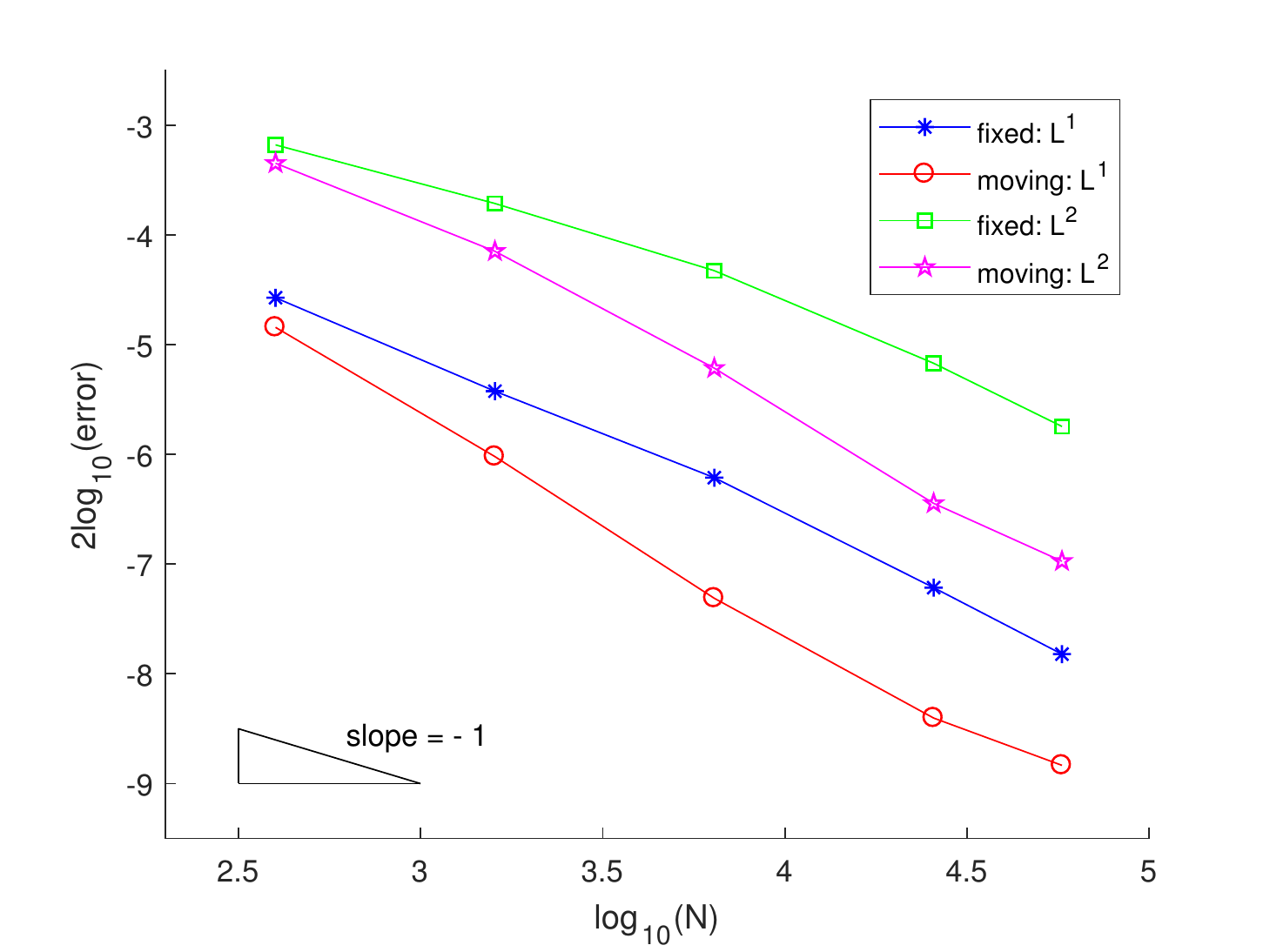}}
\subfigure[ $P^2$-DG]{
\includegraphics[width=0.45\textwidth]{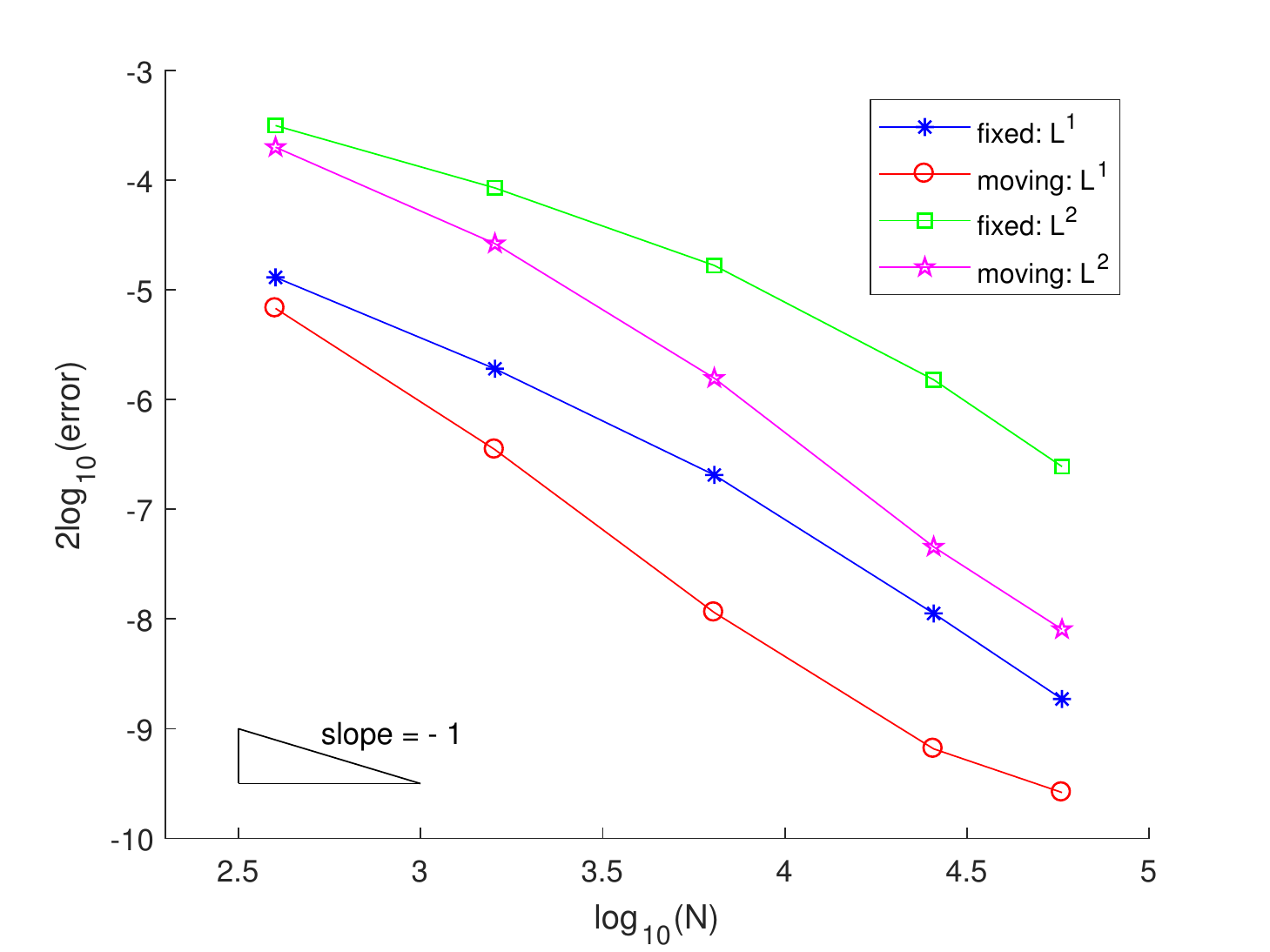}}
\caption{
\small{Example \ref{Ex4-2d}. The $L^1$ and $L^2$ norm of the error with a moving and fixed meshes.}
}\label{Fig:d2Ex4-order}
\end{figure}

The error in the $L^1$ and $L^2$ norm is plotted in Fig.~\ref{Fig:d2Ex4-order} as a function of $N$.
The convergence order is similar for both fixed and moving meshes, i.e., the order of $P^1$-DG is
about 1.7 in $L^1$ norm and 1.6 in $L^2$ norm and that of $P^2$-DG is 2.2 in both $L^1$ and $L^2$ norm.


\begin{example}\label{Ex5-2d}
(A steep transition layer example of the two-dimensional unsteady RTE for the absorbing-scattering model.)
\end{example}

\begin{figure}[H]
\centering
\subfigure[Radiative intensity on MM $N$=6400]{
\includegraphics[width=0.45\textwidth]{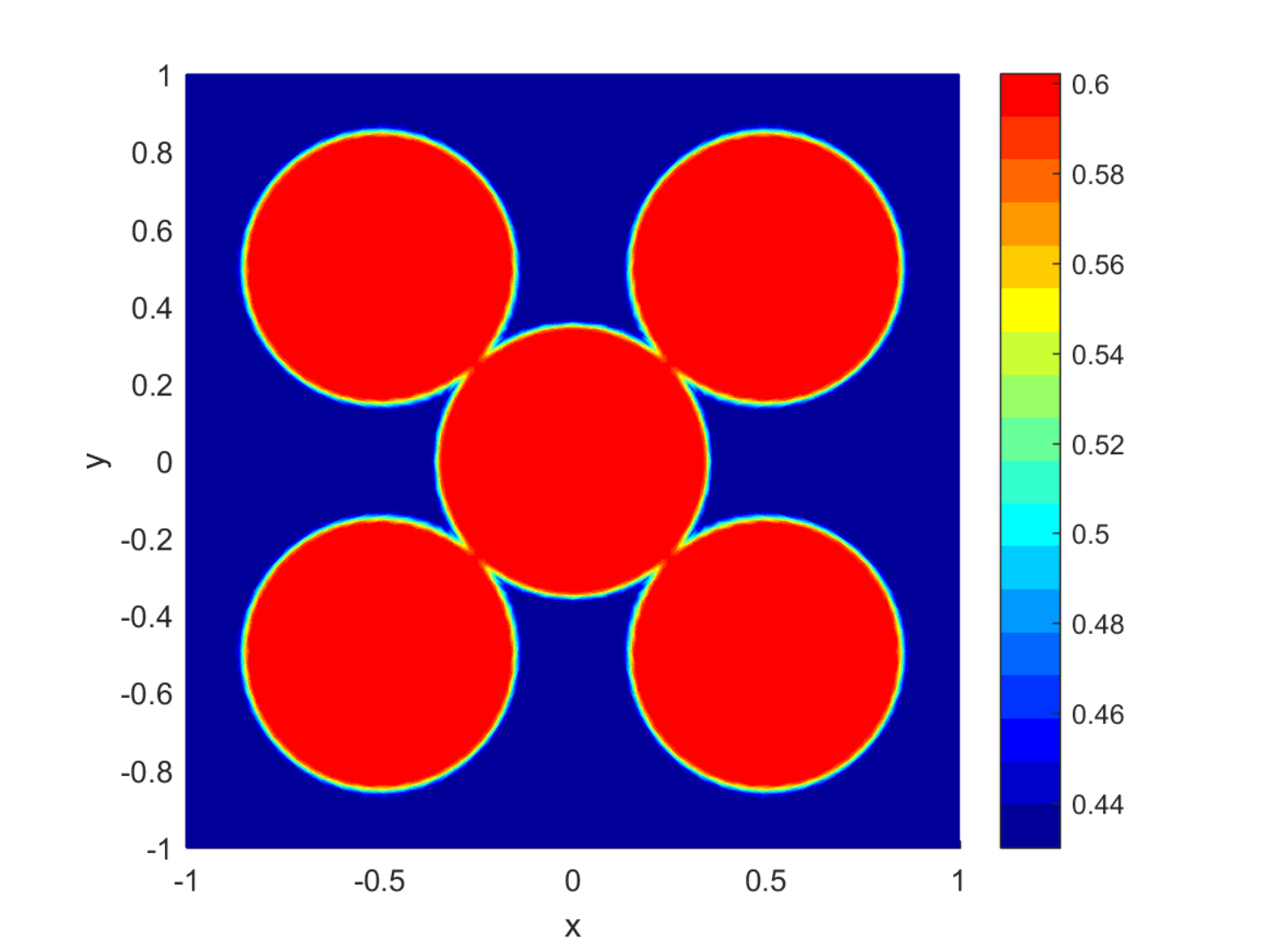}}
\subfigure[MM $N$ = 6400 at $t=0.1$]{
\includegraphics[width=0.45\textwidth]{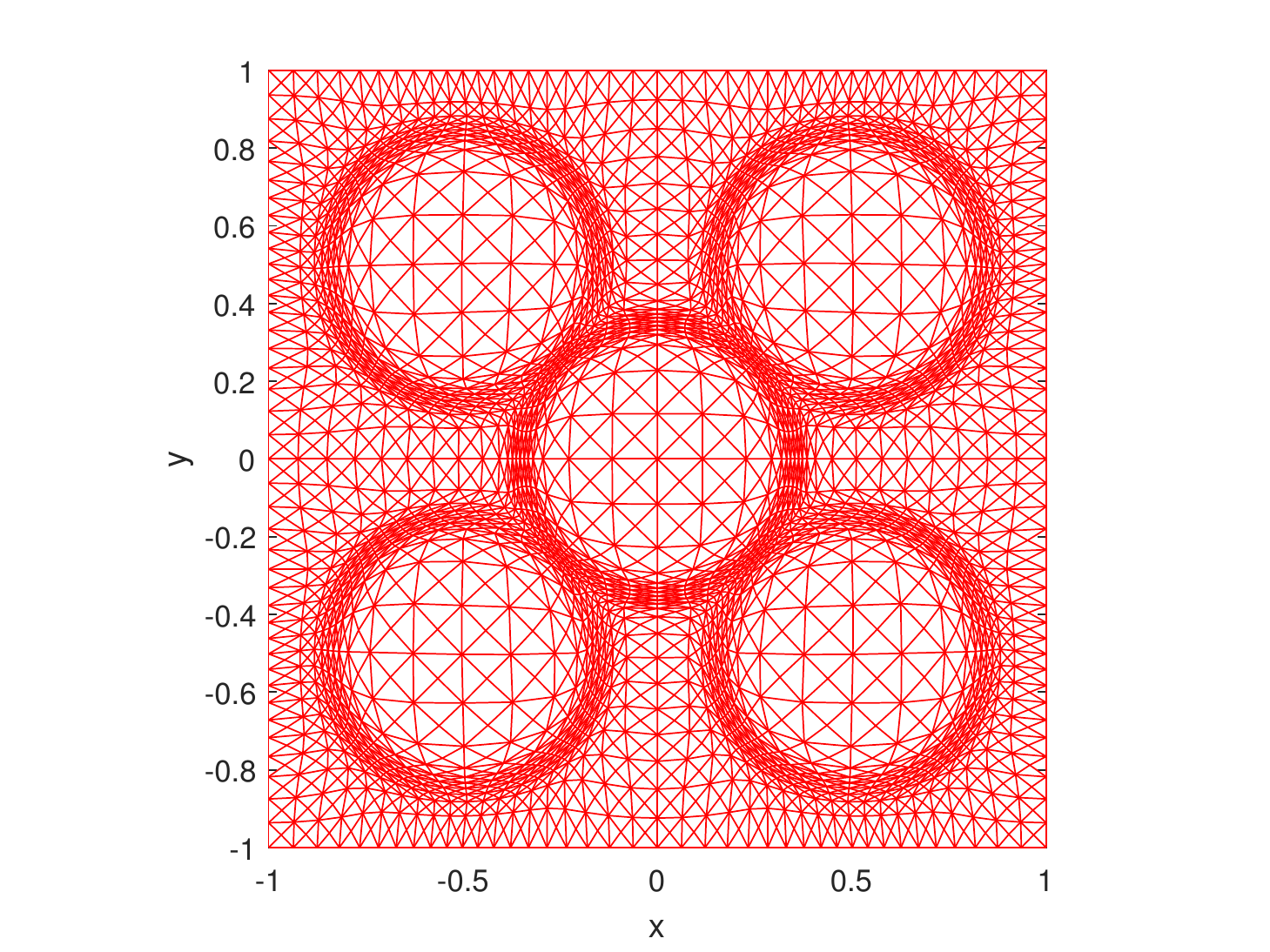}}
\subfigure[ Radiative intensity on FM $N$=6400]{
\includegraphics[width=0.45\textwidth]{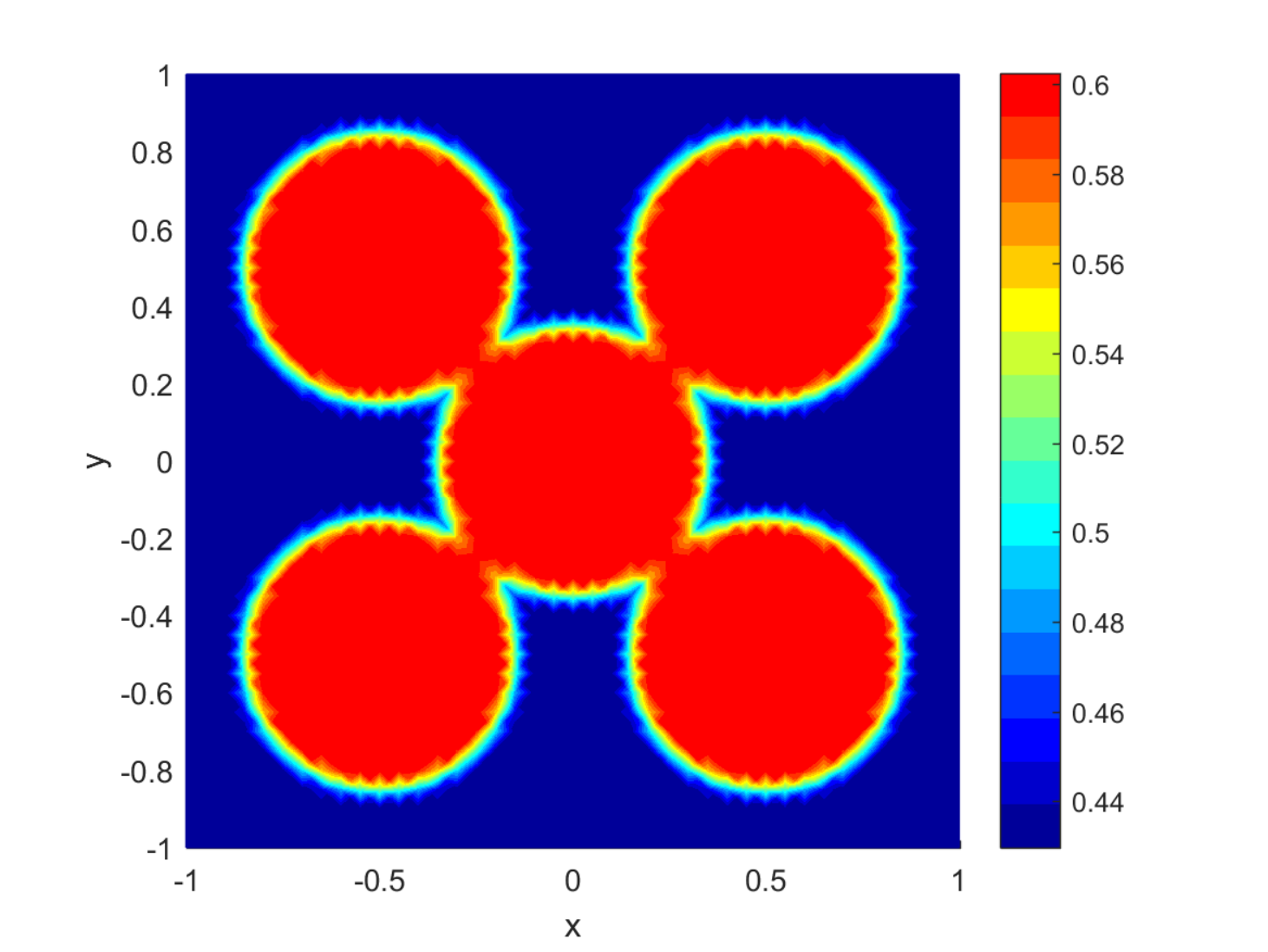}}
\subfigure[Radiative intensity on FM $N$=57600]{
\includegraphics[width=0.45\textwidth]{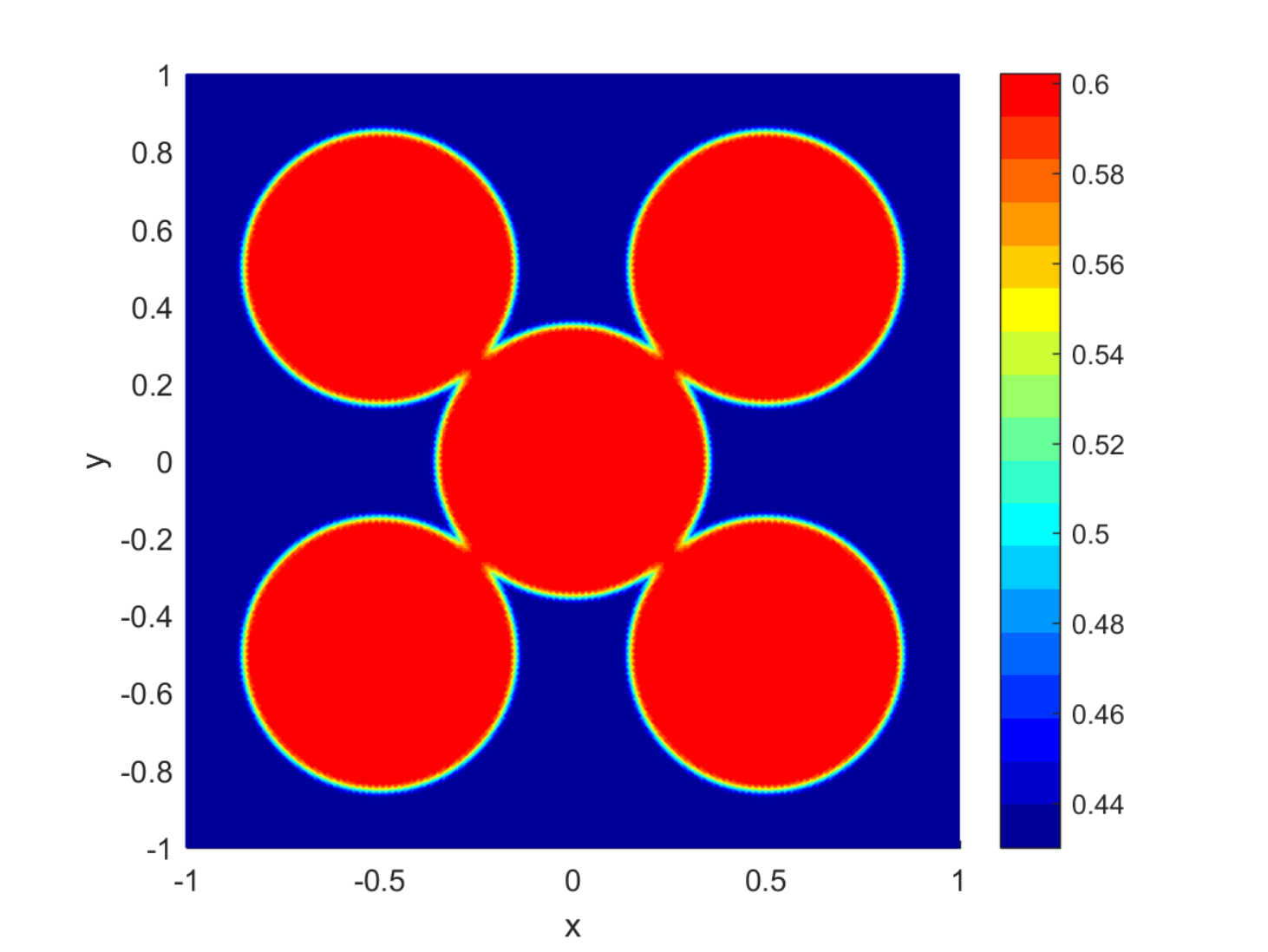}}
\caption{\small{Example \ref{Ex5-2d}. The radiative intensity contours in the direction $\Omega = (0.2578, 0.1068)$
(and mesh) at $t = 0.1$ obtained by $P^2$-DG method with the moving mesh of $N=6400$ and fixed meshes of $N=6400$ and $N=57600$.}}
\label{Fig:d2Ex5p2-1}
\end{figure}

\begin{figure}[H]
\centering
\subfigure[MM $N$=6400, FM $N$=6400]{
\includegraphics[width=0.45\textwidth]{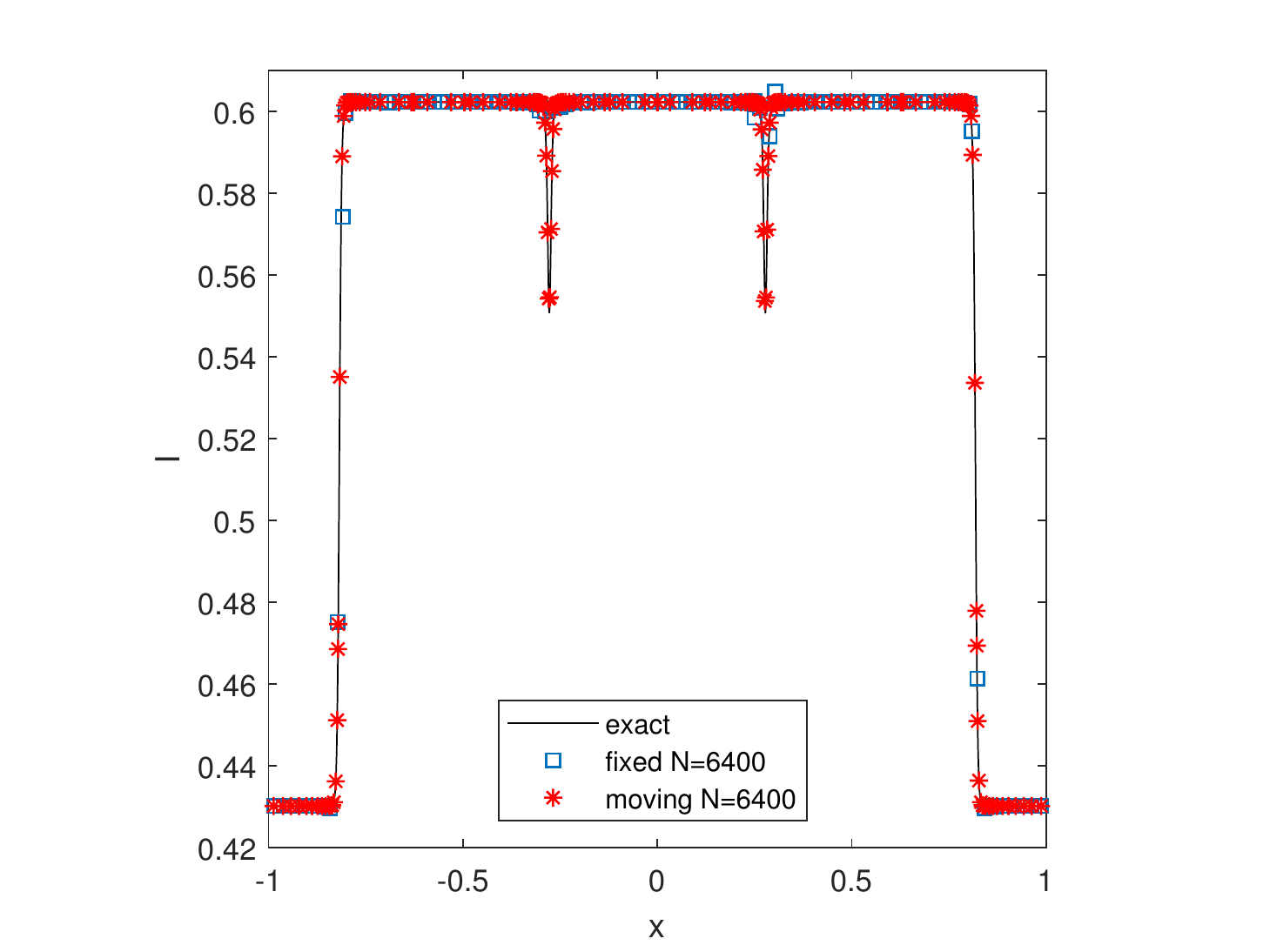}}
\subfigure[Close view of (a) near $x \in (-0.35, 0.35)$]{
\includegraphics[width=0.45\textwidth]{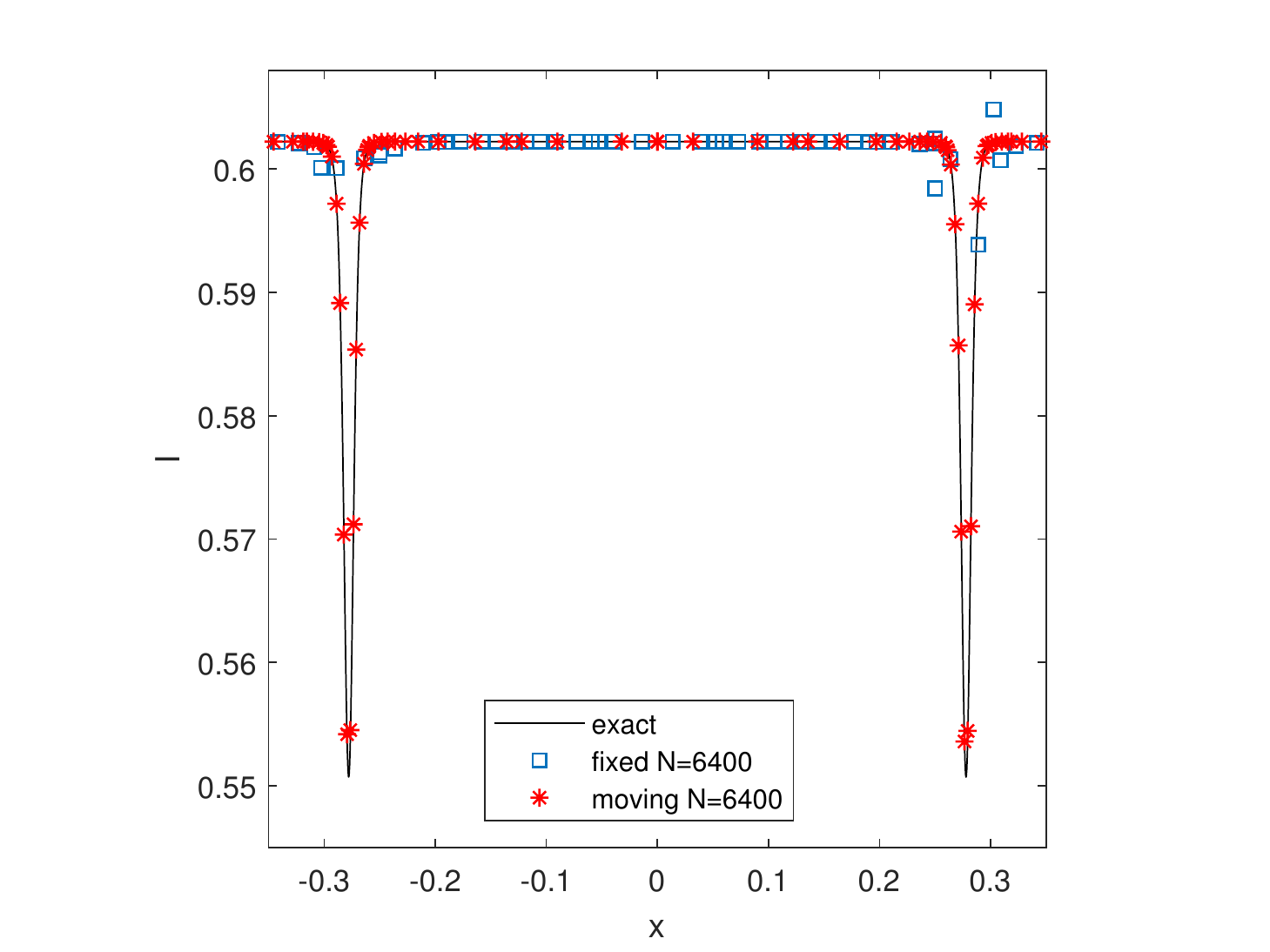}}
\subfigure[ MM $N$=6400, FM $N$=57600]{
\includegraphics[width=0.45\textwidth]{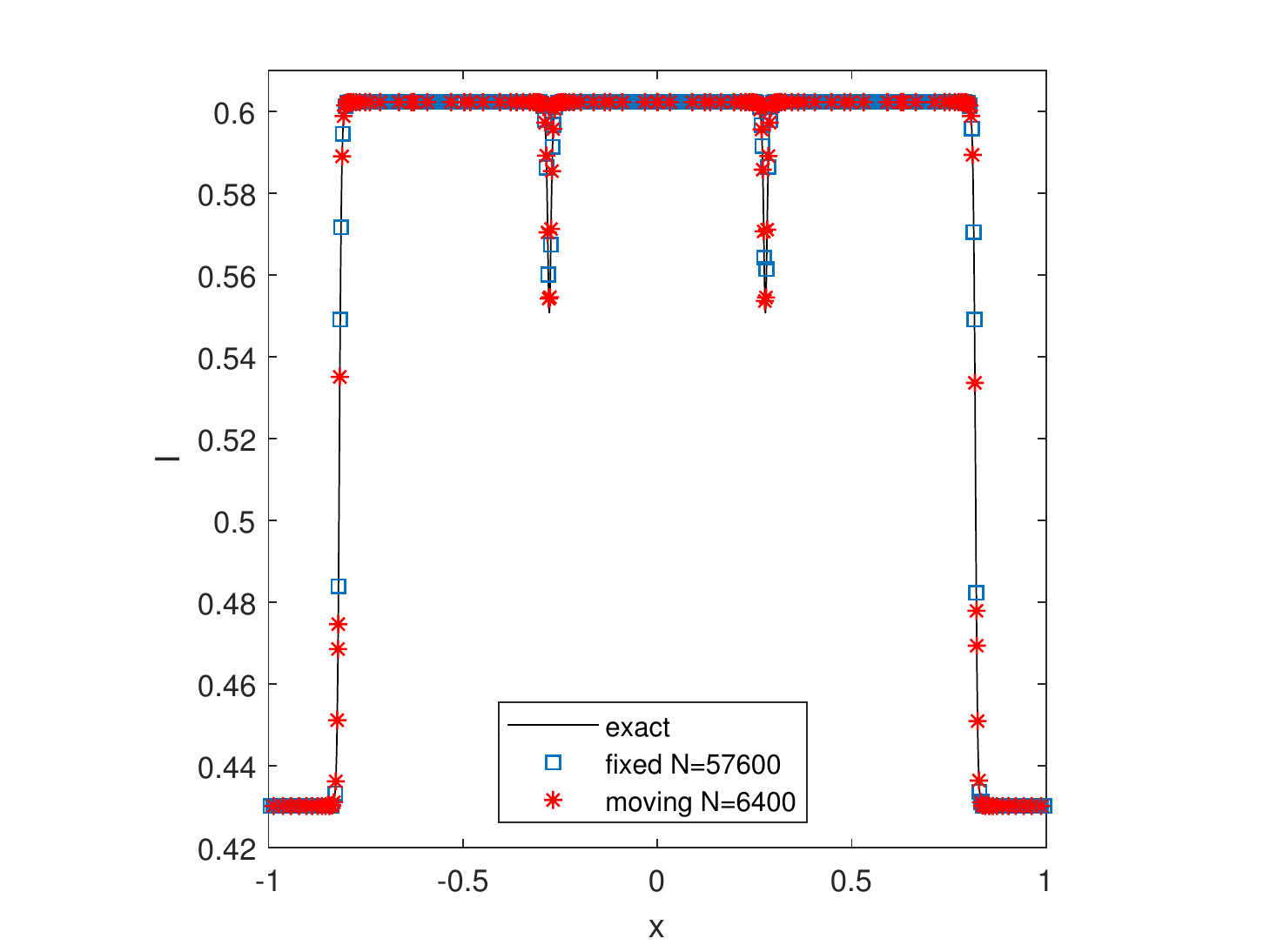}}
\subfigure[Close view of (c) near $x \in (-0.35, 0.35)$]{
\includegraphics[width=0.45\textwidth]{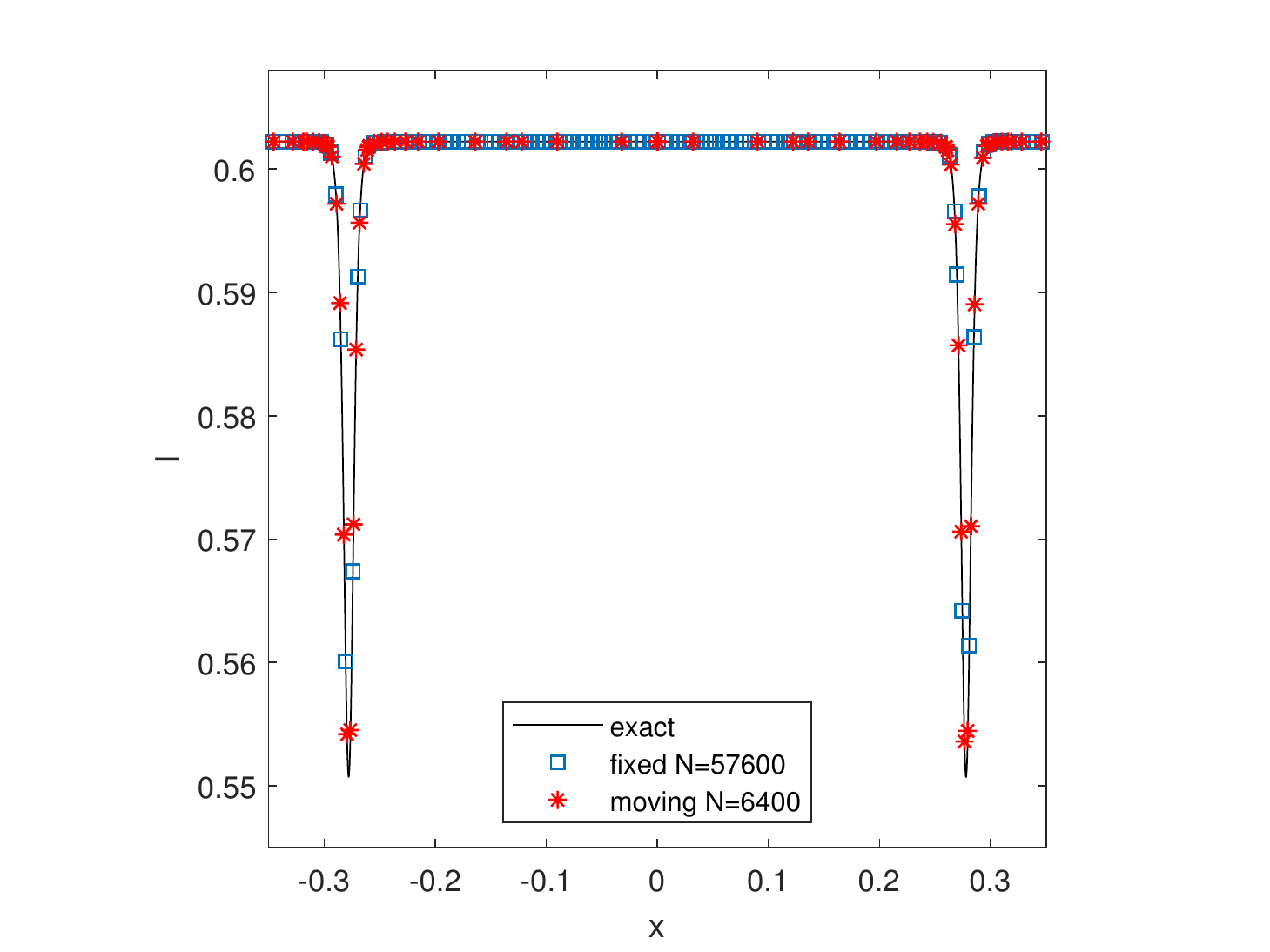}}
\caption{\small{Example \ref{Ex5-2d}. The radiative intensity in the direction $\Omega = (0.2578, 0.1068)$ cut along the line $y=0.8x$ obtained a moving mesh
of $N=6400$ is compared with those obtained with fixed meshes of $N=6400$ and $N=57600$.}}
\label{Fig:d2Ex5p2-2}
\end{figure}

\begin{figure}[H]
\centering
\subfigure[$P^1$-DG]{
\includegraphics[width=0.45\textwidth]{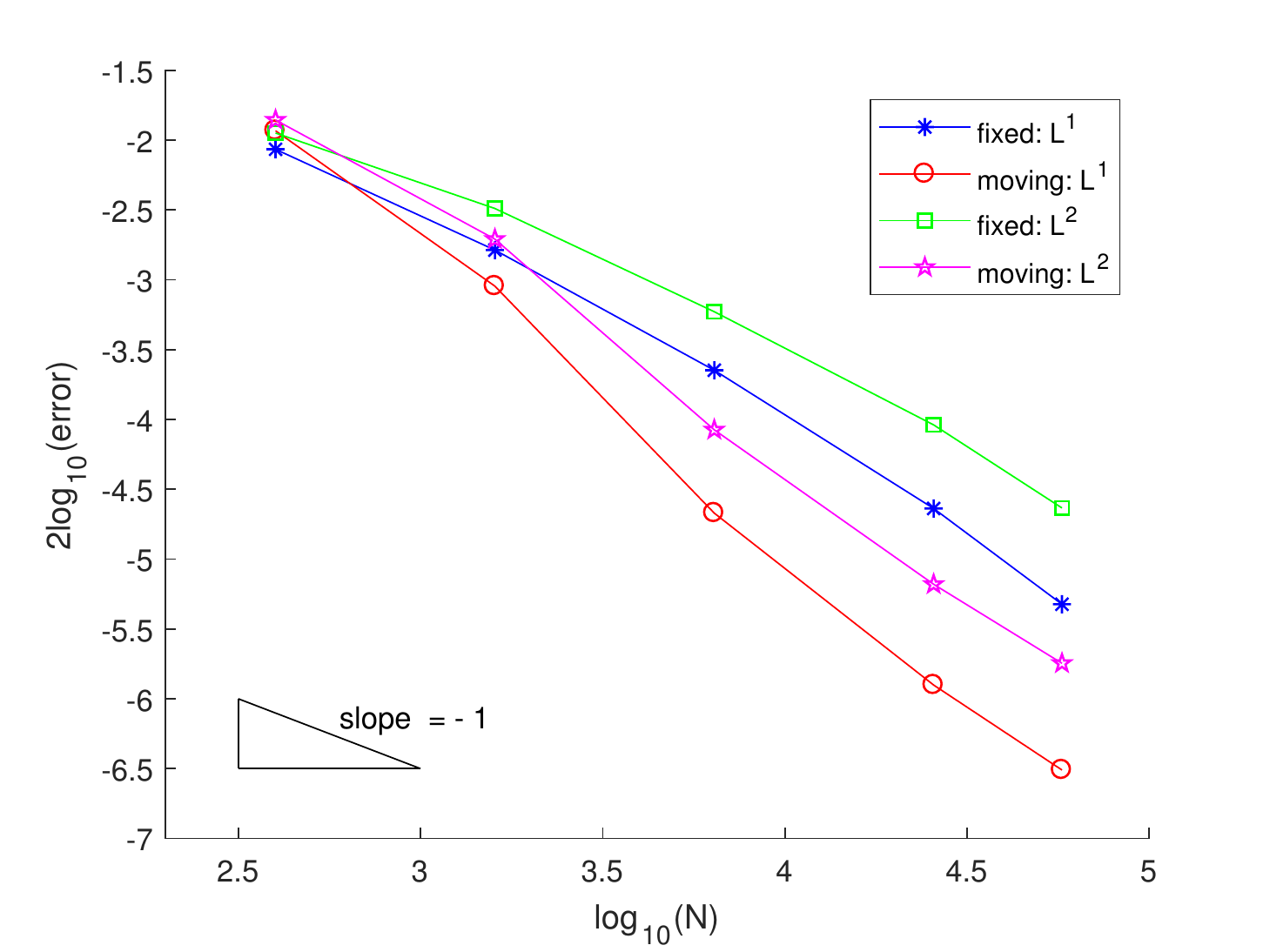}}
\subfigure[$P^2$-DG]{
\includegraphics[width=0.45\textwidth]{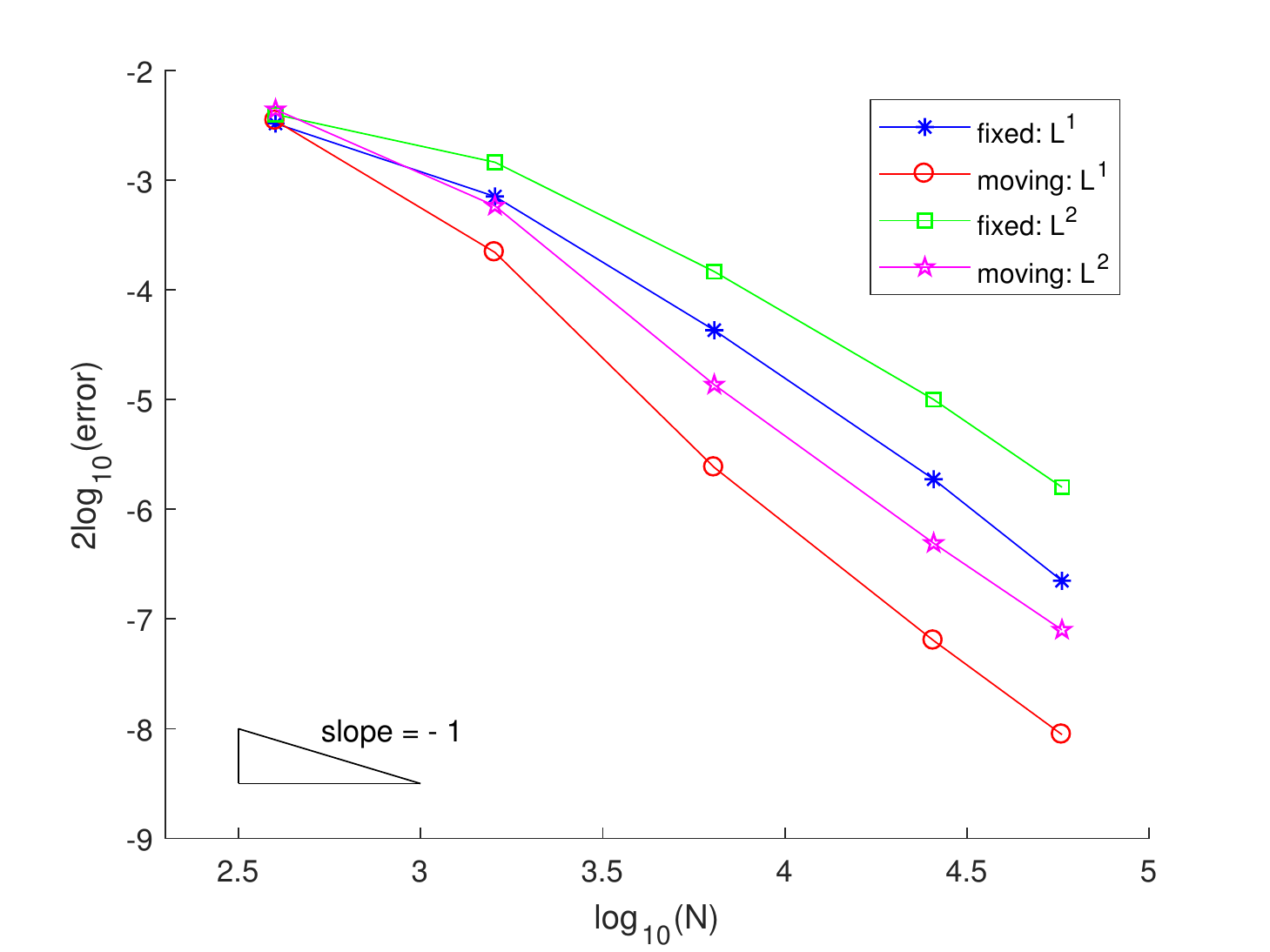}}
\caption{
\small{Example \ref{Ex5-2d}. The $L^1$ and $L^2$ norm of the error with a moving and fixed meshes.}
}\label{Fig:d2Ex5-order}
\end{figure}

\noindent
In this final example, we take $\sigma_t=33$, $\sigma_s=3$, and $c=3.0\times10^{8}$.
The computational domain is $(-1,1)\times (-1,1)$. Define
\begin{align*}
&C_{0}(x,y)=\tanh( R(x^2+y^2-\frac{1}{8}) ),
\\
&C_{1}(x,y)=\tanh( R((x-0.5)^2+(y-0.5)^2-\frac{1}{8}) ),
\\
&C_{2}(x,y)=\tanh( R((x-0.5)^2+(y+0.5)^2-\frac{1}{8}) ),
\\&C_{3}(x,y)=\tanh( R((x+0.5)^2+(y-0.5)^2-\frac{1}{8}) ),
\\& C_{4}(x,y)=\tanh( R((x+0.5)^2+(y+0.5)^2-\frac{1}{8}) ),
\end{align*}
where $R = 200$.
The source term is taken as
\begin{align*}
q(x,y,\zeta,\eta,t)=&e^t(\zeta^2+\eta^2)\Big{(}
(\frac{1}{c}+\sigma_t)\big{(}5a-\sum\limits_{i=0}^{4}C_{i}(x,y)\big{)}
\\ & - R \big{(}2(\zeta x+\eta y)(1-C_{0}(x,y)^2)
\\ &+(\zeta(2x-1)+\eta(2y-1))(1-C_{1}(x,y)^2)
\\ &
+(\zeta(2x-1)+\eta(2y+1))(1-C_{2}(x,y)^2)
\\ &+(\zeta(2x+1)+\eta(2y-1))(1-C_{3}(x,y)^2)
\\ &
+(\zeta(2x+1)+\eta(2y+1))(1-C_{4}(x,y)^2) \big{)}\Big{)}
\\ &
-\frac{2}{3}e^t\sigma_s \big{(}5a-\sum\limits_{i=0}^{4}C_{i}(x,y)\big{)},
\end{align*}
where $a = 2$.
The initial condition is $I(x,y,\zeta,\eta,0)=(\zeta^2+\eta^2)\Big{(}5a-\sum\limits_{i=0}^{4}C_{i}(x,y)\Big{)}$
and the boundary conditions are
\begin{equation*}
\begin{split}
&I(x,0,\zeta,\eta,t)=e^t(\zeta^2+\eta^2)\Big{(}5a-\sum_{i=0}^{4}C_{i}(x,0)\Big{)},\quad\eta>0  ,
\\&I(x,1,\zeta,\eta,t)=e^t(\zeta^2+\eta^2)\Big{(}5a-\sum_{i=0}^{4}C_{i}(x,1)\Big{)},\quad\eta<0,
\\&I(0,y,\zeta,\eta,t)=e^t(\zeta^2+\eta^2)\Big{(}5a-\sum_{i=0}^{4}C_{i}(0,y)\Big{)},\quad\zeta>0,
\\&I(1,y,\zeta,\eta,t)=e^t(\zeta^2+\eta^2)\Big{(}5a-\sum_{i=0}^{4}C_{i}(1,y)\Big{)},\quad\zeta<0.
\end{split}
\end{equation*}
The exact solution of the problem is
\[
I(x,y,\zeta,\eta,t)=e^t(\zeta^2+\eta^2)\Big{(}5a-\sum\limits_{i=0}^{4}C_{i}(x,y)\Big{)},
\]
which exhibits a sharp layer of five-ring shape, independent of $\zeta$ and $\eta$.
The radiative intensity contours in the directions $\Omega=(0.2578, 0.1068)$
are shown in Fig.~\ref{Fig:d2Ex5p2-1} for the $P^2$-DG method with
moving and fixed meshes. In Fig.~\ref{Fig:d2Ex5p2-2}, the radiative intensity cut along
the line $y=0.8 x $ for the direction $\Omega=(0.2578, 0.1068)$ is compared for
moving and fixed meshes. The error in $L^1$ and $L^2$ norm is plotted as a function of $N$
in Fig.~\ref{Fig:d2Ex5-order}.


\section{Conclusion}\label{sec:conclusion}

In the previous sections an adaptive moving mesh DG method has been presented for the numerical solution
of the unsteady radiative transfer equation. The RTE is an integro-differential equation modeling
the conservation of photons. It involves an integral term in the angular directions while being hyperbolic
in space. The challenges for its numerical solution include the needs to handle with its high dimensionality, the presence
of the integral term, and the development of discontinuities and sharp layers in its solution along spatial directions.
In the current work, the RTE is discretized first in the angular directions with the discrete ordinate method
and then with a DG method in space on a moving mesh. The mesh is moved adaptively
using the MMPDE strategy to provide better resolution of sharp layers or discontinuities and thus better efficiency.
The source iteration is used to avoid the coupling of the radiative intensities for all angular directions
in the integral term.

A selection of one- and two-dimensional examples have been presented to demonstrate the accuracy and efficiency
of the method. It has been shown that the method is able to automatically concentrate the mesh points
in regions of discontinuities or steep transition layers and is more efficient than its fixed mesh
counterpart. Particularly, the combination of the metric tensors for radiative intensities for different angular directions
into a single metric tensor using the matrix intersection seems to work well for the tested problems
with single or multiple sharp layers and discontinuities.
Interestingly, the results also show that the quadratic DG method has better efficiency than
the linear DG for both fixed and moving meshes.

It should be pointed out that we have not considered positivity-preserving limiters \cite{CJ, pplimiter}
nor nonoscillatory limiters such as TVB limiter \cite{bs1,bs3},
the WENO limiter \cite{WENO1,WENO2}, or HWENO limiter \cite{HWENO1,HWENO2} in the current work.
One may observe that localized spurious oscillations occur in numerical solutions containing discontinuities.
How to combine limiters with our moving mesh DG method for the unsteady RTE or more general
integro-differential equations will be an interesting research topic for the near future.
Other future work will include extending the method to the numerical solution of the RTE coupled with
with the Euler equations, the material equation or the energy equation for real situations.




\end{document}